\documentclass[onefignum,onetabnum]{siamart171218}
\usepackage{amsmath}%
\usepackage{amsfonts}%
\usepackage{amssymb}%
\usepackage{graphicx}
\usepackage{float} 
\usepackage{mathdots} 
\usepackage{todonotes}
\usepackage{showlabels}
\usepackage{mathtools}
\DeclarePairedDelimiter{\ceil}{\lceil}{\rceil}
\usepackage{algorithm}
\floatname{algorithm}{Algorithm 1}
\usepackage[noend]{algpseudocode}

\newtheorem{dfn}{Definition}[section]

\newtheorem{rmk}{Remark}[section]

\usepackage{mathtools}

\newcommand{\ep}{\varepsilon_{\mathcal{T}}}
\newcommand{\R}{\mathbb R}

\newcommand{\new}[2][black]{\emph{\textcolor{#1}{#2}}}
\newcommand{\revfirst}[2][black]{{\textcolor{#1}{#2}}}
\newcommand{\revsecond}[2][black]{\textcolor{#1}{#2}}
\newcommand{\rev}[2][black]{\textcolor{#1}{#2}}
\DeclareMathOperator*{\argmax}{arg\,max}
\DeclareMathOperator*{\argmin}{arg\,min}

\title{An efficient DP algorithm on a tree-structure\\for finite horizon optimal control problems\thanks{
MF and LS would like to thank the support obtained by the 2018 INDAM -GNCS research project {\em Metodi numerici per problemi di controllo multiscala e applicazioni}.}}
\author{Alessandro Alla\footnote{PUC-Rio, Rua Marques de Sao Vicente, 225, G\'avea - Rio de Janeiro, RJ - Brasil - 22451-900, alla@mat.puc-rio.br}, Maurizio Falcone\footnote{Sapienza Università di Roma - Piazzale Aldo Moro 5, 00185 Roma, falcone@mat.uniroma1.it}, Luca Saluzzi\footnote{Gran Sasso Science Institute - Viale F. Crispi 7, 67100 L'Aquila, luca.saluzzi@gssi.it}}
\begin{document}

\maketitle

\begin{abstract}
The classical Dynamic Programming (DP) approach to optimal control problems is based on the characterization of the value function as the unique viscosity solution of a Hamilton-Jacobi-Bellman (HJB) equation. The DP scheme for the numerical approximation of viscosity solutions of Bellman equations is typically based on a time discretization which is projected on a fixed state-space grid. The time discretization can be done by a one-step scheme for the dynamics and the projection on the grid typically uses a local interpolation. Clearly the use of a grid is a limitation with respect to possible applications in high-dimensional problems due to the {\em curse of dimensionality}.

Here, we present a new approach for finite horizon optimal control problems where the value function is computed using a DP algorithm with a tree structure algorithm (TSA) constructed by the time discrete dynamics. In this way there is no need to build a fixed space triangulation and to project on it. The tree will guarantee a perfect matching with the discrete dynamics and drop off the cost of the space interpolation allowing for the solution of very high-dimensional problems. Numerical tests will show the effectiveness of the proposed method.

\end{abstract}

\begin{keywords}
dynamic programming, Hamilton-Jacobi-Bellman equation, optimal control, tree structure
\end{keywords}
\begin{AMS}
49L20, 49J15, 49J20, 93B52
\end{AMS}

\section{Introduction}
The Dynamic Programming (DP) approach has been introduced and developed by Richard Bellman in the '50s in a series of pioneering papers (see e.g. \cite{B57}). Since then it has been applied to many problems in deterministic and stochastic optimal control although its real application has been up to now  limited to low dimensional problems. Via the Dynamic Programming Principle (DPP) one can obtain  a characterization of the value function as the unique viscosity solution of a nonlinear partial differential equation (the Hamilton-Jacobi-Bellman (HJB) equation) and then use the value function  to  get a synthesis of a feedback control law. This is the major advantage over the approach based on the Pontryagin Maximum Principle (PMP) \cite{BGP56,PBGM62} that gives necessary conditions for the characterization of the open-loop optimal control  and of  the corresponding optimal trajectory. As it is well known, the DP approach suffers from the {\em curse of dimensionality} since one has to solve a nonlinear partial differential equation (PDE) whose dimension is the same of the dynamical system. This has always been the main obstacle to apply that  theory to real industrial applications despite the large number of theoretical results established  for many classical control problems via the DP  approach (see e.g.  the monographies by Bardi and Capuzzo-Dolcetta \cite{BCD97} on deterministic control problems and by Fleming and Soner \cite{FS93} on stochastic control problems). Even in low dimension this  is a challenging problem  since the value function associated to  the control problem  (i.e. the viscosity solution of the HJB equation) is known to be only Lipschitz continuous also when the  dynamics and the running costs are regular functions.  The numerical analysis of low order numerical methods is now rather complete even for a state space in $\R^d$ and several methods have been proposed to solve the HJB equation using a number of different techniques including finite differences, semi-Lagrangian, finite  volumes and finite elements. We refer the interested reader  to the monographies by Sethian \cite{S99}, Osher and Fedkiw \cite{OF03}, Falcone and Ferretti \cite{FF13} for an extensive discussion of some of these methods and for an extended list of references on numerical methods.  
All the above mentioned methods are based on a space discretization which requires the construction of a space grid (or triangulation). For higher dimensional problems the method needs a huge amount of memory allocations and makes the problem unfeasible for a dimension $d>5$ on a standard computer. Several efforts have been made to mitigate the {\em curse of dimensionality}.  Although a detailed description of these contributions goes beyond the scopes of this paper, we want to mention \cite{FLS94} for a domain decomposition method with overlapping between the subdomains and \cite{CFLS94} for similar results without overlapping. It is important to note that in these papers the method is applied to subdomains with a rather simple geometry (see the book by Quarteroni and Valli \cite{QV99} for a general introduction to this technique) to pass down conditions to the boundaries. More recently another way to decompose the problem has been proposed in \cite{NK07} who have used a patchy decomposition based on Al'brekht method (see e.g. \cite{A61}). Later in \cite{CCFP12} the patchy idea has been implemented  taking into account an approximation of the underlying optimal dynamics  to obtain subdomains which are almost invariant with respect to the optimal dynamics, clearly in this case the geometry of the subdomains can be rather complex but the transmission conditions at the internal boundaries can be eliminated saving on the overall complexity of the algorithm. More recently other decomposition techniques for optimal control problems and games have been proposed in \cite{F16} where the parallel algorithm is based on the construction of independent sub-domains and in \cite{F18} where a parallel version of the HowardÕs algorithm is proposed and analyzed. In general, domain decomposition methods  reduce a huge problem into subproblems of manageable size and allow to mitigate the storage limitation distributing the computation over several processors. However, the approximation schemes used in every subdomain are rather standard. 

Another improvement can be obtained using efficient acceleration methods  for the computation of the value function in every subdomain. In the framework of optimal control problems an efficient acceleration technique based on the coupling between value and policy iterations has been recently proposed and studied in \cite{AFK15}. The construction of a DP algorithm for  time dependent problems has been addressed in \cite{FG99} where also a-priori error estimates have been studied.  An adaptation of similar methods for high-dimensional problems has been proposed later in \cite{CFF04}.

However, we also mention that high-dimensional problems often imply a huge amount of data and are too complex to be solved even by a direct approach based on domain decomposition (this approach is typically feasible below dimension 10). A reasonable solution to attack high-dimensional problems is to apply first model order reduction techniques (e.g. Proper Orthogonal Decomposition \cite{V13}) to have a low dimensional version of the dynamics. Thus, if the reduced system of coordinates for the dynamics has a low number of dimension (e.g. $d\approx 5$) the problem can be solved via the DP approach. 
Model reduction techniques are based on orthogonal projections where the choice of the basis functions is non trivial, e.g. it requires to compute some reference trajectories corresponding to a priori given control strategies and compute the basis via an SVD. At the end of this step, the set of controlled trajectories will be represented as a linear combination of the basis functions. Whenever we are able to compute accurate projectors we drastically reduce the dimension of the control problem, say $\ell\ll d$ but we lose the physical meaning of the projected dynamical system. \revfirst{This makes it difficult} to define a reasonable choice of the numerical domain $\Omega$ and the easiest solution is to choose $\Omega$ as a rather large box in $\R^\ell$.
%
We refer, among others, to the pioneer work on the coupling between model reduction and HJB approach \cite{KVX04} and the recent \cite{AFV17} work which provides a-priori error estimates for the aforementioned coupling method.
We also mention a sparse grid approach in \cite{GK17} where the authors apply HJB to the \revfirst{control of the wave equation} and a spectral elements approximation in \cite{KK18} which allows to solve the HJB equation up to dimension $12$.\\
Despite these efforts and the mathematical elegance of the DP approach, its impact in industrial applications is limited by this bottleneck and the
solution of many optimal control problems has been accomplished instead via 
open-loop control. More information on the topic can be found in the monographies by Hinze, Pinnau, Ulbrich Ulbrich \cite{HPUU09} and by Tr\"oltzsch \cite{T10}. 

The aim of this paper is to eliminate the space discretization and the construction of a grid to reduce the memory allocations and improve the applicability of the DP approach. This can be done for the finite horizon problem via the construction of a tree-structure that will account for the controlled dynamics. For numerical purposes, we will assume that the system has a finite number of controls at every time step $t_n$ and, to simplify the presentation, we are keeping this number $M$ constant during the evolution although the extension to a variable number $M_n$ is straightforward. Under these hypotheses starting from a point $x$ we can reach $M$ points in the state space according to the discrete time dynamics. So a single starting point will produce a tree $\mathcal{T}$ of order $O(M^{\overline{N}+1})$ points in $\overline{N}$ time steps and the number of points is exponentially increasing as expected. Note that we will not compute the value function by the DP algorithm on that tree: exploiting the Lipschitz continuity of the value function in the space variable (see Section \ref{sec:focp}) we are going to prune the tree identifying the nodes that are "very close". The {\em pruning step} of the algorithm will be governed by a pruning parameter $\ep$ and at every step many branches will be cut away so the final complexity will be drastically reduced.\\
 Working on the tree has several {\em advantages}:
 \begin{itemize}
\item[(i)] we do not need to define a priori a numerical domain $\Omega$ where we want to solve the problem,  the original tree is constructed by the controlled dynamics;
\item[(ii)] we do not need to build a space grid and to make a space interpolation on the grid nodes, therefore we do not introduce the interpolation error;
\item[(iii)] \revfirst{the pruned tree allows to deal with high-dimensional problems.}
\end{itemize}
In conclusion, with respect to the standard space discretization we can drop the interpolation step that is rather expensive in high-dimension and we do not need the classical assumptions at the boundary of $\Omega$ which classically requires to have an invariant dynamics or to impose boundary conditions (Dirichlet, Neumann or state constraint). 

Via the tree structure algorithm (TSA) we eliminate these difficulties at least for the finite horizon problem and we can directly solve the discrete time HJB equation for $d=1000$ \revfirst{without any particular assumption on the structure of the problem as in model reduction context. This will be shown in Section \ref{sec:nt}.}

The paper is organized as follows: in Section~\ref{sec:focp} we recall some basic facts about the time approximation of the finite horizon problem via the DP approach, we introduce our notation and prove that the discrete time value function is Lipschitz continuous in space. Section~\ref{sec:tree} is devoted to present the construction of the tree-structure related to the controlled dynamics. In Section~\ref{sec:hints}, we present some hints on the actual implementation of the method, in particular the pruning technique used to cut off the branches of the tree in order to reduce the global complexity of the algorithm. Some numerical tests are presented and analyzed in Section~\ref{sec:nt}. We give our conclusions and perspectives in Section \ref{sec:con}.

\section{Finite horizon optimal control problems via dynamic programming principle}\label{sec:focp}
In this section we will summarize the basic results that will constitute the building blocks for our new algorithm. The essential features will be briefly sketched, and more details can be found in \cite{BCD97,FF13} and the references therein. Let us present the method for the classical {\it finite horizon problem}. Let the system be driven by

\begin{equation}\label{eq}
\left\{ \begin{array}{l}
\dot{y}(s)=f(y(s),u(s),s), \;\; s\in(t,T],\\
y(t)=x\in\R^d.
\end{array} \right.
\end{equation}
We will denote by $y:[t,T]\rightarrow\R^d$ the solution, by $u$ the  control $u:[t,T]\rightarrow\R^m$, by $f:\R^d\times\R^m\times[t,T]\rightarrow\R^d$ the dynamics and by
\[\mathcal{U}=\{u:[t,T]\rightarrow U, \mbox{measurable} \}
\]
the set of admissible controls where $U\subset \R^m$ is a compact set. 
We assume that there exists a unique solution for \eqref{eq} for each $u\in\mathcal{U}$. 

The cost functional for the finite horizon optimal control problem will be given by 
\begin{equation}\label{cost}
 J_{x,t}(u):=\int_t^T L(y(s,u),u(s),s)e^{-\lambda (s-t)}\, ds+g(y(T))e^{-\lambda (T-t)},
\end{equation}
where $L:\R^d\times\R^m\times [t,T]\rightarrow\R$ is the running cost, \revfirst{$g:\R^d\rightarrow\R$ is the final cost} and $\lambda\geq0$ is the discount factor.\\
The goal is to find a state-feedback control law $u(t)=\Phi(y(t),t),$ in terms of the state \revfirst{variable} $y(t),$ where $\Phi$ is the feedback map. To derive optimality conditions we use the well-known DPP due to Bellman. We first define the value function for an initial condition $(x,t)\in\R^d\times [t,T]$:
\begin{equation}
v(x,t):=\inf\limits_{u\in\mathcal{U}} J_{x,t}(u)
\end{equation}
which satisfies the DPP, i.e. for every $\tau\in [t,T]$ :
\begin{equation}\label{dpp}
v(x,t)=\inf_{u\in\mathcal{U}}\left\{\int_t^\tau L(y(s),u(s),s) e^{-\lambda (s-t)}ds+ v(y(\tau),\tau) e^{-\lambda (\tau-t)}\right\}\;.
\end{equation}
Due to (\ref{dpp}) we can derive the HJB for every $x\in\R^d$, $s\in [t,T)$: 
\begin{equation}\label{HJB}
\left\{
\begin{array}{ll} 
&-\dfrac{\partial v}{\partial s}(x,s) +\lambda v(x,s)+ \max\limits_{u\in U }\left\{-L(x, u,s)- \nabla v(x,s) \cdot f(x,u,s)\right\} = 0 \;, \\
&v(x,T) = g(x)\;.
\end{array}
\right.
\end{equation}
Suppose that the value function is known, by e.g. \eqref{HJB}, then it is possible to  compute the optimal feedback control as:
\begin{equation}\label{feedback}
u^*(t):=  \argmax_{u\in U }\left\{-L(x,u,t)- \nabla v(x,t) \cdot f(x,u,t)\right\}. 
\end{equation}
Equaton \eqref{HJB} is a nonlinear PDE of the first order which is hard to solve analitically although a general theory of weak solutions is available in e.g. \cite{BCD97}. 
Rather, we can solve equation \eqref{HJB} numerically by means of finite difference or semi-Lagrangian methods. 
In the current work we recall the semi-Lagrangian method. One usually starts the numerical method by discretizing in time the \revfirst{underlying control problem} with a time step $\Delta t: = [(T-t)/\overline N]$ where $\overline{N}$ is the number of temporal time steps \revfirst{ and then projects} the semi-discrete scheme on a grid obtaining the fully discrete scheme:
\begin{equation}
\left\{\begin{array}{ll}\label{SL}
V_i^{n}=\min\limits_{u\in U}[\Delta t\,L(x_i, u, t_n)+e^{-\lambda \Delta t}I[V^{n+1}](x_i+\Delta t f(x_i, u, t_n))], \\
\qquad\qquad\qquad\qquad\qquad\qquad\qquad\qquad\qquad\qquad\qquad\qquad n= \overline{N}-1,\dots, 0,\\
V_i^{\overline{N}}=g(x_i) \qquad\qquad\qquad\qquad\qquad\qquad\qquad\qquad\qquad\quad x_i\in \Omega,
\end{array}\right.
\end{equation}
where $t_n=t+n \Delta t,\, t_{\overline N} = T$, $\Omega$ is the numerical domain and $x_i$ is an element of its discretization, $V^n_i:=V(x_i,t_n)$ and  $I[\cdot]$ is an interpolation operator which is necessary to compute the value of $V^n$ at the point $x_i+\Delta t\, f(x_i,u,t_n)$ 
(in general, this point will not be a node of the grid). The interested reader will find in \cite{FG99} a detailed presentation of the scheme and a priori error estimates for its numerical approximation. 
We note that it is possible to show that the value function $v(x,t)$ is Lipschitz continuous on compact sets provided that  $f, L$ and $g$ are Lipschitz continuous with constant $L_f,L_L,L_g>0$ respectively. It is possible to extend the result for the numerical value function $V(x,t)$ as explained in the following proposition. The proof follows closely from the continuous version in \cite[Prop. 3.1]{BCD97}.
\begin{proposition}\label{prp:lip}
Let us suppose the functions $\revsecond{f(\cdot,u,t)}, \revsecond{L(\cdot,u,t)}$ and $\revsecond{g(\cdot)}$ are Lipschitz continuous uniformly with respect to the other variables. Then, the\,numerical value function $V^n(x)$ is Lipschitz in $x$


\begin{equation}\label{num:lip}
|V^n(x)-V^n(y)| \le
\left\{\begin{array}{ll}
 |x-y|\left(\frac{L_L}{L_f-\lambda} (e^{(T-t_n)(L_f-\lambda)}-1)+L_g  e^{(T-t_n) (L_f-\lambda)} \right), \\
\qquad\qquad\qquad\qquad\qquad\qquad\qquad\qquad\qquad\qquad \mbox{ for }L_f>\lambda,\\
 |x-y|\left(L_L(T-t_n)+L_g  e^{(T-t_n) (L_f-\lambda)} \right), \\
\qquad\qquad\qquad\qquad\qquad\qquad\qquad\qquad\qquad\qquad \mbox{ for }L_f\le\lambda,
\end{array}\right.
\end{equation} 
$\rev{\forall \, x, y \in \mathbb{R}^d}$ and $n=0,\ldots, \overline{N}$.
\end{proposition}

\begin{proof*}
In the case $n=\overline{N}$, we have that $V^{\overline{N}}(x)=g(x)$, then the estimate follows directly from the hypothesis on $g$.\\
In the case $n <\overline{N}$, \rev{we fix $\overline{x},\overline{y} \in \mathbb{R}^d$} and  consider the following quantity $V^n(\rev{\overline{x}})-V^n(\rev{\overline{y}})$:
\begin{align}\label{ineq}
V^n(\rev{\overline{x}})-V^n(\rev{\overline{y}}) \le \; &e^{-\lambda \Delta t}V^{n+1}(\rev{\overline{x}}+ \Delta t f(\rev{\overline{x}},u^n_*,t_n))+ \Delta t L(\rev{\overline{x}},u^n_*,t_n)\nonumber\\
 &\quad -e^{-\lambda \Delta t}V^{n+1}(\rev{\overline{y}}+ \Delta t f(\rev{\overline{y}},u^n_*,t_n))- \Delta t \, L(\rev{\overline{y}},u^n_*,t_n)\nonumber\\
&\le e^{-\lambda \Delta t}( V^{n+1}(\rev{\overline{x}}+ \Delta t f(\rev{\overline{x}},u^n_*,t_n))-V^{n+1}(\rev{\overline{y}}+ \Delta t f(\rev{\overline{y}},u^n_*,t_n)))&\\
&+ \Delta t\, L_L |\rev{\overline{x}}-\rev{\overline{y}}|, \nonumber
\end{align}

provided that
\rev{
$$
u^n_*= \argmin_{u\in U }\left\{e^{-\lambda \Delta t}V^{n+1}\left(\overline{y}+ \Delta t f(\overline{y},u,t_n)\right)+ \Delta t L(\overline{y},u,t_n) \right\}.
$$
}
To achieve the desired estimate \eqref{num:lip}, we need to iterate \eqref{ineq} \rev{starting from $\overline{x}$ and $\overline{y}$ at time $t_n$}. Let us first define 
the whole tree paths \rev{$\{x^m \}_m$ and $\{y^m \}_m$ as
$$x^{m}:= x^n+ \Delta t \sum_{j=n}^{m-1} f(x^j, u_*^j,t_j), \qquad y^m:= y^n + \Delta t \sum_{j=n}^{m-1} f(y^j, u_*^j,t_j),$$
where
$$
u^{j}_*= \argmin_{u\in U }\left\{e^{-\lambda \Delta t}V^{j+1}\left(y^j+ \Delta t  f(y^j,u,t_j)\right)+ \Delta t L(y^j,u,t_{j}) \right\}, \quad j=n,\ldots,m-1. 
$$
}
\rev{
By the discrete Gr\"onwall's lemma, it is easy to prove the following estimate for Euler schemes starting from $x^n=\overline{x}$ and $y^n=\overline{y}$
\begin{equation}
|x^{n+k}-y^{n+k}| \le |x^n-y^n| e^{k \Delta t L_f}=|\overline{x}-\overline{y}| e^{k \Delta t L_f} , \quad k=0,\ldots,\overline{N}-n. 
\label{gron}
\end{equation}
}
\rev{
Then, iterating \eqref{ineq} we obtain}
\rev{
\begin{align}
V^n(\overline{x})-V^n(\overline{y}) &\le \Delta t\, L_L  \sum_{k=0}^{\overline{N}-n-1} e^{-\lambda  k\Delta t}|x^{n+k}-y^{n+k}| + e^{-\lambda (T-t_n)}|g(x^{\overline{N}})-g(y^{\overline{N}}) |\nonumber\\
&\le \Delta t\, L_L  \sum_{k=0}^{\overline{N}-n-1} e^{-\lambda  k\Delta t}|x^{n+k}-y^{n+k}| +L_g e^{-\lambda (T-t_n)} |x^{\overline{N}}-y^{\overline{N}} |\nonumber\\
& \le|\overline{x}-\overline{y}|\left(\Delta t L_L\sum_{k=0}^{\overline{N}-n-1} e^{  k\Delta t (L_f-\lambda)} +L_g  e^{(T-t_n) (L_f-\lambda)} \right),
\label{est}
\end{align}
} \rev{
where we used \eqref{gron} and the Lipschitz continuity of $g$.}
%

If $L_f > \lambda$, then by \eqref{est} and \rev{the equality $(\overline{N}-n) \Delta t= T-t_n$}, we get 
\rev{
\begin{align}
\begin{aligned}
V^n(\overline{x})-V^n(\overline{y}) &\le |\overline{x}-\overline{y}|\left(\Delta t L_L \frac{e^{(T-t_n)(L_f-\lambda)}-1}{e^{\Delta t(L_f-\lambda)}-1} +L_g  e^{(T-t_n) (L_f-\lambda)} \right)\\
&\le |\overline{x}-\overline{y}|\left(\frac{L_L}{L_f-\lambda} (e^{(T-t_n)(L_f-\lambda)}-1) +L_g  e^{(T-t_n) (L_f-\lambda)} \right),
\end{aligned}
\end{align}
}
whereas if $L_f \le \lambda$, noticing that $e^{  k\Delta t (L_f-\lambda)} \le 1$, we directly obtain
\begin{equation}
V^n(\rev{\overline{x}})-V^n(\rev{\overline{y}}) \le |\rev{\overline{x}}-\rev{\overline{y}}|\left(L_L(T-t_n)+L_g  e^{(T-t_n) (L_f-\lambda)} \right).
\end{equation}
Analogously, it is possible to obtain the same estimate for $V^n(\rev{\overline{y}})-V^n(\rev{\overline{x}})$ which leads to the desired result.

$\qed$
\end{proof*}
In the next section we will take advantage of the estimate \eqref{ineq} to guarantee the feasibility of our proposed method.
The numerical approximation of the feedback control \eqref{feedback} follows directly from the SL-scheme \eqref{SL} and reads
\begin{align*}
u_*^n(x)=\argmin\limits_{u\in U}[\Delta t\,L(x, u, t_n)+e^{-\lambda \Delta t}I[V^{n+1}](x+ &\Delta t  f(x, u, t_n))].
\end{align*}

\section{HJB on a tree structure}\label{sec:tree}
The DP approach for the numerical approximation of viscosity solutions of the HJB equation is typically based on a time discretization which is projected on a fixed state-space grid of the numerical domain. The choice of the numerical domain is already one bottleneck of the method. In fact, although the theory is valid in the whole space $\R^d$ for computational reasons we need to restrict to a compact set in $\R^d$ which should be large enough to include all the possible trajectories. That also yields the selection of some boundary conditions which are not trivial.

In this section we will provide a novel algorithm which does not require a state-space grid and therefore avoids (i) the choice of the numerical domain, (ii) the computation of polynomial interpolation, (iii) the selection of boundary conditions and finally (iv) we can solve the problem for larger dimension, such as $d\gg5$ (in Section \ref{sec:nt} \revfirst{we provide an example} in dimension $1000$). Note that dimension $5$ was the maximum dimension for SL-schemes based on a grid on a standard computer (see e.g. \cite{AFV17}). 


 
\begin{figure}[htbp]
\centering
\includegraphics[scale=0.6]{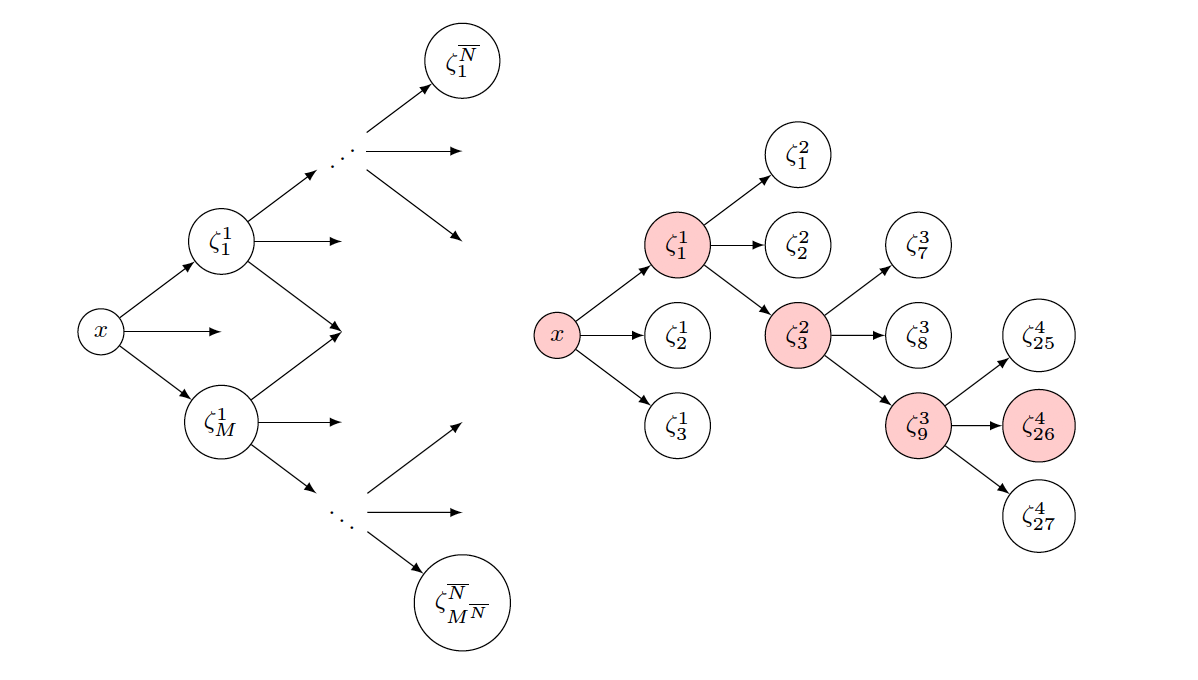}
\caption{Example of the tree $\mathcal{T}$ (left), path to reach $\zeta_{26}^4$ starting from the initial condition $x$ with $U=\{u_1,u_2,u_3\}$ (right).}
\label{fig:tree}
\end{figure}

\paragraph{\bf Construction of the tree data structure}
We build the nodes tree $\mathcal{T}$ starting from a given initial condition $x$ and following directly the dynamics in \eqref{eq} discretized by e.g. Euler method. Since we only discretize in time, we set a temporal step $\Delta t$ which divides the interval $[t,T]$ into $\overline{N}$ subintervals. We note that $\mathcal{T}:=\cup_{j=0}^{\overline{N}} \mathcal{T}^j $, where each $\mathcal{T}^j$ contains the nodes of the tree correspondent to time $t_j$. The first level $\mathcal{T}^0 = \{x\}$ is simply given by the initial condition $x$. To compute the second level, and the other levels we suppose to discretize the control domain $U$ with step-size $\Delta u$. We denote that the control set $U$ is a subset in $\mathbb{R}^m$, in particular we will consider $U$ as a hypercube, discretized in all directions with constant step-size $\Delta u$, obtaining $U^{\Delta u}=\{u_1,...,u_M \}$. To ease the notation in the sequel we continue to denote by $U$ the discrete set of controls.  Then, starting from the initial condition $x$, we consider all the nodes obtained following the dynamics \eqref{eq} discretized using e.g. an explicit Euler scheme with different discrete controls $u_j \in \mathbb{U} $
$$\zeta_j^1 = x+ \Delta t \, f(x,u_j,t_0),\qquad j=1,\ldots,M.$$ Therefore, we have $\mathcal{T}^1 =\{\zeta_1^1,\ldots, \zeta^1_M\}$. We note that all the nodes can be characterized by their $n-$th {\em time level}, as in the following definition.
\begin{dfn}
The general $n$-th level of the tree will be composed by $M^n$  nodes denoted by 
$$\mathcal{T}^n = \{ \zeta^{n-1}_i + \Delta t f(\zeta^{n-1}_i, u_j,t_{n-1}) \}_{j=1}^{M}\quad i = 1,\ldots, M^{n-1}.$$ 
\end{dfn}
We show in the left panel of Figure \ref{fig:tree} the structure of the whole tree $\mathcal{T}$.
All the nodes of the tree can be shortly defined as
 $$\mathcal{T}:= \{ \zeta_j^n  \}_{j=1}^{M^n},\quad n=0,\ldots \overline{N},$$ 
where the nodes $\zeta^n_i$ are the result of the dynamics at time $t_n$ with the controls $\{u_{j_k}\}_{k=0}^{n-1}$:
\begin{align*}
\zeta_{i_n}^n &= \zeta_{i_{n-1}}^{n-1} + \Delta t f(\zeta_{i_{n-1}}^{n-1}, u_{j_{n-1}},t_{n-1})\\
&= x+ \Delta t \sum_{k=0}^{n-1} f(\zeta^k_{i_k}, u_{j_k},t_k), 
\end{align*}
with $\zeta^0 = x$, \rev{$i_k = \ceil[\bigg]{\dfrac{i_{k+1}}{M}}$ and $j_k\equiv i_{k+1} \mbox{mod } M$, where $\ceil{\cdot}$ is the ceiling function.}
We note that $\zeta_i^k \in \R^d, i=1,\ldots, M^k$. On the right panel of Figure \ref{fig:tree} we show the path to reach for instance $\zeta^4_{26}$ if the control set contains only three elements. We, again, would like to emphasize that the domain is not chosen a priori, but constructed following the dynamics.

%

In what follows we provide two remarks about the properties of the tree $\mathcal{T}$ under some particular assumptions on the dynamics $f$.
\begin{rmk}
Let us suppose that the dynamics is \revsecond{affine} with respect to $u$ and that $u\in [u_{min},u_{max}]\subset\R$, e.g. the following decomposition holds true
$$f(x,u,t)=f_1(x,t)u+f_2(x,t).$$
Then, all the nodes in $\mathcal{T}^n$ will lie on the segment with extremal points given by the controls at the boundary $\partial U = \{u_1= u_{min}, u_{M}=u_{max}\}$. Specifically,
$$\mbox{if } z\in\mathcal{T}^n \mbox{ this implies } z\in [\zeta^n_{i_1},\zeta^n_{i_{M}}]$$
where $\zeta^n_{i_1}$ and $\zeta^n_{i_{M}}$ are obtained by using the control $u_1$ and $u_{M}$ respectively.



\end{rmk}

\begin{rmk}
Let us suppose that the dynamics is monotone with respect to $u\in [u_{min},u_{max}]\subset\R$:
\begin{align*}
\min_{\widetilde{u} \in \{u_{min}, u_{max} \}} f_j(x, \widetilde{u},t) \le f_j(x,u,t) \le \max_{\widetilde{u} \in \{u_{min}, u_{max} \}} f_j(x,\widetilde{u},t),\\
\qquad\qquad \forall u \in [u_{min},u_{max}],\, j =1,\ldots, d.
\end{align*}
Then the nodes of the tree will belong to a box with vertices given by the coordinates of the nodes obtained with the extremal controls $u_{min}$ and $u_{max}$ as follows:
$$\min_{\overline{i} \in \{i_1,i_{M} \} }\zeta^n_{\overline{i}} \leq\zeta^n_i \leq \max_{\overline{i} \in \{i_1,i_{M} \} } \zeta^n_{\overline{i}}\,, \quad i \in \{i_1,\ldots, i_{M}\},$$
where the last inequality holds component-wise.
%
\end{rmk}

\paragraph{\bf Approximation of the value function}

The numerical value function $V(x,t)$ will be computed on the tree nodes in space, whereas in time it will approximated as a piecewise constant function, $i.e.$
$$
V(x,t)=V^n(x) \quad \forall x, \mbox{ and } t \in [t_n,t_{n+1}),
$$
where $t_n=t+ n \Delta t$.

We note that we start to approximate the value function once the tree $\mathcal{T}$ has been already built. Then, we will be able to approximate the value function $V^{n}(x_i+\Delta t f(x_i,u,t_n))$ in \eqref{SL} without the use of an interpolation operator on a grid. The reason is that we build our domain according to all the possible directions of the dynamics for a discrete set of controls and, as a consequence, all the nodes $x_i+\Delta t f(x_i,u,t_n)$ will belong to the grid. It is now straightforward to evaluate the value function.
The TSA defines a grid $\mathcal{T}^n=\{\zeta^n_j\}_{j=1}^{M^n}$ for $n=0,\ldots, \overline{N}$,
we can approximate \eqref{HJB} as follows: 
\begin{equation}
\begin{cases}
V^{n}(\zeta^n_i)= \min\limits_{u\in U} \{e^{-\lambda \Delta t} V^{n+1}(\zeta^n_i+\Delta t f(\zeta^n_i,u,t_n)) +\Delta t \, L(\zeta^n_i,u,t_n) \}, \\
\qquad\qquad \qquad\qquad \qquad\qquad \qquad\qquad \qquad\qquad  \zeta^n_i \in \mathcal{T}^n\,, n = \overline{N}-1,\ldots, 0, \\
V^{\overline{N}}(\zeta^{\overline{N}}_i)= g(\zeta_i^{\overline{N}}), \\\qquad\qquad \qquad\qquad \qquad\qquad \qquad\qquad \qquad\qquad   \zeta_i^{\overline{N}} \in \mathcal{T}^{\overline{N}}.
\end{cases}
\label{HJBt2}
\end{equation}

We note that the minimization is computed by comparison on the discretized set of controls $U$. We refer to \cite{B73,KKK16} \revfirst{ for a more sophisticated approach} to compute the minimum in \eqref{SL}. \\

\begin{rmk}\label{dyn_auto}
If the dynamics \eqref{eq} is autonomous the evolution of the dynamics will not depend explicitly on $t_n$ and the problem can be simplified since the argument of the minimization in \eqref{HJBt2} will be
$$
e^{-\lambda \Delta t}  V^{n+1}(\zeta+\Delta t f(\zeta,u)) +\Delta t \, L(\zeta,u,t_n).
$$
At the time $t_n$ we have $n$ levels of the tree on the left and $\overline{N} -n$ levels on the right (till $t_{\overline{N}}$). Since the computation is going backward, to compute the value function at time $t_n$, we need to do $\overline{N}-n$ steps in time starting from the final condition at time T. Once we know $V^n$ this information can also be interpreted as a final condition for the sub-tree $\cup_{k=0}^n \mathcal{T}^k$ and, since the dynamics is autonomous, we can proceed backward computing $V^{n-1}$  for the nodes belonging to all the $k-$th time levels, for $k \le n-1$. Indeed the nodes $\zeta+\Delta t f(\zeta,u)$ do not depend explicitly on the time and they can be involved in the computation of the value function at different time steps (this is not the case for a non-autonomous dynamics). Thus, we will proceed as follows: first we impose the final cost $g$ on the whole tree, then we start computing the value function backward. This procedure leads to a more extensive knowledge of the value function on the tree.

\end{rmk}

\section{Hints on the algorithm}\label{sec:hints}
\setcounter{algorithm}{1}
In this section we will provide further details on the implementation of the method proposed in Section \ref{sec:tree}. We will explain how to reduce the number of tree nodes to make the problem feasible, compute the feedback control and recall the whole procedure.
\paragraph{\bf Pruning the tree}
The proposed method mitigates the curse of dimensionality and it allows to deal with problems in $\mathbb{R}^d$ with $d\gg5$, which is absolutely not feasible with the classical approach. However, we still have dimensionality problem related to the amount of nodes in the tree $\mathcal{T}$. In fact, given $M>1$ controls and $\overline{N}$ time steps, the cardinality of the tree is
$$
|\mathcal{T}| = \sum_{i=0}^{\overline{N}} M^{i}= \frac{M^{\overline{N}+1}-1}{M-1},
$$
which is infeasible due to the huge amount of memory allocations, if $M$ or $\overline{N}$ are too large. 
Therefore, we suggest to select the nodes of the TSA neglecting those very close to each other, assuming that the value function will not be completely different on those nodes, e.g.
$$\zeta^n_i\approx \zeta^{\new{n}}_j \Longrightarrow V(\zeta^n_i)\approx V(\zeta^{\new{n}}_j).$$

This is a realistic assumption since the numerical value function is Lipschitz continuous as explained in Proposition \ref{prp:lip}. We can introduce the {\em pruning rule.}
\begin{dfn}[Pruning rule]
Two given nodes $\zeta^n_i$ and $\zeta^{\new{n}}_j$ can be merged if 
\begin{equation}\label{tol_cri}
\Vert \zeta^n_i-\zeta^{\new{n}}_j \Vert \le \ep, \quad \mbox{ with } n = 0,\ldots, \overline{N}, 
\end{equation}
for a given threshold $\ep>0$.
\end{dfn}

Specifically, if during the construction of the tree, a node $\zeta^{\new{n-1}}$ has as a son a new node $\zeta_j^{n}$ which verifies \eqref{tol_cri} with a certain $\zeta_i^n$, then we will not add the new node to the tree and we will connect the node $\zeta^{\new{n-1}}$ with $\zeta_i^n$.
We cut the node which verifies the criteria before going on with the construction of the tree, in this way we avoid the sub-tree coming out from this node, saving a huge amount of memory.

\begin{figure}[htbp]
\centering
       \includegraphics[scale=0.5]{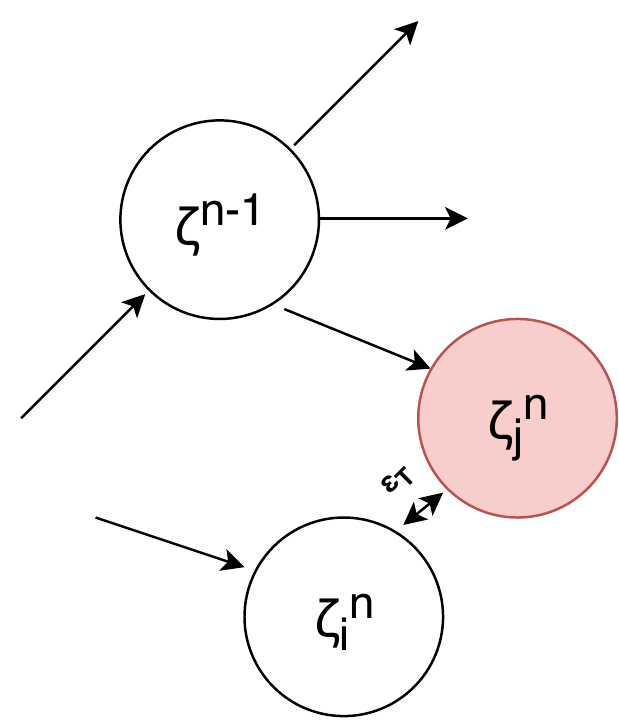} 
       \hspace{1.5cm}
       \includegraphics[scale=0.5]{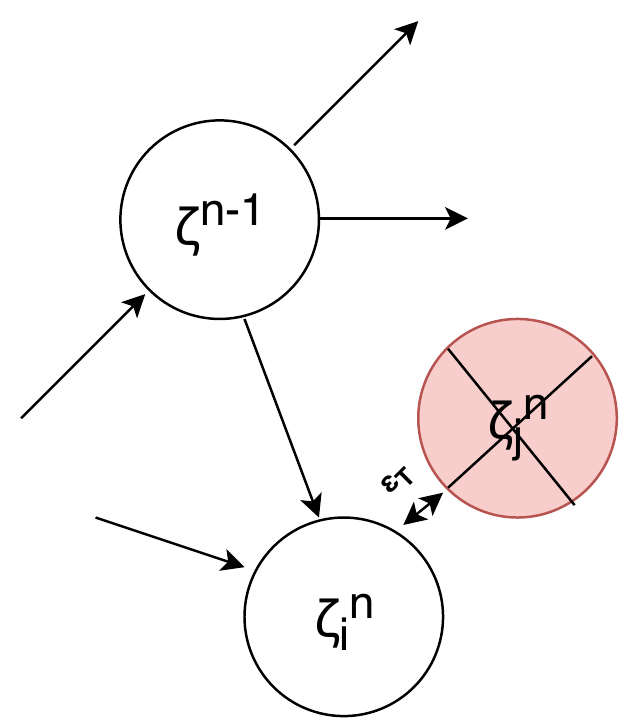} 
       \caption{Pruning technique throughout the construction of the tree it might happen the two nodes are very close (left) and we link those node in order to prune the tree. } 
       \label{fig:prun}
	\end{figure}

The cut of the tree works as follows: during the construction of the $\new{n}$-th level, the new node will be compared with the previous nodes \revfirst{already computed at the same level $n$}. If the new node $\zeta^{\new{n}}_j$, whose father is $\zeta^{\new{n-1}}$, satisfies the condition \eqref{tol_cri} with a node $\zeta^n_i$, the new node will not be added to the tree and the adjacency matrix will be uploaded, connecting the node $\zeta^{\new{n-1}}$ to the node $\zeta^n_i$. Figure \ref{fig:prun} provides a graphic idea about the application of the pruning criteria.
The choice of the tolerance plays an important role: if $\ep$ is very small, the algorithm will be very slow, whereas if it is too large, we will not obtain an accurate approximation. A reasonable choice turns out to be $\ep = \Delta t^2$, as shown in Section \ref{sec:nt}. \revfirst{The interested reader will find a rigorous proof of this heuristic statement in \cite{SAF18} together with convergence results of the proposed method.}
\new{\begin{rmk}[Pruning rule in the autonomous case]
If the dynamics is autonomous, as explained in Remark \ref{dyn_auto}, we can extend the computation of the value function at time $t_n$ even for nodes belonging to the subtree $\cup_{k=0}^n \mathcal{T}^k$. Therefore, we can extend the pruning criteria \eqref{tol_cri} as follows. Two given nodes $\zeta^n_i$ and $\zeta^m_j$ can be merged if 
\begin{equation}\label{tol_cri_auto}
\Vert \zeta^n_i-\zeta^m_j \Vert \le \ep, \quad \mbox{ with } n,m = 0,\ldots, \overline{N}, 
\end{equation}
for a given threshold $\ep>0$.
\end{rmk}
\begin{rmk}[Efficient Pruning]
\revfirst{
The computation of the distances among all the nodes would be very expensive, especially for high dimensional problems. Hence, we need an efficient algorithm to compute the distances quickly. One possible strategy is the Principal Analysis Component (\cite{P01},\cite{J02}). Our aim is to project the data onto a lower dimensional linear space such that the variance of the projected data is maximized. This can be done e.g. computing  the Singular Value Decomposition of the data matrix and taking the first basis. Once we project the data, the distances will be computed in a lower dimension space and this turns out to accelerate the algorithm.}
\end{rmk}}

\paragraph{\bf Feedback reconstruction and closed-loop control}
During the computation of the value function, we store the control indices corresponding to the argmin in \eqref{HJBt2}. Then starting from $\zeta^0_*= x$, we follow the path of the tree to build the optimal trajectory $\{\zeta^n_*\}_{n=0}^{\overline{N}}$ in the following way
\begin{equation} \label{feed:tree}
u_{n}^{*}:=\argmin\limits_{u\in U} \left\{ e^{-\lambda \Delta t}V^{n+1}(\zeta^n_*+\Delta t f(\zeta^n_*,u,t_n)) +\Delta t \, L(\zeta^n_*,u,t_n) \right\},
\end{equation}
\begin{equation*} 
\zeta^{n+1}_* \in \mathcal{T}^{n+1} \; s.t. \; \zeta^n_* \rightarrow^{u_{n}^{*}} \zeta^{n+1}_*,
\end{equation*}
for $n=0,\ldots, \overline{N}-1$, where the symbol $\rightarrow^u$ stands for the connection of two nodes by the control $u$.  We note that this is possible because in the current work we assume to consider the same discrete control set $U$ for both HJB equation \eqref{HJBt2} and feedback reconstruction \eqref{feed:tree}.

\paragraph{\bf Algorithm}
In what follows we summarize the whole algorithm including the construction of the tree, the selection of the nodes and, finally, the approximation of the value function.

\begin{algorithm}[H]
\caption{TSA algorithm with pruning}
\begin{algorithmic}[1]
\State $\mathcal{T}^0 \gets x$
\label{alg_1}
\For{$n=1,...,\overline{N}$}
 \For{$u_j \in U$, $\zeta^{n-1} \in \mathcal{T}^{n-1}$}
 \State{$\zeta_{new}=\zeta^{n-1}+ \Delta t f(\zeta^{n-1},u_j,t_{n-1})$}
\If{$\Vert \zeta_{new} - \zeta \Vert > \ep, \forall \zeta \in \mathcal{T}$ }
\State{$\mathcal{T}^n \gets \zeta_{new}$}
\State{$\zeta^{n-1} \rightarrow^u \zeta_{new} $}
\Else
   \State{$\overline{\zeta}= arg\,min_{\zeta \in \mathcal{T}}\Vert \zeta_{new} - \zeta \Vert$}
   \State{$\zeta^{n-1} \rightarrow^u \overline{\zeta}$}
\EndIf
\EndFor
\EndFor
\State $V^{\overline{N}}(\zeta)=g(\zeta), \forall \zeta \in \mathcal{T}^{\overline{N}}$
\For{$n=\overline{N}-1,...,0$}
\State $V^{n}(\zeta^n)= \min\limits_{\zeta^{n+1}: \zeta^n \rightarrow^u \zeta^{n+1} } \{e^{-\lambda \Delta t} V^{n+1}(\zeta^{n+1}) +\Delta t \, L(\zeta^n, u, t_n) \},\quad \zeta^n \in \mathcal{T}^n$. 
\EndFor
\end{algorithmic}
\end{algorithm}

As one can see in Algorithm \ref{alg_1}, we first start the construction of the tree $\mathcal{T}$ from 1 to step 10. We note that the pruning criteria is involved in the steps 5-10 of Algorithm \ref{alg_1}. Clearly, a very small tolerance will not allow any selection of the nodes and we will work with a full tree. Finally, in step 11-12-13 we compute the approximation of the value function. In the last step, the computation of the value function $V^{n}(\zeta^n)$ can be extended to the nodes in the tree $\cup_{k=0}^{n}\mathcal{T}^k$ in the case of autonomous dynamics.

\section{Numerical tests}\label{sec:nt}
In this section we are going to apply the proposed algorithm to show the effectiveness of the method. 

We will present five test cases. In the first we are able to compute the analytical solution of the HJB equation and, therefore, to compute the error with our method compared to the classical approach, see e.g. \cite{FG99}. The second test concerns the well-known Van der Pol equation and the third is \revfirst{ about non-autonomous dynamics}. Finally we present the results for two different linear PDEs which shows the power of the method even for large-scale problems.

The numerical simulations reported in this paper are performed on a laptop with 1CPU Intel Core i5-3,1 GHz and 8GB RAM. The codes are written in C++. 

\subsection{Test 1: Comparison with exact solution of the value function}

\begin{figure}[htbp]
\centering
\includegraphics[scale=0.25]{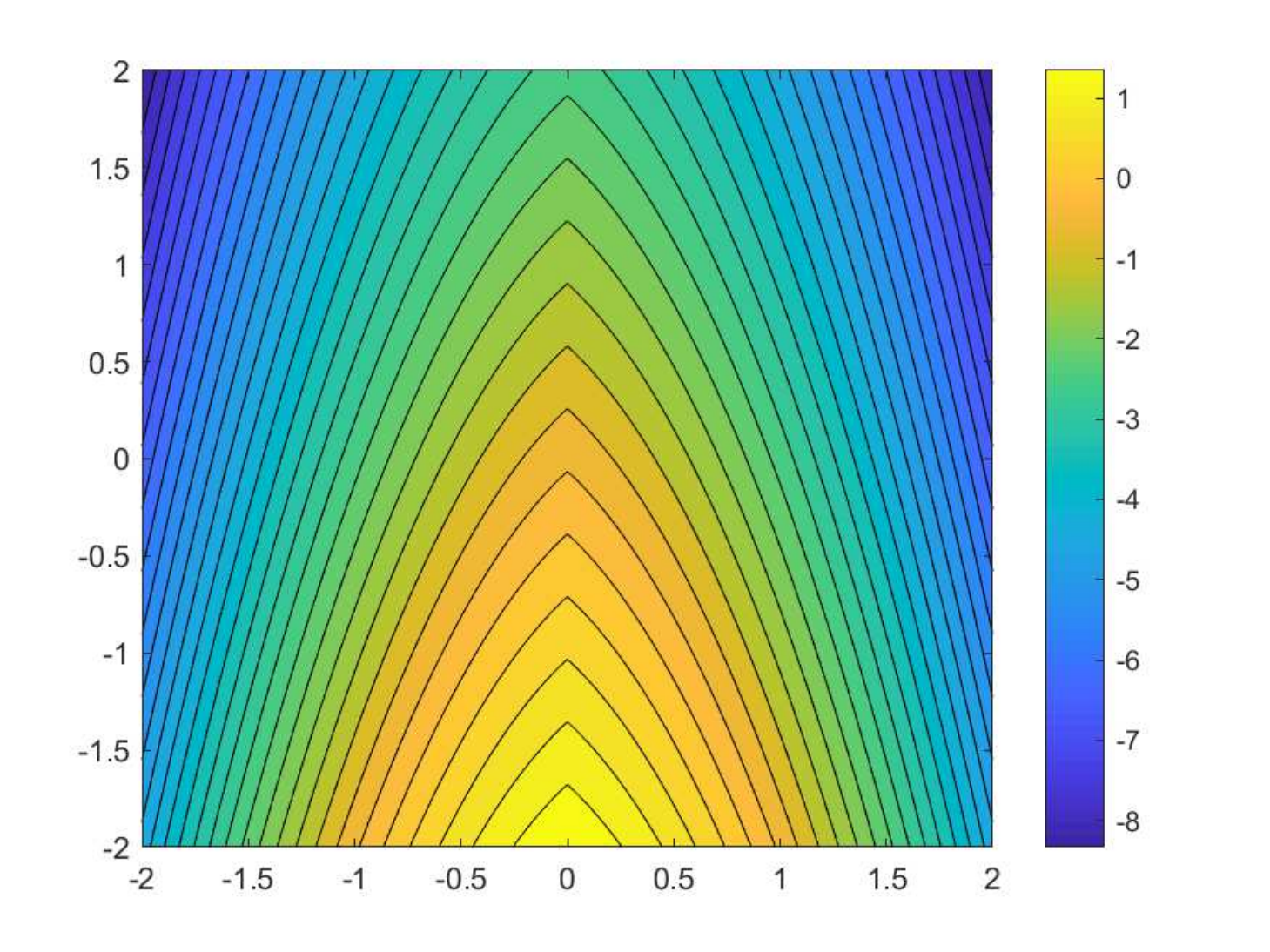}
\includegraphics[scale=0.25]{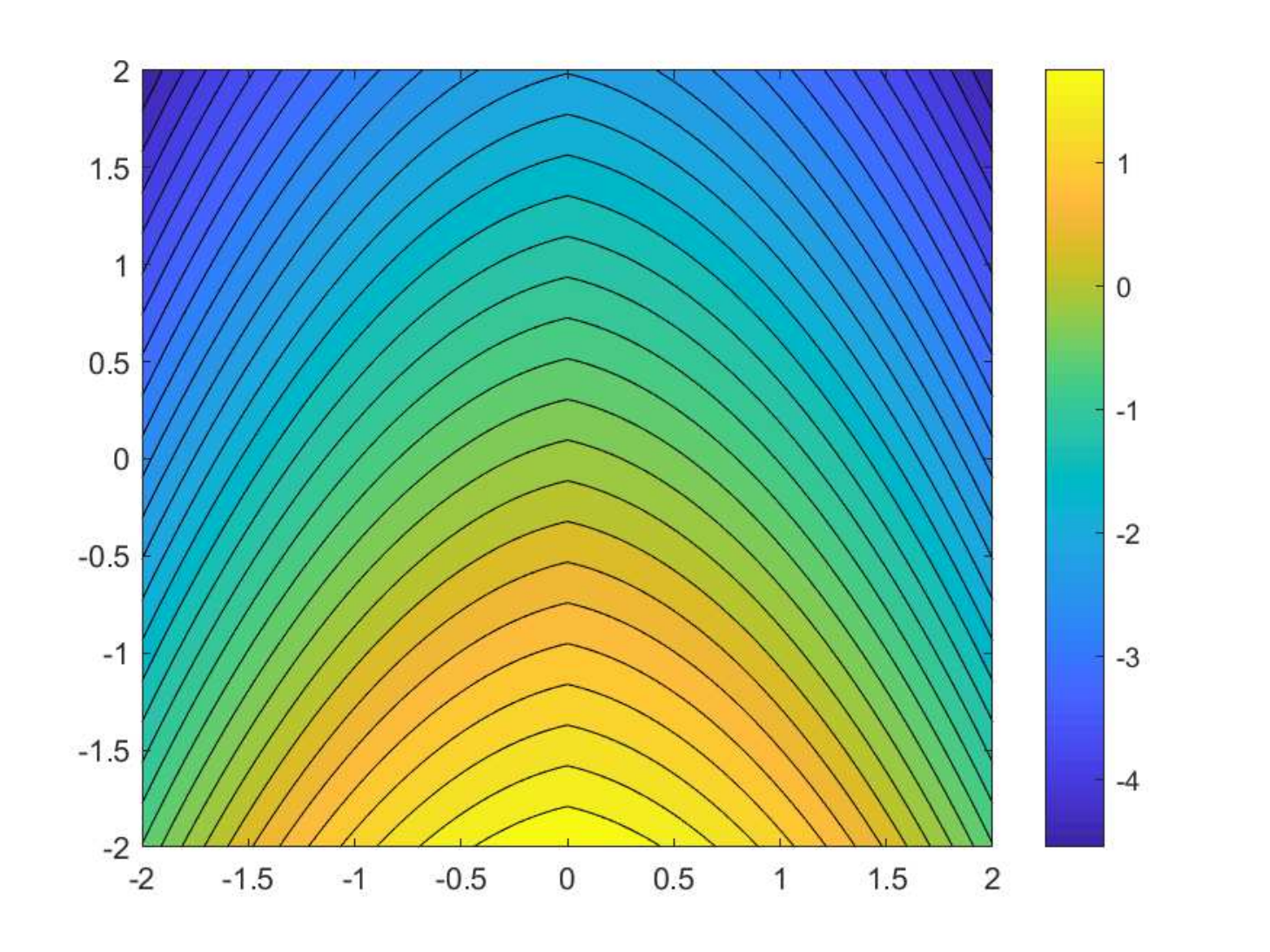}
\includegraphics[scale=0.25]{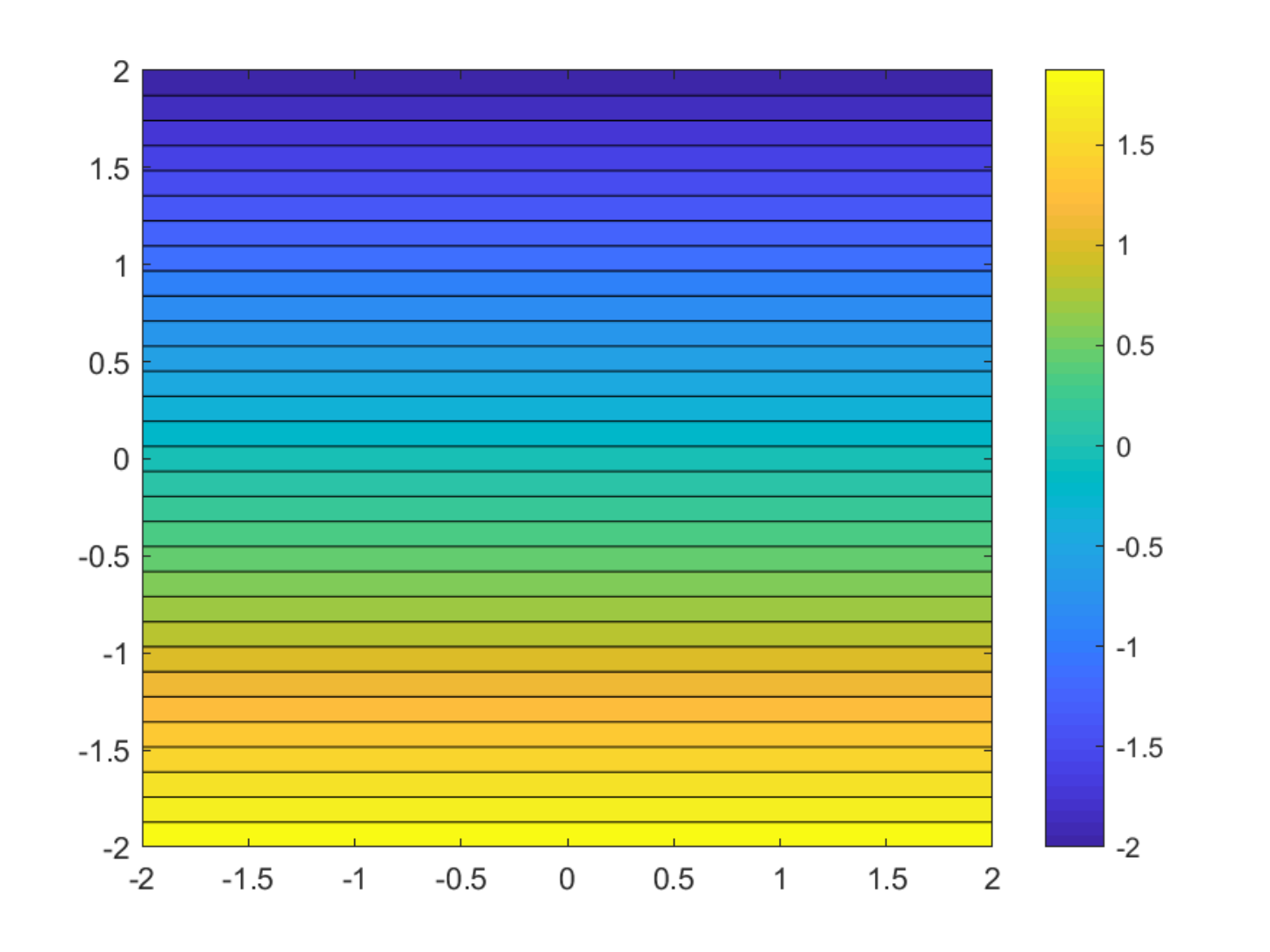}
\caption{Test 1: Contour lines for \eqref{VF_given} with $t=0$ (left), $t=0.5$ (middle) and $t=1$ (right).}
\label{fig1:vf}
\end{figure}

In the first example we consider the following dynamics in \eqref{eq}
\begin{equation}
	f(x,u)=
		  \begin{pmatrix}
         u \\
         x_1^2
          \end{pmatrix},\, u \in U\equiv [-1,1],
       \label{exact}   
\end{equation}
where $x=(x_1,x_2)\in\R^2.$
The cost functional in \eqref{cost} is: 
\begin{equation}
L(x,u,t) = 0,\qquad g(x) =-x_2, \qquad \lambda = 0,
\end{equation}
where we only consider the terminal cost $g$. 
The corresponding HJB equation is
\begin{equation}\label{HJB1}
\begin{cases}
-V_t+|V_{x_1}|-x_1^2V_{x_2}=0  & (x,t) \in \R^2\times[0,T], \\
V(x,T)=g(x)\;, & 
\end{cases}
\end{equation}
where its unique viscosity solution reads
\begin{equation}\label{VF_given}
V(x,t)=-x_2-x_1^2(T-t)-\frac{1}{3}(T-t)^3-|x_1|(T-t)^2.
\end{equation}
Furthermore, we set $T=1$. Figure \ref{fig1:vf} shows the contour lines of the value function $V(x,t)$ in \eqref{VF_given} for time instances $t=\{0,0.5,1\}$. 

In this example, we compare the classical approach with the TSA algorithm proposed in Algorithm \ref{alg_1} using both strategies: (i) no selection of the nodes and (ii) applying criteria \eqref{tol_cri} to select the nodes as explained in Section \ref{sec:hints}. To perform a fair comparison we projected the value function computed with the classical method into the tree nodes. We note that it will not modify the accuracy of the classical approach since the interpolation has to be performed also on a structured grid. 
We compare the different approximations according to $\ell_2-$relative error with the exact solution on the tree nodes
$$
\mathcal{E}_2(t_n)= \sqrt{ \frac{ \sum\limits_{x_i \in \mathcal{T}^n } |v(x_i,t_n)-V^n(x_i)|^2}{  \sum\limits_{x_i \in \mathcal{T}^n }|v(x_i,t_n)|^2}},
$$
where $v(x_i,t_n)$ represents the analytical solution and $V^n(x_i)$ its numerical approximation.

\begin{figure}[htbp]
 \centering
        \includegraphics[scale=0.4]{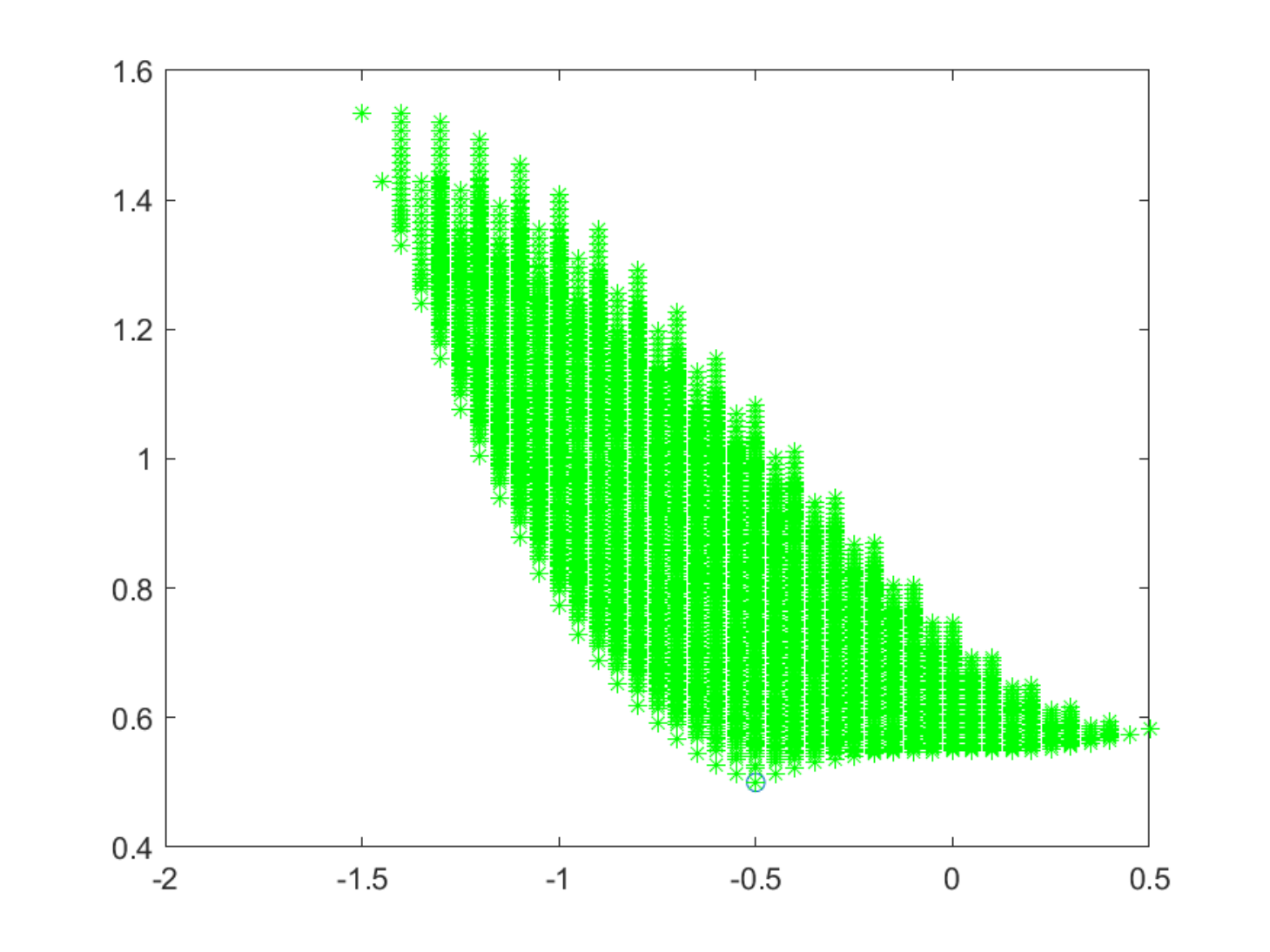} 
          \includegraphics[scale=0.4]{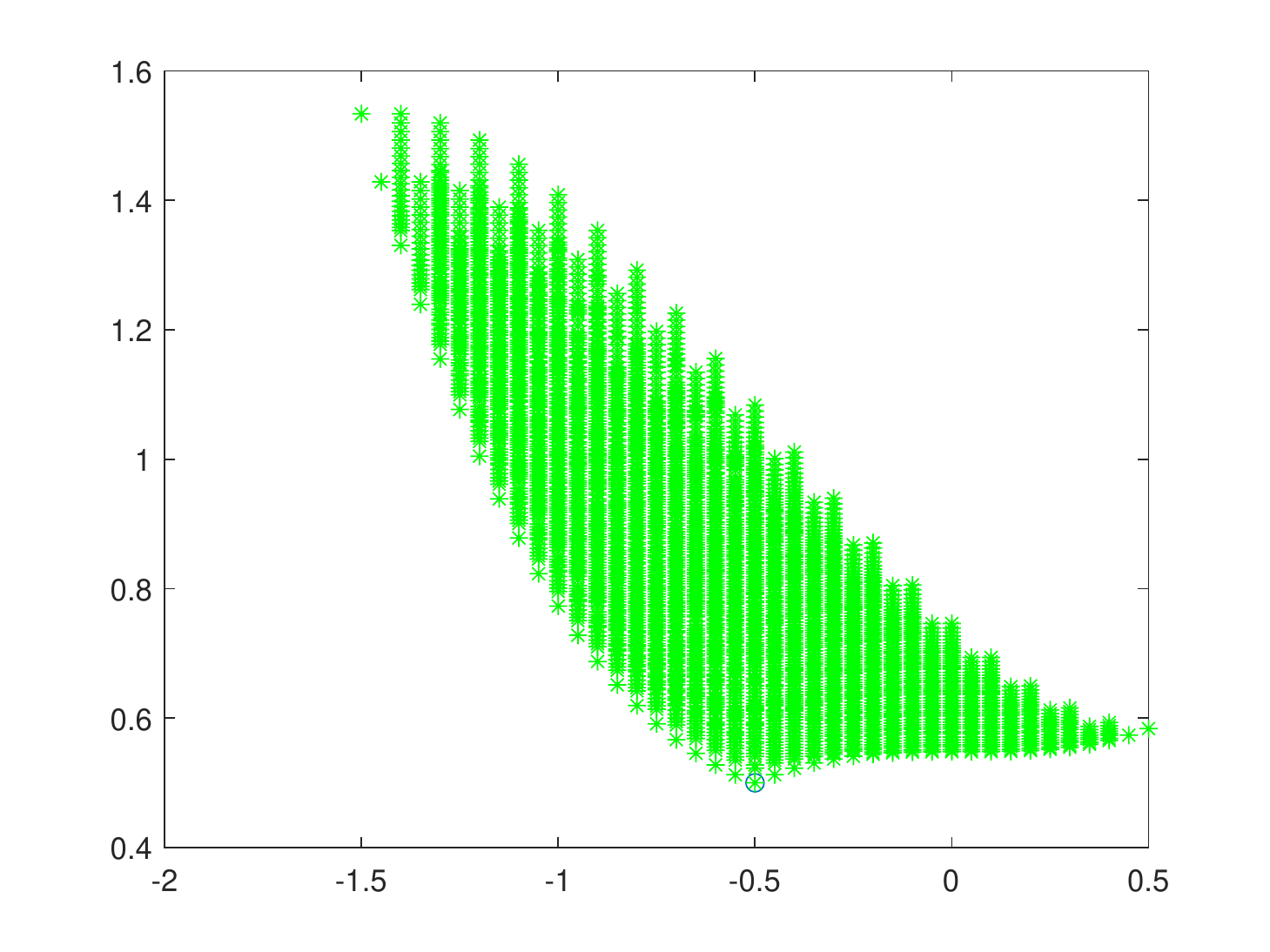} 
       \caption{Test 1:  Tree nodes without tolerance (left) and with tolerance equal to $\Delta t^2$ (right) for $x =(-0.5,0,5)$.   }
       \label{fig1:tree}
	\end{figure}

In Figure \ref{fig1:tree}, we show all the nodes of the tree $\mathcal{T}$ for the initial condition  $x=(-0.5,0.5),$ $\Delta t=0.05$ and different choices of $\ep =\{0, \Delta t^2\}$. We note that there is a huge difference between the cardinality of the trees, that is $|\mathcal{T}| = 2097151$ when the tolerance is not applied whereas we have $|\mathcal{T}| = 3151 $ for $\ep = \Delta t^2$. 

\begin{figure}[htbp]
 \centering
        \includegraphics[scale=0.4]{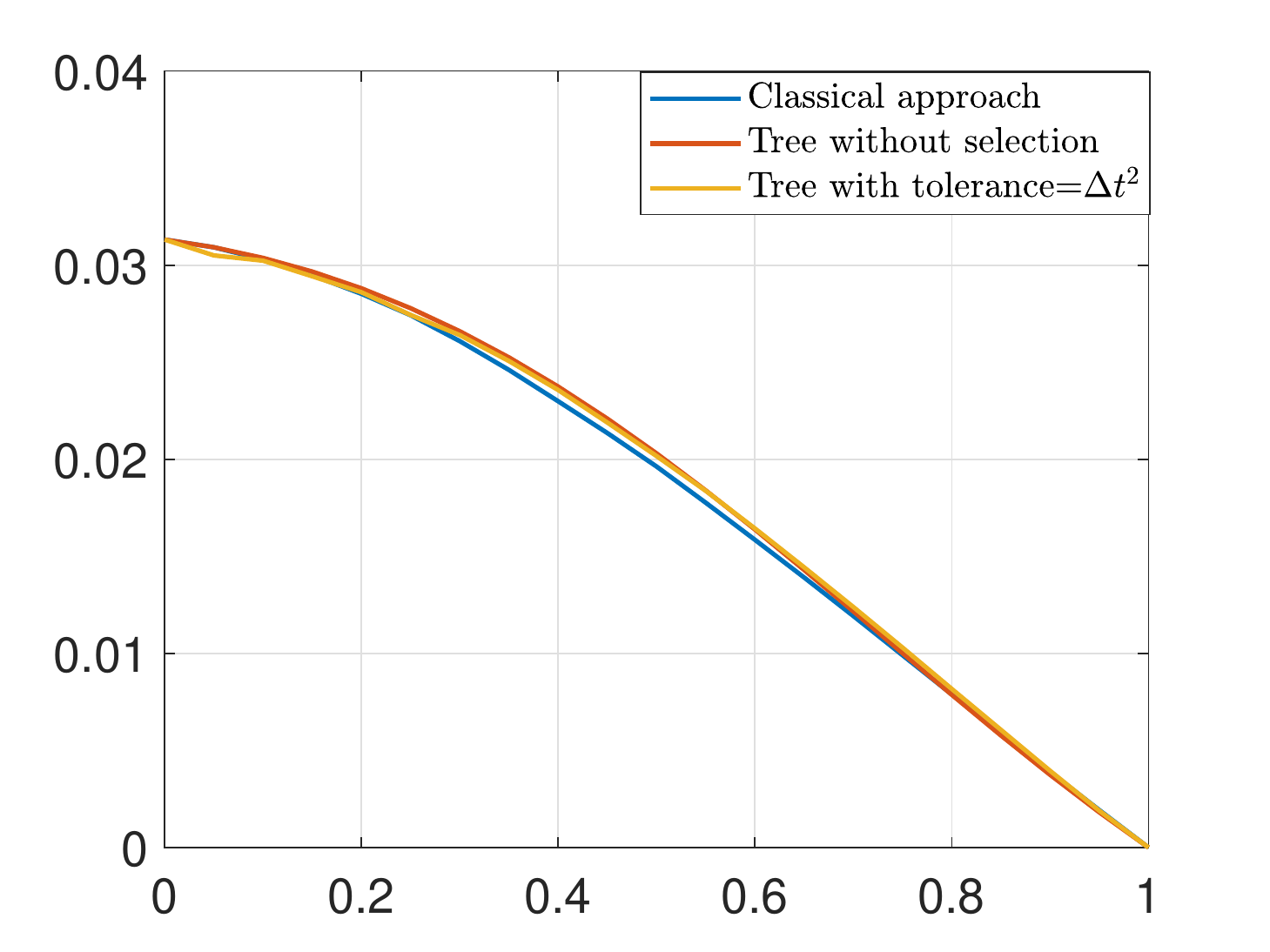} 
        \includegraphics[scale=0.4]{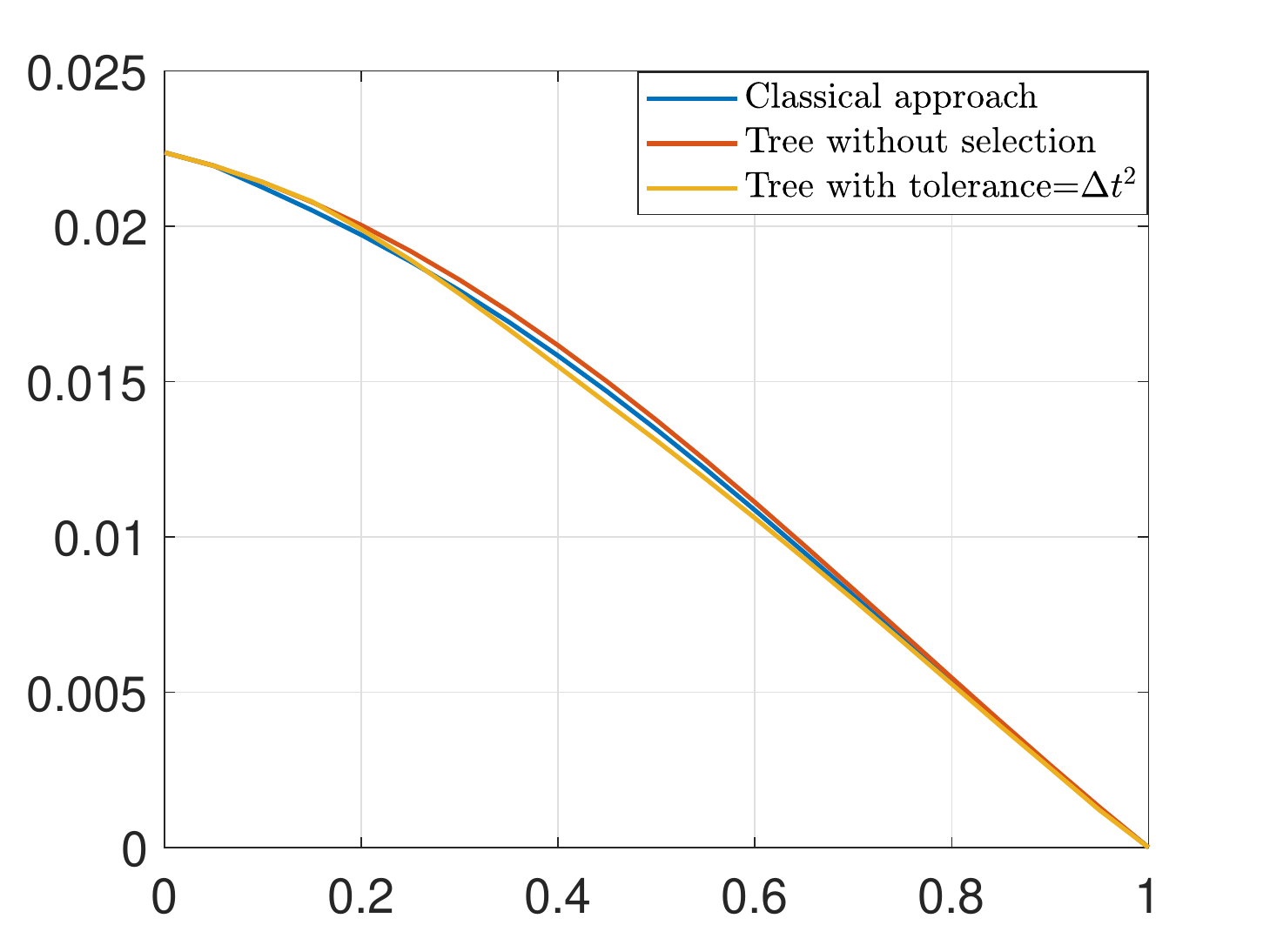} 
       \caption{Test 1: Comparison of the different methods with initial datum $(-0.5,0.5)$ (left) and with initial datum $(1,1)$ (right) for each time instance ($x$-axis). }  
       \label{fig1:err}
	\end{figure}

In Figure \ref{fig1:err}, we show the behaviour of the error $\mathcal{E}_2$ for two different initial conditions $x$. We note that its behaviour is very similar using both the classical approach and the TSA with or without the pruning criteria \eqref{tol_cri} for the nodes. As already mentioned, we would like to stress that the domain for the solution of the classical approach is chosen as large as possible to avoid that the boundary conditions are active, whereas with TSA we do not have this kind of problem, since the domain of the tree constructed according to the vector field. We note that to compute the value function in the classical approach we use the following step size: $\Delta x= \Delta t=0.05$.

In Table \ref{test1:tab_nosel} we show the error decay decreasing the temporal step size $\Delta t$ for $x=(-0.5,0.5)$ and $\ep=0$ (i.e no pruning criteria has been applied). We compute the error as follows:

$$Err_{2,2}= \sqrt{ \Delta t\sum_{n=0}^{\overline{N}} \mathcal{E}^2_2(t_n)}, \quad Err_{\infty,2}=\max_{n=0,...,\overline{N}} \mathcal{E}_2(t_n) $$
and the order
$$
Order_{2,2}=log_2 \left( \frac{Err_{2,2}(\Delta t)}{Err_{2,2}(\Delta t/2)} \right), \quad Order_{\infty,2}=log_2 \left( \frac{Err_{\infty,2}(\Delta t)}{Err_{\infty,2}(\Delta t/2)} \right).
$$

We note that the order of convergence is linear as the order of the method used to discretize the dynamics \eqref{eq}, e.g. forward Euler scheme. This feature will be analyzed in a follow-up paper where we would like to provide error estimate for our proposed algorithm.

 
	\begin{table}[htbp]
		\centering
		\begin{tabular}{cccccccc}
					\hline 
					$\Delta t$ & $|\mathcal{T}|$ & CPU & $Err_{2,2}$ &  $Err_{\infty,2} $ & $Order_{2,2}$  & $Order_{\infty,2}$  \\
			\hline
	\phantom{x}0.2 & \phantom{xxxx}63 &  0.05s & 0.090 & 0.122 & & \\ 
			
		\phantom{x}0.1 &	\phantom{xx}2047 &   0.35s & 0.044 & 0.062 &  1.04 & 0.98\\ 
			 
			0.05 &	2097151 &  1.1s &  0.022 & 0.031 & 1.02 & 0.99 \\ 
		\end{tabular}
   	\caption{Test 1: Error analysis and order of convergence of the TSA without pruning rule.}
	\label{test1:tab_nosel}
	\end{table}
However, the case without selection is quite unfeasible for more than 20 time steps since it requires to store a huge amount of nodes of order $O(M^{21})$, whereas with the selection we can obtain an impressive improvement. The results are shown in Table \ref{test1:tab_sel} where we can see, although the pruning of the nodes, we are still able to achieve an order of convergence close to $1$.

\begin{figure}[htbp]
 \centering
   \includegraphics[scale=0.4]{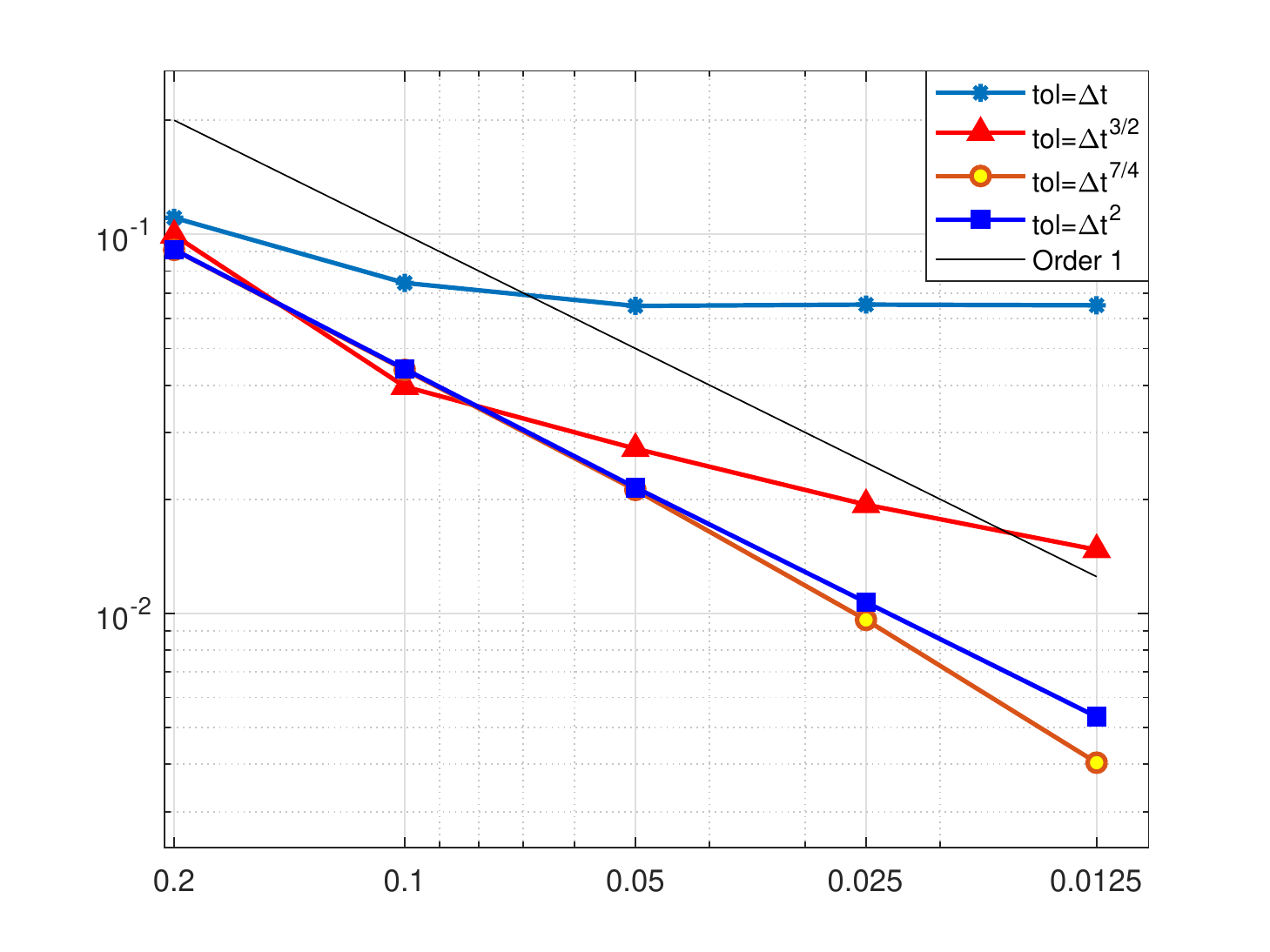}
  \includegraphics[scale=0.4]{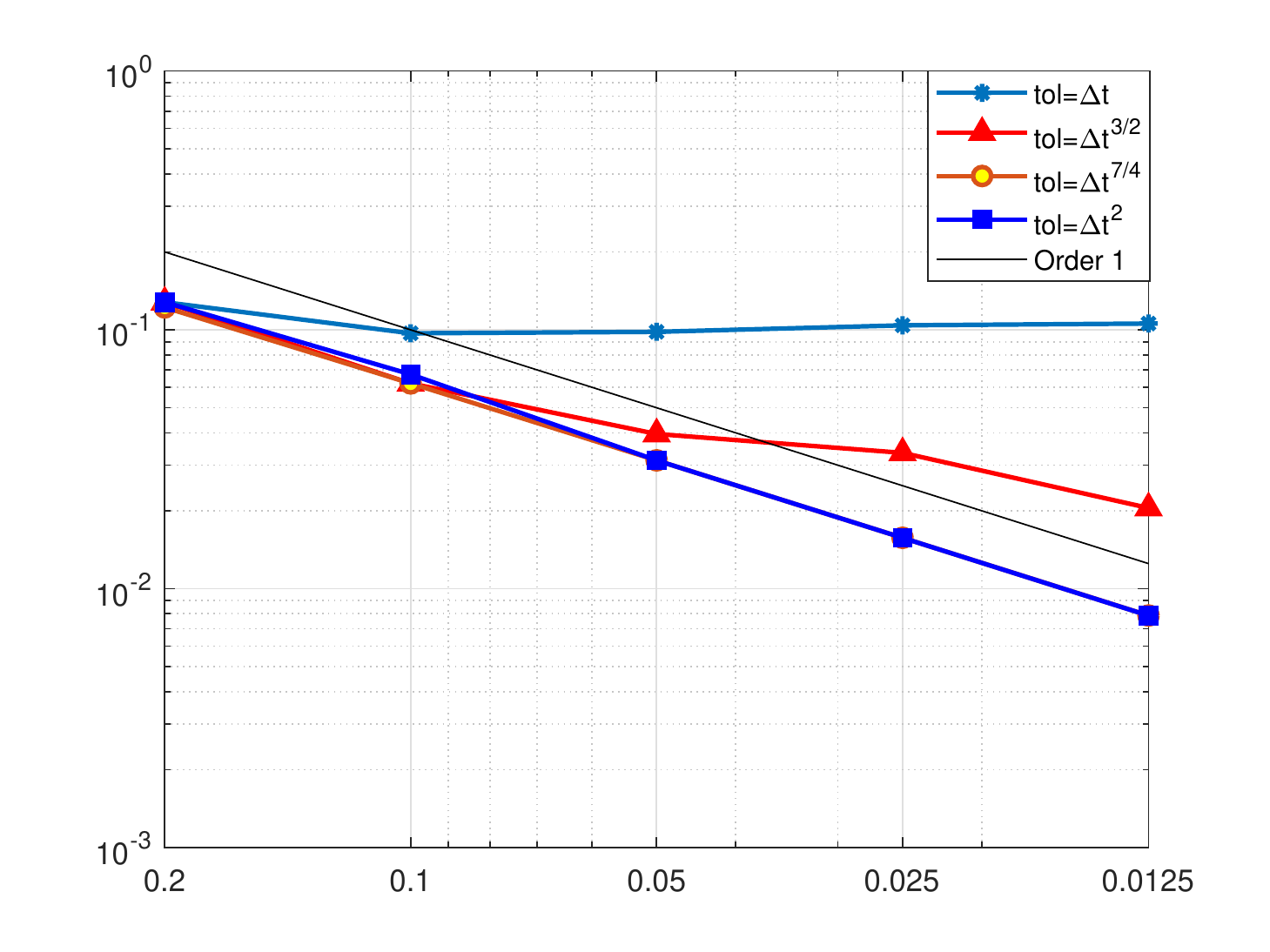}           
  \caption{Test 1: Comparison of the error $Err_{2,2}$ (left) and the error $Err_{\infty,2}$ (right) for the pruned TSA with different tolerances $\ep$ as a function of $\Delta t$}  
       \label{fig1:testsol}
	\end{figure}

	\begin{table}[hbtp]
		\centering
		\begin{tabular}{cccccccc}
					\hline $\Delta t$ & $|\mathcal{T}|$ & CPU & $Err_{2,2}$ &  $Err_{\infty,2} $ & $Order_{2,2}$  & $Order_{\infty,2}$  \\
					\hline
	\phantom{xxx}0.2 & \phantom{xxx}42 &  \phantom{xxx}0.05s & 0.091 & 0.122 & & \\ 
	\phantom{xxx}0.1 &	\phantom{xx}324 &  \phantom{xxx}0.08s & 0.044 & 0.062 & \phantom{x}1.05 & \phantom{x}0.98 \\ 
	\phantom{xx}0.05 &	\phantom{xx}3151 &  \phantom{xxxx}0.1s & 0.021 & 0.031 &  \phantom{x}1.04 & \phantom{x}0.99 \\ 
		\phantom{x}0.025 &	\phantom{x}29248 &  \phantom{xxxx}0.5s & 0.011 & 0.016 & 1.005  & 0.994 \\ 
			0.0125 &	252620 &  \phantom{xxx} 10s & 0.005 & 0.008 & 1.004  & 0.997 \\ 			
		\end{tabular}
   	\caption{Test 1: Error analysis and order of convergence of the TSA with $\ep= \Delta t^2$ and $T=1$.}
	\label{test1:tab_sel}
	\end{table}
%
%
	\begin{table}[hbtp]
		\centering
		\begin{tabular}{cccccccc}
					\hline $\Delta t$ & $|\mathcal{T}|$ & CPU & $Err_{2,2}$ &  $Err_{\infty,2} $ & $Order_{2,2}$  & $Order_{\infty,2}$  \\
					\hline
	\phantom{xxx}0.2 & \phantom{xxx}1420 &  \phantom{xxx}0.2s & 0.124 & 0.088 & & \\ 
	\phantom{xxx}0.1 & 	\phantom{xxx}15231 &  \phantom{xx}0.11s & 0.061 & 0.045 & \phantom{x}1.02 & \phantom{x}0.98 \\ 
	\phantom{xx}0.05 &  \phantom{xx}141142 &  \phantom{xxxxx}4s &  0.030 & 0.022 &  \phantom{x}1.03 & \phantom{x}1.01 \\ 
		\phantom{x}0.025 &   \phantom{x}1204637 &  \phantom{xxx}147s & 0.015 & 0.011 & 1.009  & 1.002 \\ 
			0.0125 &  	10037898 &  \phantom{x} 7171s & 0.007 & 0.006 & 1.009  & 1.004 \\ 			
		\end{tabular}
   	\caption{\revfirst{Test 1: Error analysis and order of convergence of the TSA with $\ep= \Delta t^2$ and $T=3$.}}
	\label{test1:tab_sel3}
	\end{table}
The tolerance $\ep$ has been set equal to $\Delta t^2$ to keep the same order of convergence of the algorithm as the one without pruning. This is shown in Figure \ref{fig1:testsol}, where we compare the orders of the method with different tolerances. We note that we need to reduce the tolerance to $\Delta t^2$ to ensure linear convergence.

\revfirst{Furthermore, the pruned TSA allows to approximate HJB equation with rather small $\Delta t$ and large horizon, e.g. $T=3$ in a fast way, as shown in Table \ref{test1:tab_sel3} keeping the order of convergence found in the previous case.}
Finally, for the sake of completeness we would like to mention that similar convergence results have been achieved even for other initial conditions $x$.
\subsection{Test 2: Van der Pol oscillator}
\begin{figure}[htbp]
\centering
\includegraphics[scale=0.4]{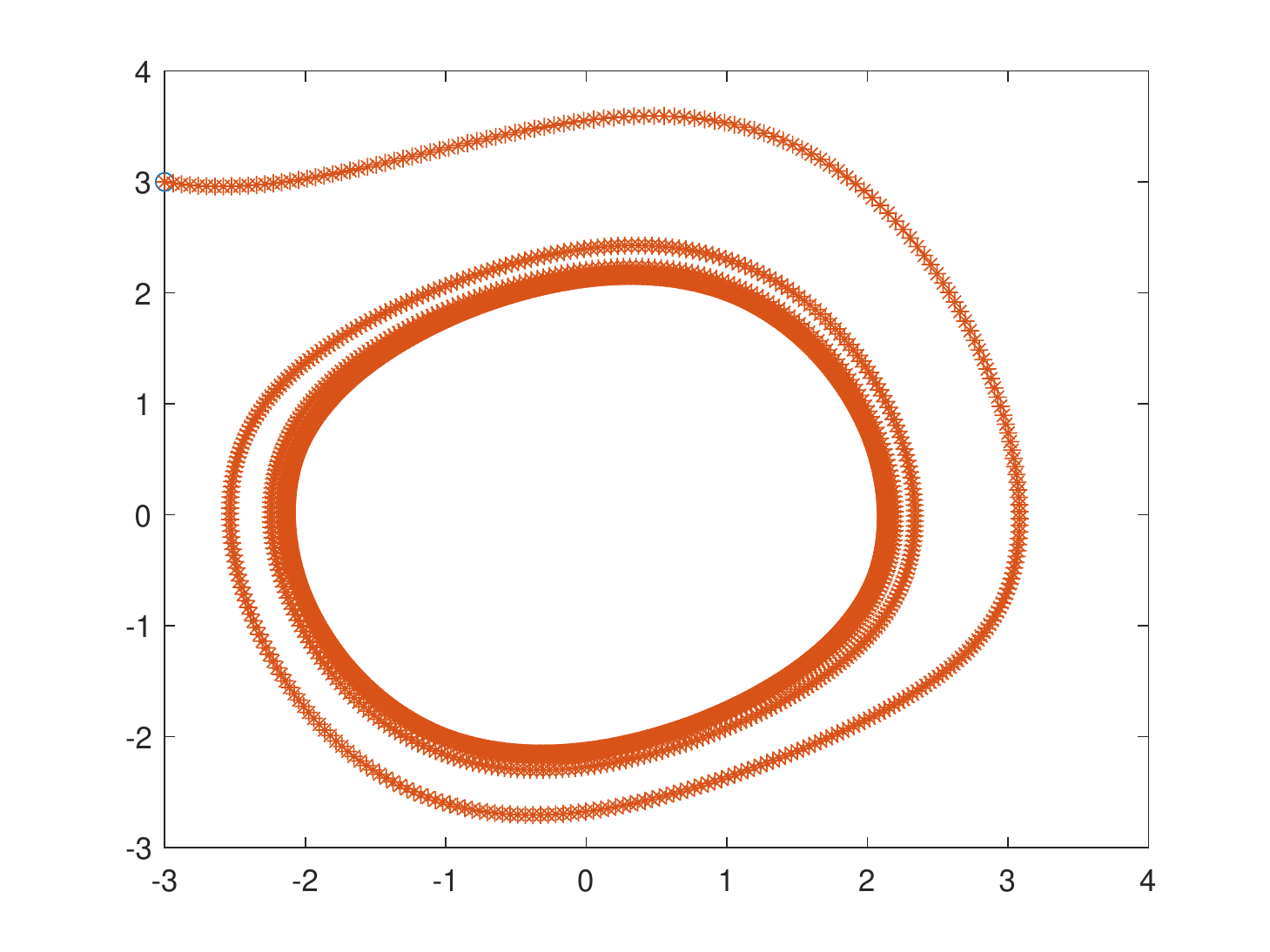}
\includegraphics[scale=0.4]{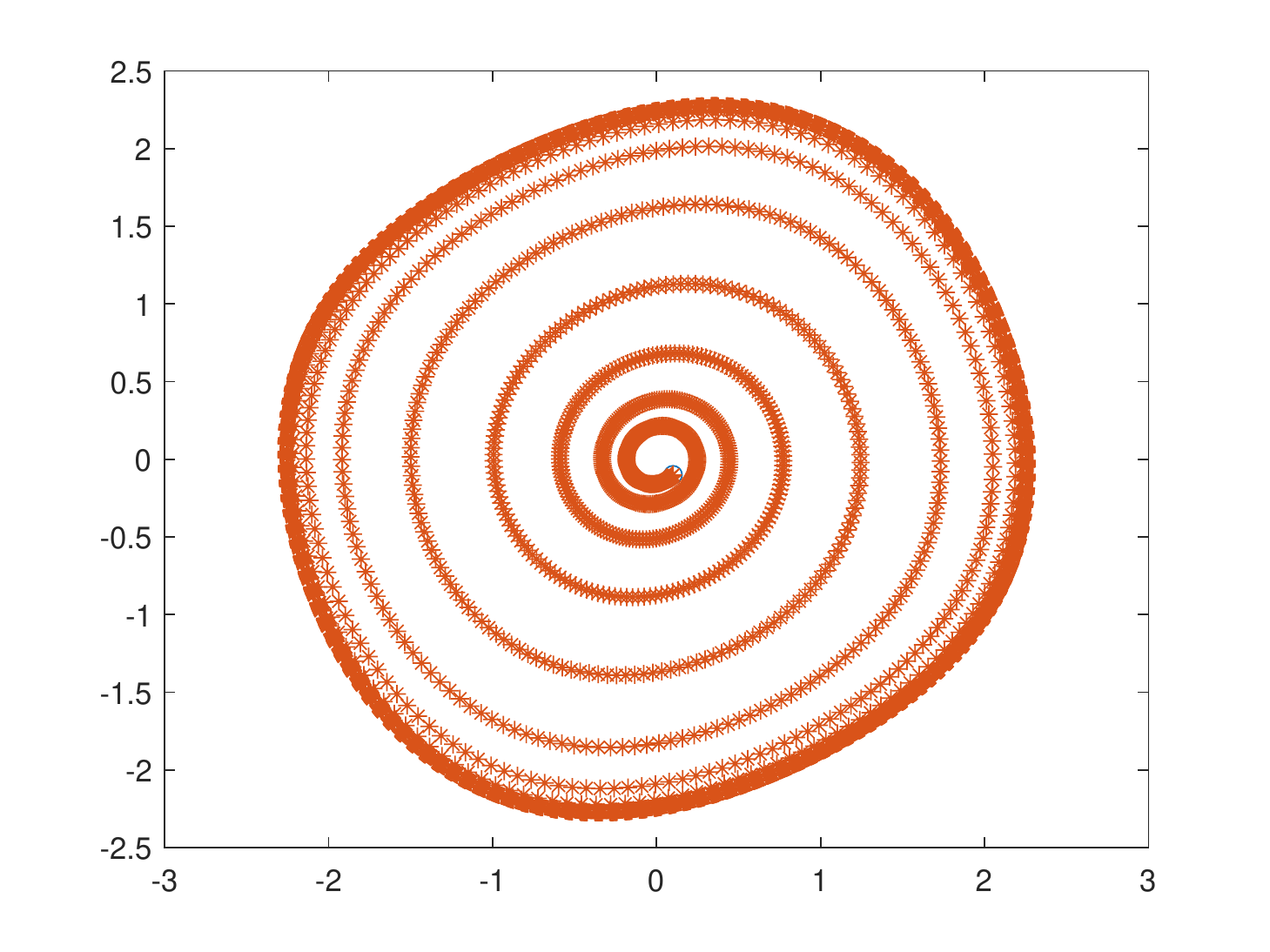}
\caption{Test 2: Cycle limit for Van der Pol oscillator with initial point (-3,3) (left) and with initial point (0.1,-0.1) (right)}
\label{fig:vdp}
\end{figure}

In the second test case we consider the Van der Pol oscillator. The dynamics in \eqref{eq} is given by
\begin{equation}\label{eq:vdp}
	f(x,u)=
		  \begin{pmatrix}
         x_2 \\
         \omega(1-x_1^2)x_2-x_1+u
          \end{pmatrix}\quad u\in U\equiv [-1, 1].
\end{equation}
We note that the origin is a repulsive point for the uncontrolled dynamics in \eqref{eq:vdp}, e.g. $u=0$, if $\omega \in (0,2]$. For this example we consider $\omega=0.15$ in \eqref{eq:vdp}. It is well-known that Van der Pol oscillator is characterized by its cycle limit as shown in Figure \ref{fig:vdp} with two different initial conditions.
%
%
			
In this example we want to minimize the following cost functional:
\begin{equation}\label{cost:vdp}
J_{x,t}(u)= \int_t^T \left( \delta_1 \|y(s)\|_2^2 + \gamma |u(s)|^2 \right) \, ds + \delta_2  \|y(T)\|_2^2,
\end{equation}
where $\delta_1,\delta_2,\gamma$ are positive constants.

\paragraph{Case 1}

We consider the minimization of the terminal cost in \eqref{cost:vdp}, e.g. $\delta_1=\gamma=0$ and $\delta_2=1$. 
Let us consider $x=(-1,1)$, $\Delta t=0.05$ and $T=1$.
The error is computed with respect to the classical approach with a fine grid ($\Delta t=\Delta x=0.002$).
 
\begin{figure}[htbp]
\centering
       \includegraphics[scale=0.25]{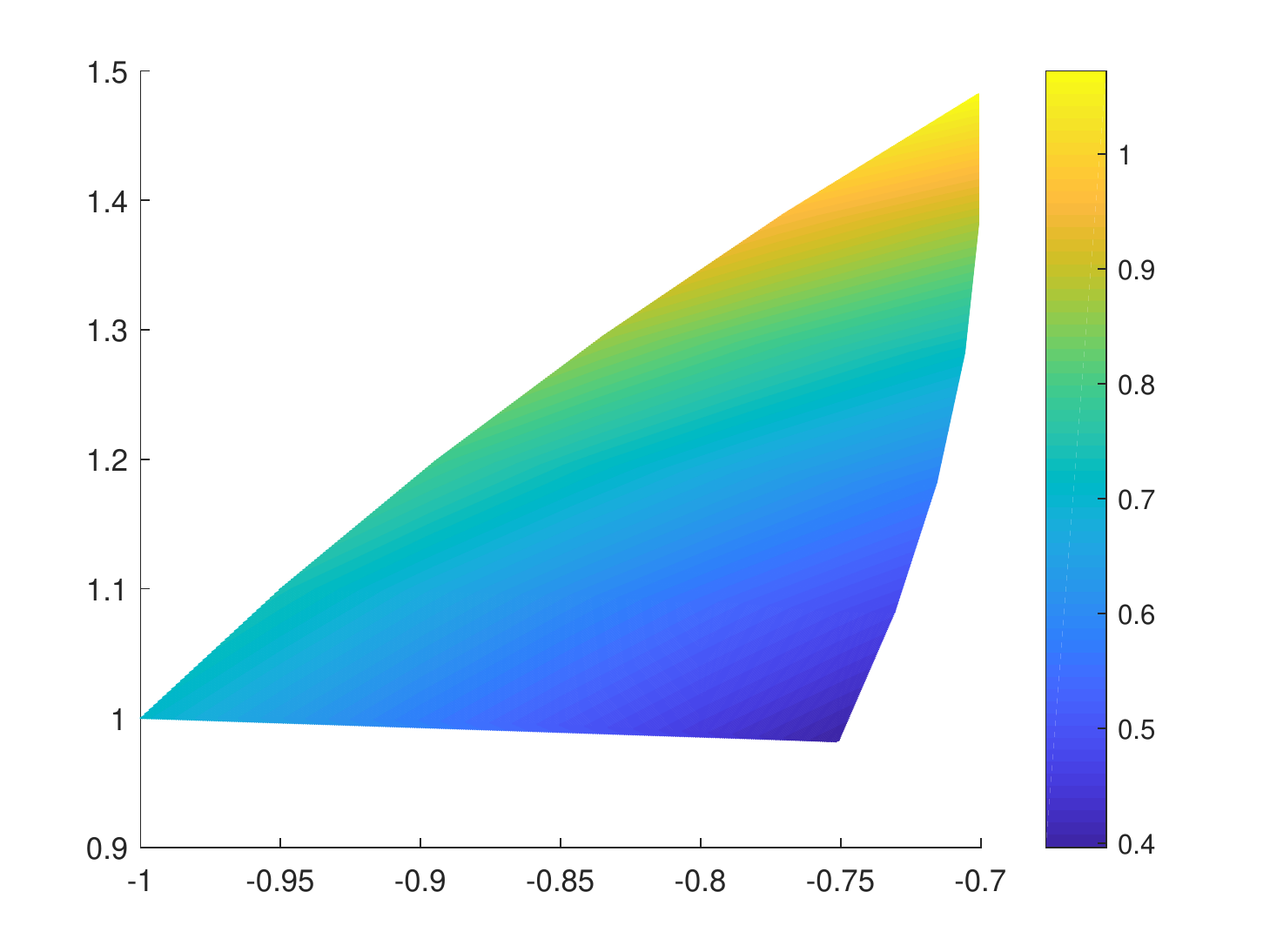} 
         \includegraphics[scale=0.25]{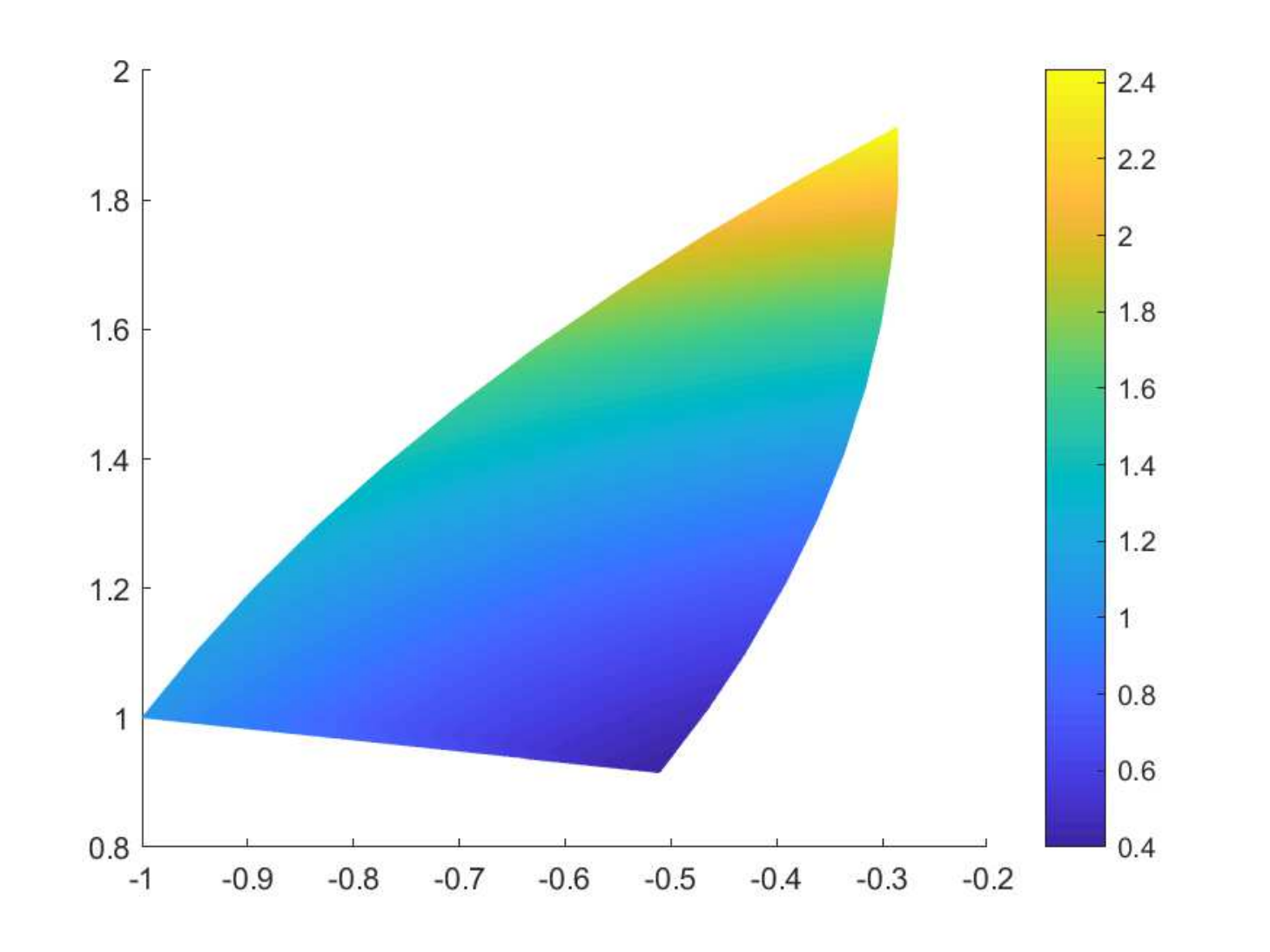}
         \includegraphics[scale=0.25]{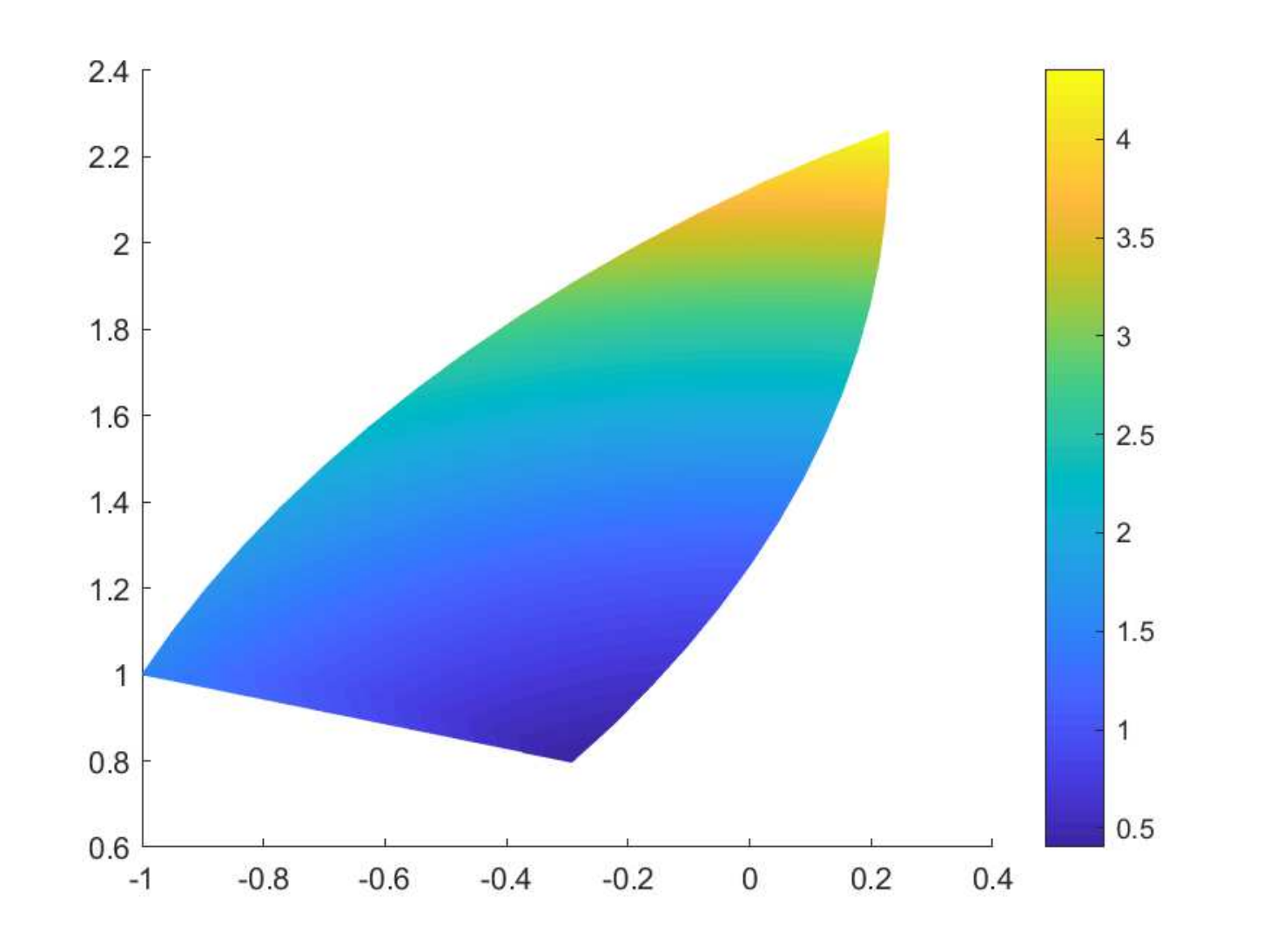}
       \includegraphics[scale=0.25]{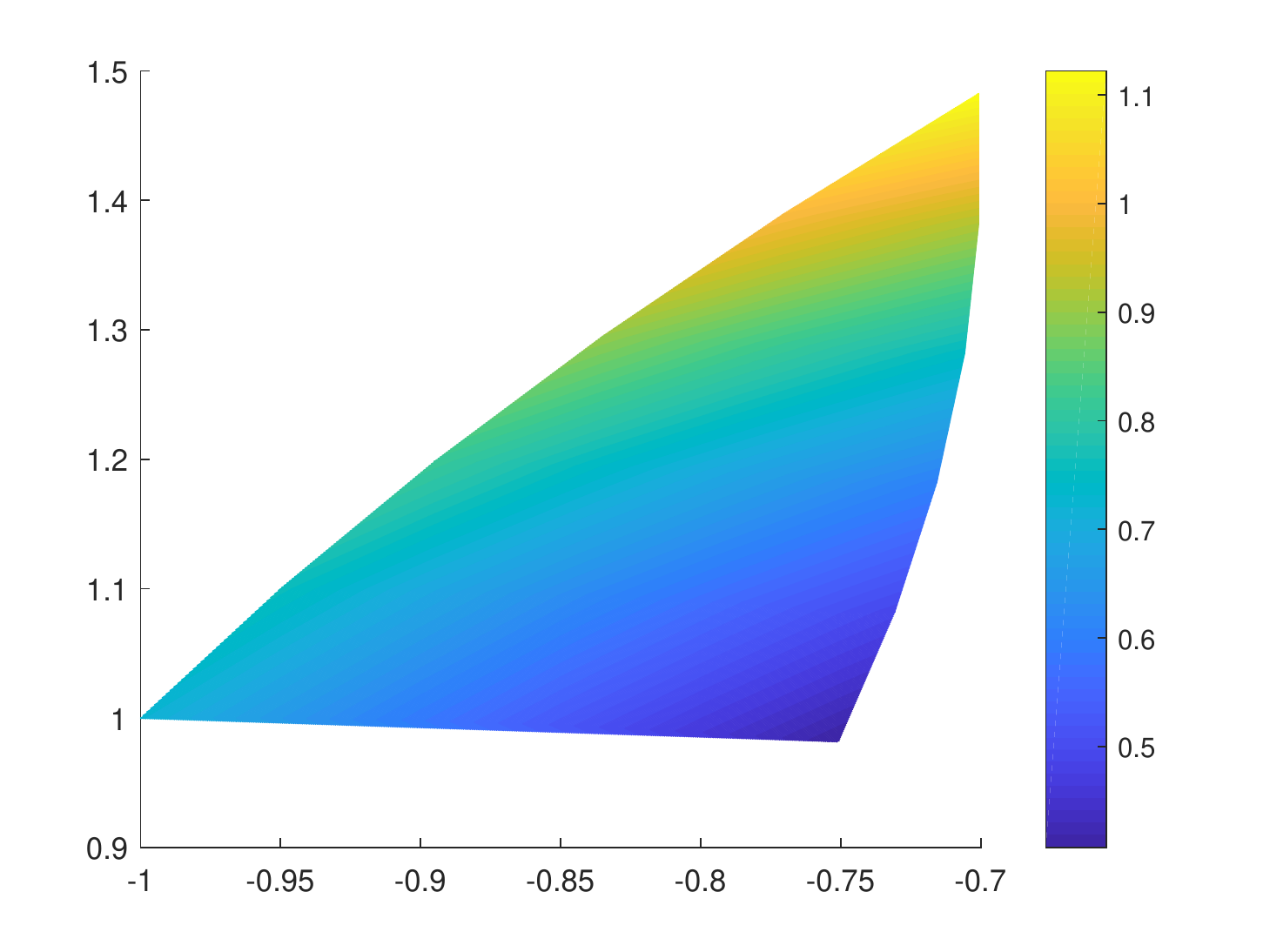} 
         \includegraphics[scale=0.25]{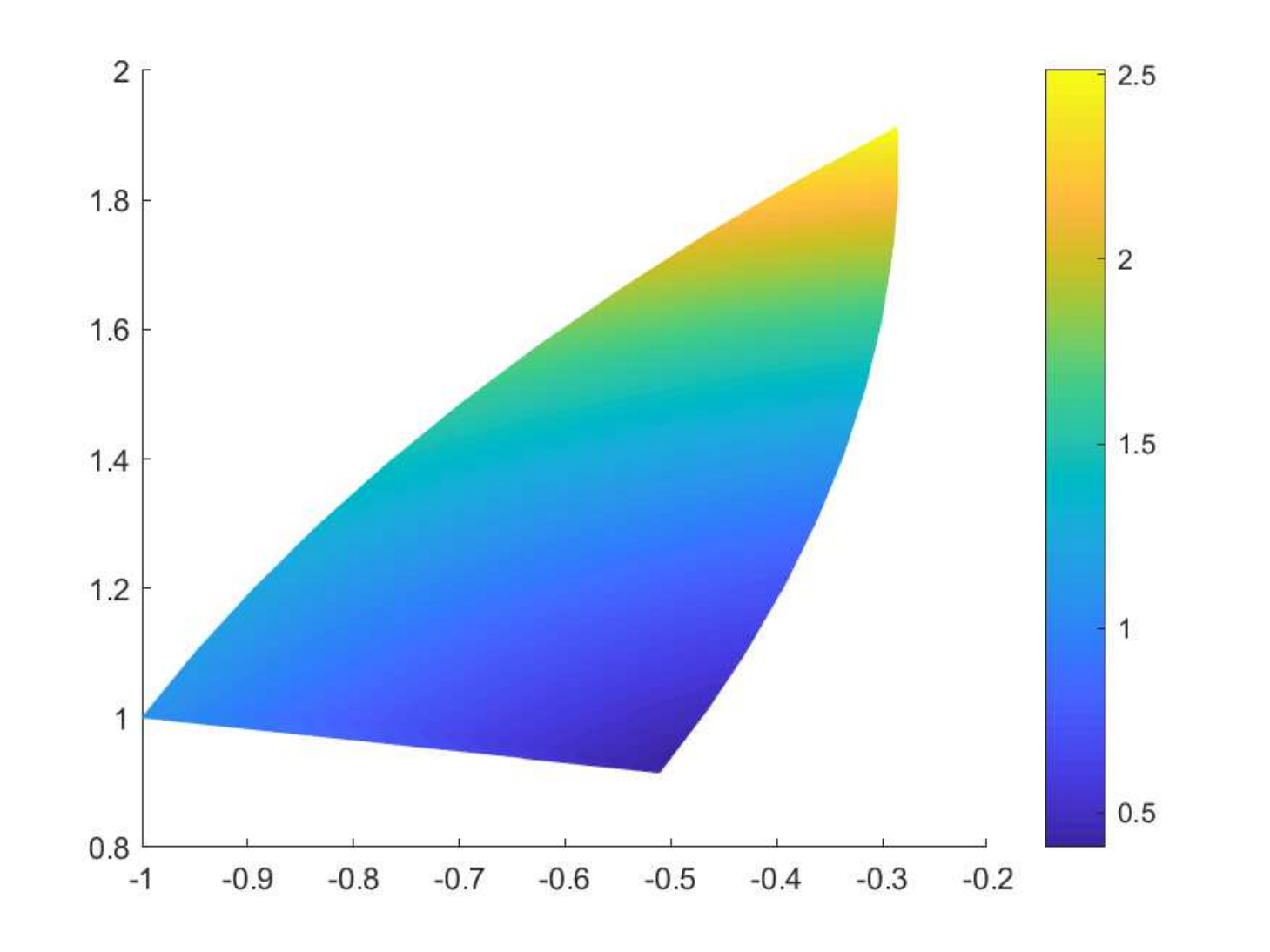}
         \includegraphics[scale=0.25]{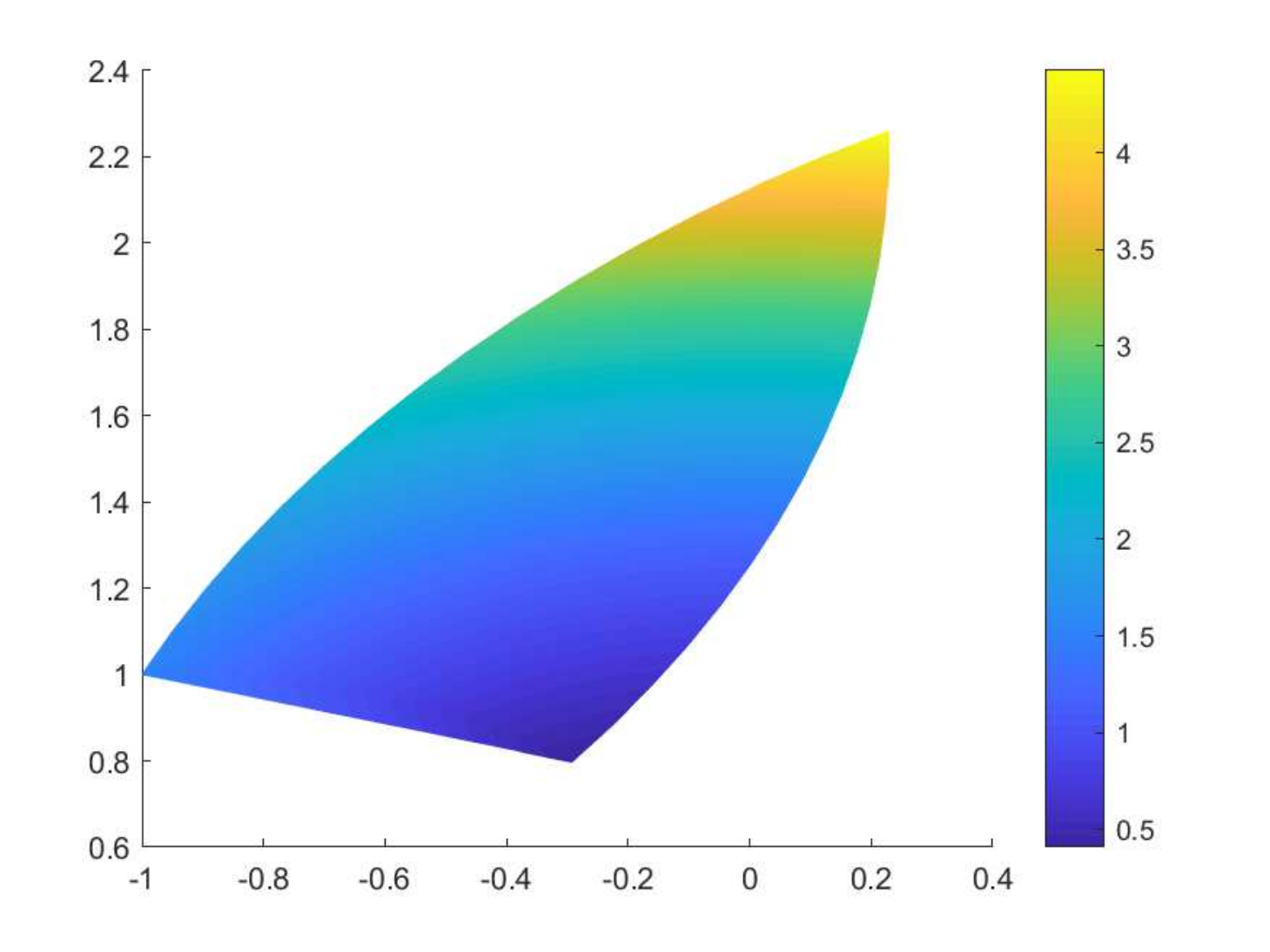}
       \caption{Test 2:  Value function with the classical approach (top) on tree nodes at time $t=0.25$ (left), $t=0.5$ (middle) and $t=0.75$ (right). Value function with the TSA (bottom) on tree nodes at time $t=0.25$ (left), $t=0.5$ (middle) and $t=0.75$ (right)}
       \label{test2:con}
	\end{figure}

We will consider Euler scheme with $U= \{-1,1 \}$ and the tolerance is set equal to $\ep=\Delta t^2$ with $|\mathcal{T}| = 37030$.
In Figure \ref{test2:con} we compare the contour lines of the value function computed by the classical approach with a fine grid and the TSA. We note the approximations show the same behaviour. Furthermore, we mention that the contour line of the value functions are obtained by using MATLAB function {\tt tricontour}, based on a Delaunay's triangulation of the scattered data. We remark that we can compute the value function $V^n(\zeta)$ for $\zeta\in\cup_{k=0}^n \mathcal{T}^k$ since the dynamics is autonomous.

\begin{figure}[htbp]
	\centering
	  \includegraphics[scale=0.4]{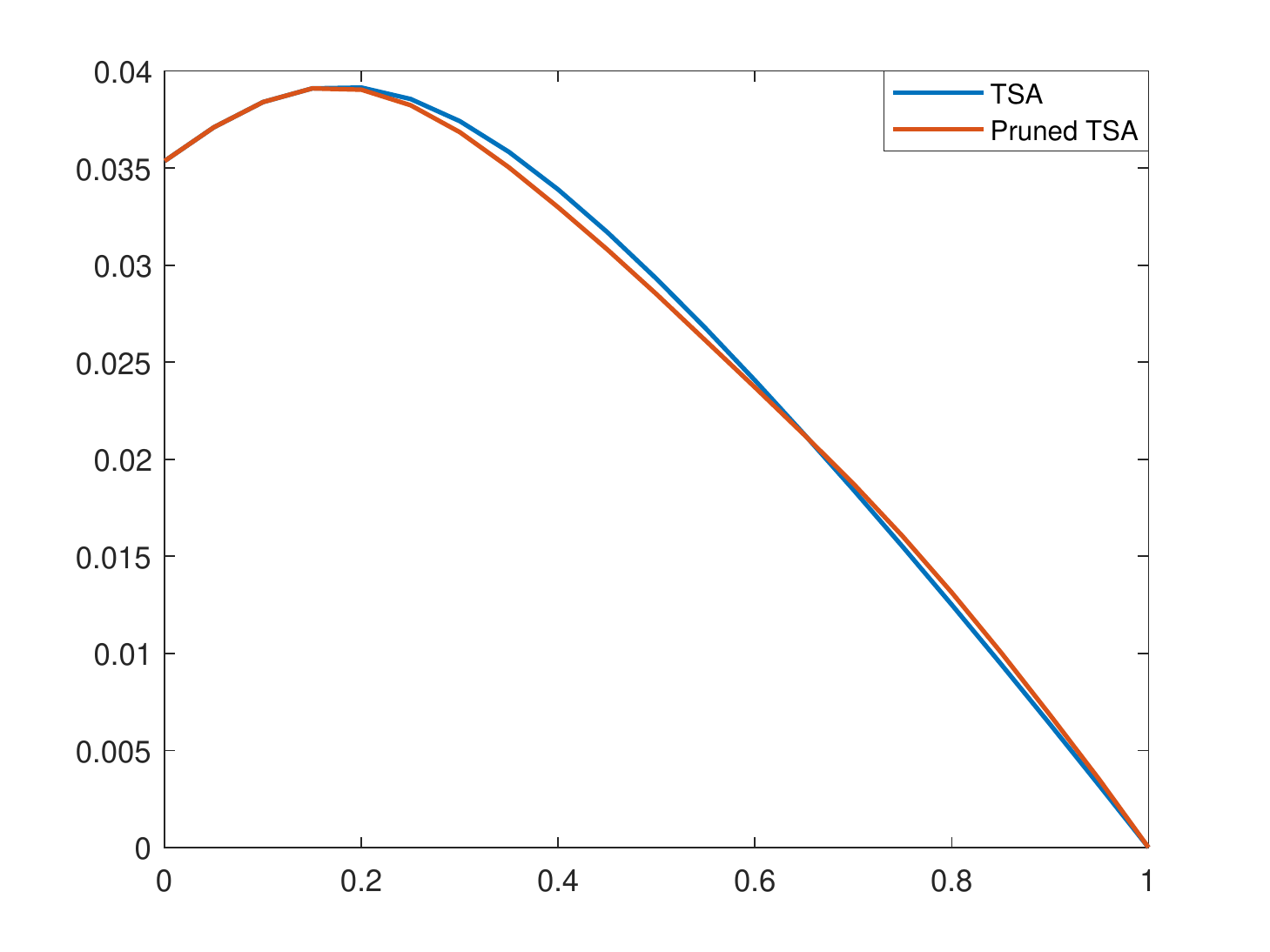} 
	   \includegraphics[scale=0.4]{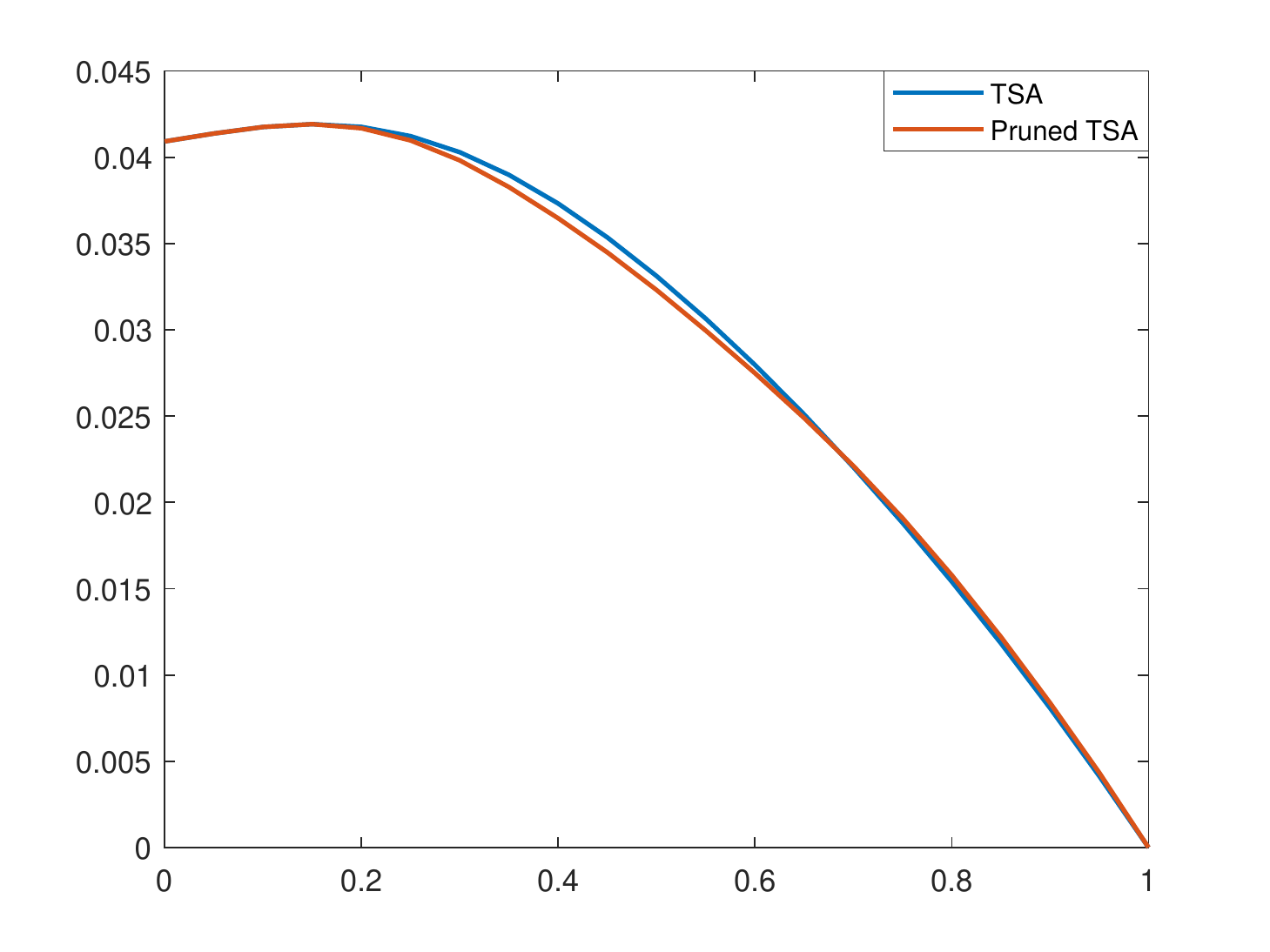} 
       \caption{Test 2: Error in time with TSA without pruning and with pruning with tolerance $\ep=\Delta t^2$ for Case 1 (left) and Case 2 (right) with respect to a value function computed with the classical approach with a very fine grid.}
  \label{fig2:err}
	\end{figure}

The quality of the numerical approximation is confirmed by the error shown in Figure \ref{fig2:err}. As we can see, pruning the nodes does not influence the error. For each time step the error is below to $0.05$ which leads to an accurate approximation of the value function.

\begin{figure}[htbp]
\centering
       \includegraphics[scale=0.25]{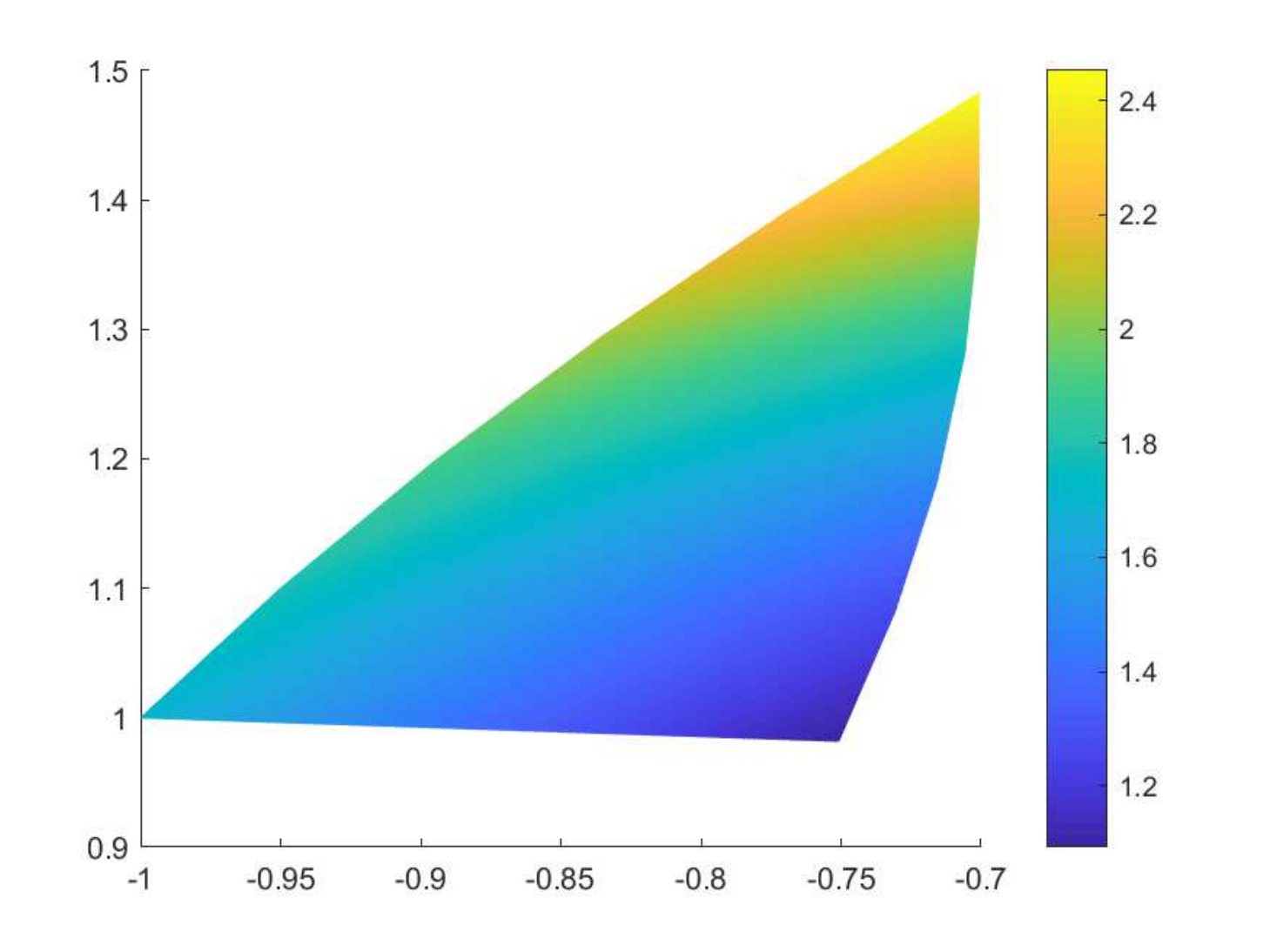} 
         \includegraphics[scale=0.25]{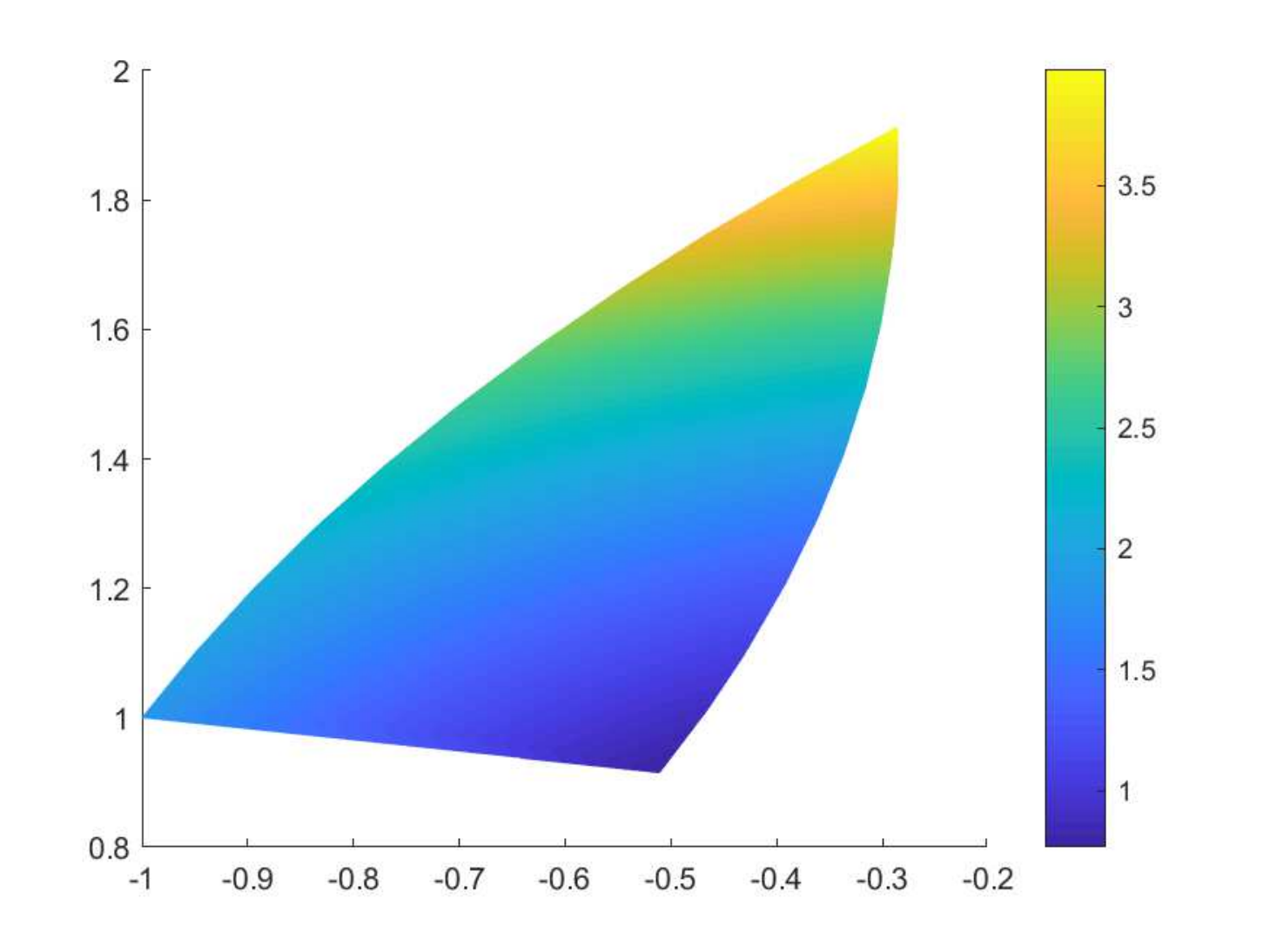}
         \includegraphics[scale=0.25]{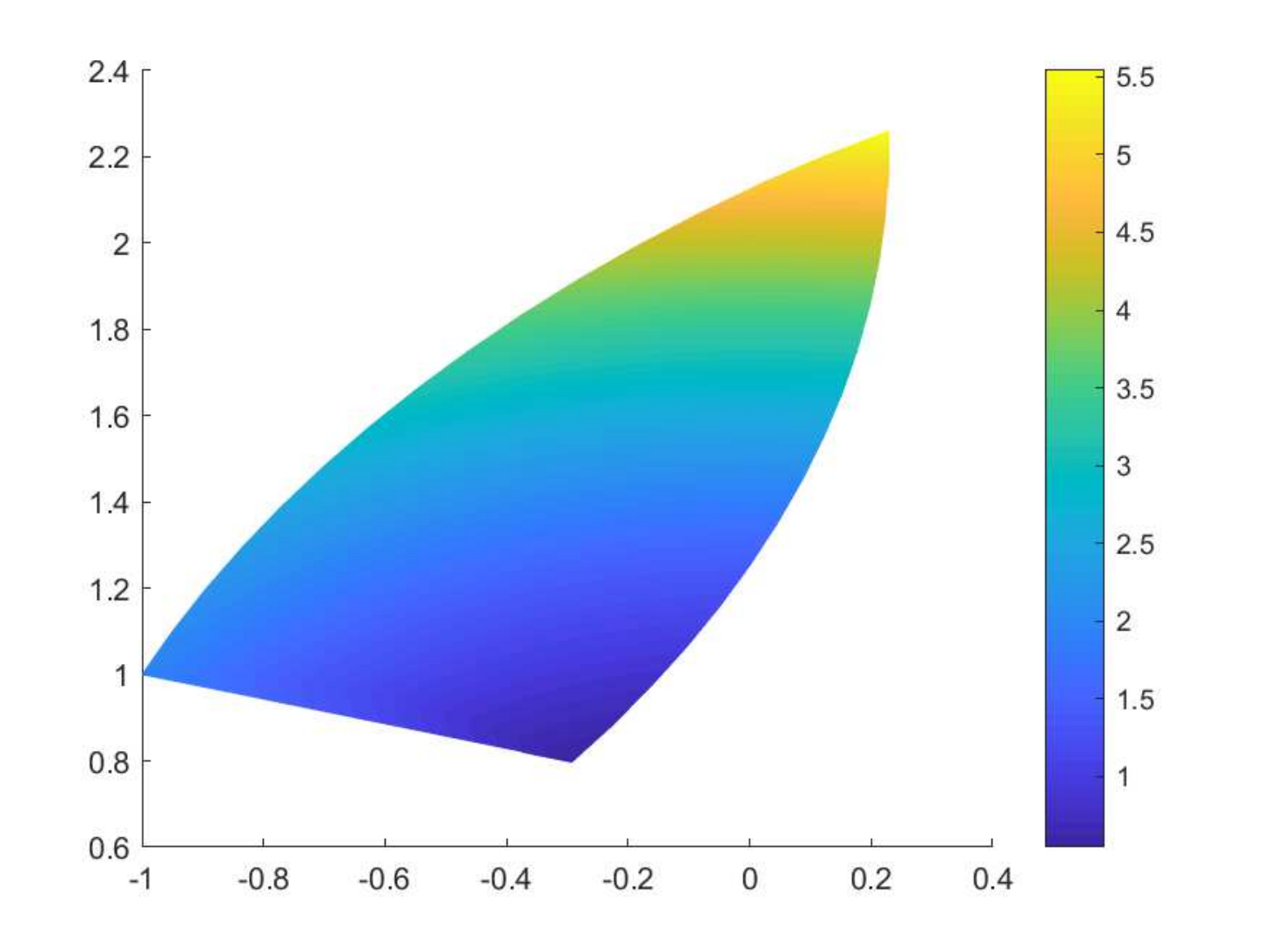}
       \includegraphics[scale=0.25]{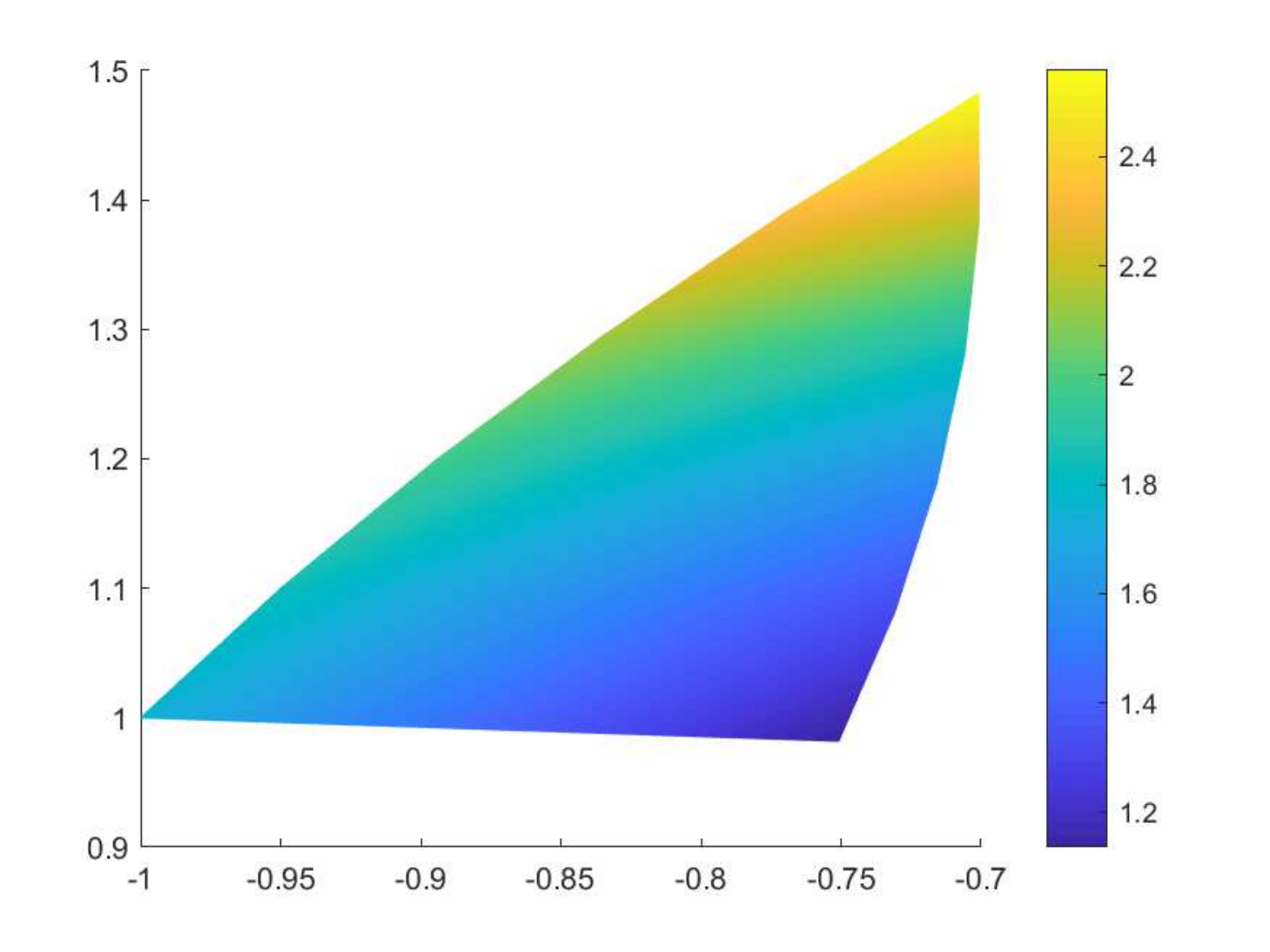} 
         \includegraphics[scale=0.25]{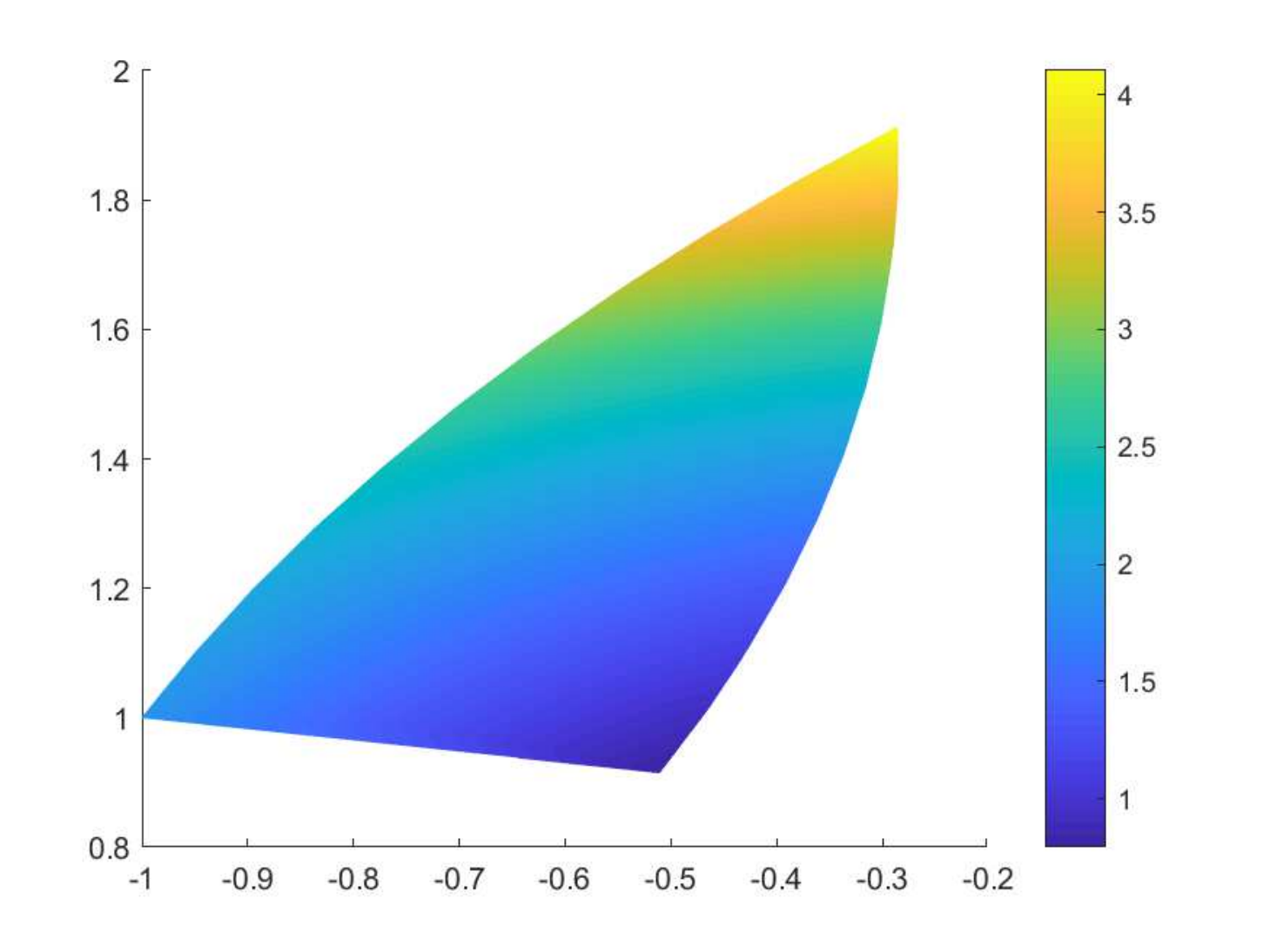}
         \includegraphics[scale=0.25]{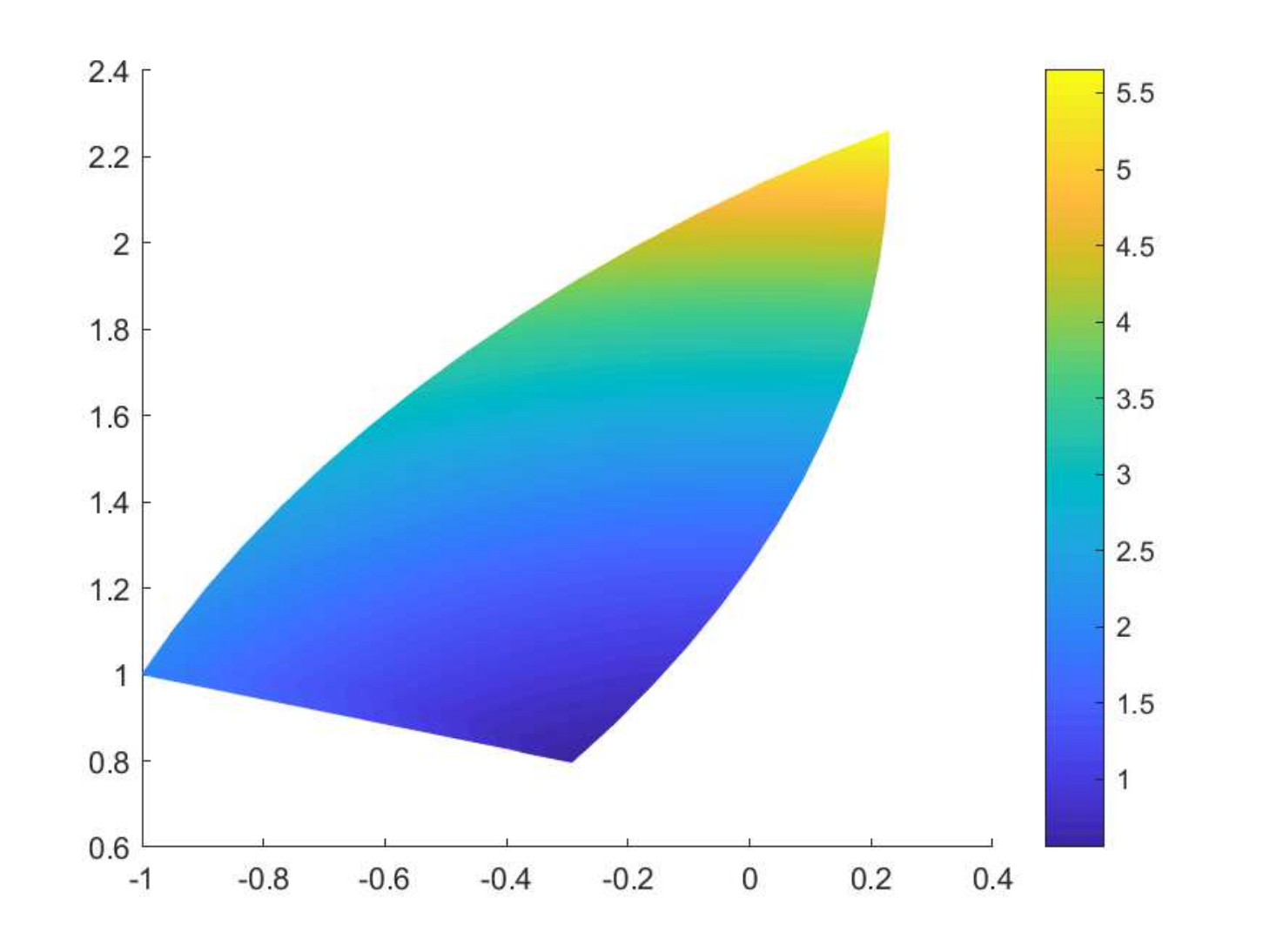}
       \caption{Test 2: Value function with the classical approach (top) on tree nodes at time $t=0.25$ (left), $t=0.5$ (middle) and $t=0.75$ (right).  Value function with the TSA (bottom) on tree nodes at time $t=0.25$ (left), $t=0.5$ (middle) and $t=0.75$ (right)}
       \label{fig2:con2}
	\end{figure}

\paragraph{Case 2}
We consider the minimization of the cost functional in \eqref{cost:vdp} with $\delta_1=\delta_2=1$ and $\gamma=0.01$. Furthermore we set the same initial condition, discretization step and tolerance as in the previous case. The contour lines of the value function are shown in Figure \ref{fig2:con2}.

We note that the results are very similar to the previous case. Our approach is robust with respect to different cost functionals and initial conditions. The right panel of Figure \ref{fig2:err} shows the error for each time step considering the tree algorithm with and without nodal selection.

\paragraph{\revfirst{Case 3}}

\begin{figure}[htbp]
\centering
       \includegraphics[scale=0.4]{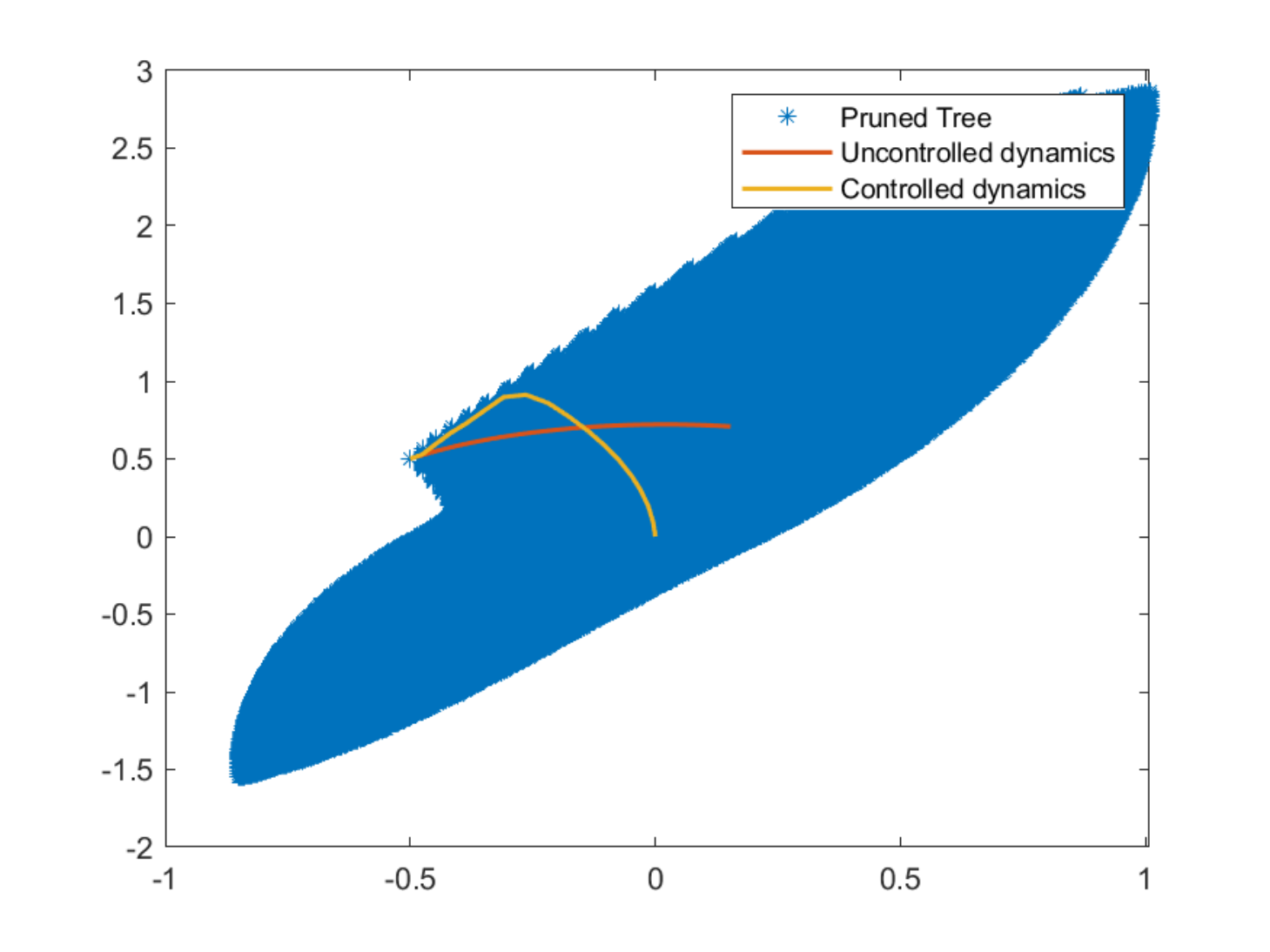} 
         \includegraphics[scale=0.4]{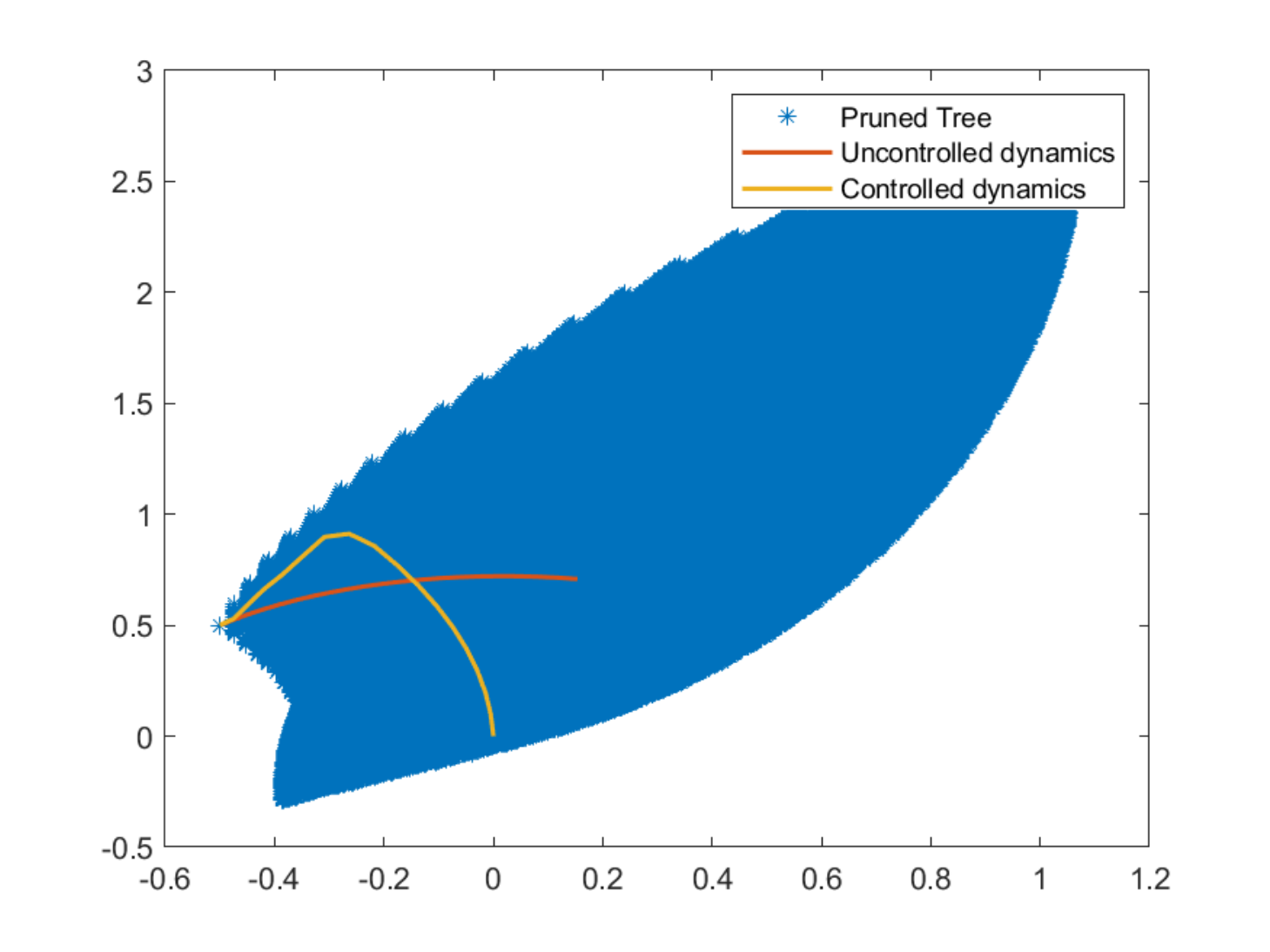}
       \caption{\revfirst{Test 2: Pruned tree with the uncontrolled and controlled dynamics with $U=[-2,0]$ (left) and with $U=[-1,1]$ (right)}}
       \label{fig2:con3}
	\end{figure}

\revfirst{ In the last case we deal with a two dimensional control space, considering the parameter $\omega$ in \eqref{eq:vdp} as a control, e.g. $\omega\in U$. Therefore, we consider as control variables $(\omega, u)\in U\times U$ in \eqref{eq:vdp}. 
 In the cost functional \eqref{cost:vdp} we consider again $\delta_1=\gamma=0.1$ and $\delta_2=1$, with $x=(-0.5,0.5)$, $\Delta t = 0.05$ and  $T=1$. We consider two different choices for the control set: $U=[-2,0]$ and $U=[-1,1]$. The control set is discretized with step-size $\Delta u =0.2$, obtaining altogether $100$ discrete controls for both examples. In Figure \ref{fig2:con3} we show the results in both situations. We can observe that the tree has a different shape due to the different control space. Here, we have set the pruning criteria with $\ep = \Delta t^2$. Finally, we note that in both situations we are able to steer the solution to the origin.}
\subsection{Test 3: Damped harmonic oscillator with sinusoidal driving force}
In this third example we consider a non-autonomous dynamical system: a damped oscillator driven by a sinusoidal external force. The dynamics in \eqref{eq} is given by
\begin{equation}\label{eq:dho}
	f(x,u,t)=
		  \begin{pmatrix}
         x_2 \\
         -\omega x_2-\omega^2 x_1+\sin(\omega t) +u
          \end{pmatrix}\quad u\in U\equiv [-1, 1].
\end{equation}
for $x = (x_1,x_2)\in\R^2$. In this example, we aim to show that our approach works also with non-autonomous dynamics. In this case we can not compute the value function $V^n(\zeta)$ on the sub-tree $\cup_{k=0}^n \mathcal{T}^k$, but only at the $n-$th time level $\mathcal{T}^n$ \new{and we will apply the pruning rule \eqref{tol_cri}}. The uncontrolled dynamics ($e.g.\, u=0$) converges asymptotically to the cycle limit:
$$
\overline{x}_1(t)=\frac{1}{\omega^2}\sin(\omega t + \pi/2), \; \overline{x}_2(t)=\frac{1}{\omega}\cos(\omega t + \pi/2) \;.
$$

We used the same cost functional of the previous case with $\delta_1=\gamma=0.1$, $\delta_2=1$. The parameters are set as follows: $\omega=\pi/2,\, x=(-0.5,0.5),\, U=\{-1,0,1 \}, \Delta t=0.05,\, T=1,\, \ep=\Delta t^2$. The cardinality of tree in this case is $32468$.
\begin{figure}[htbp]
	\centering
	  \includegraphics[scale=0.4]{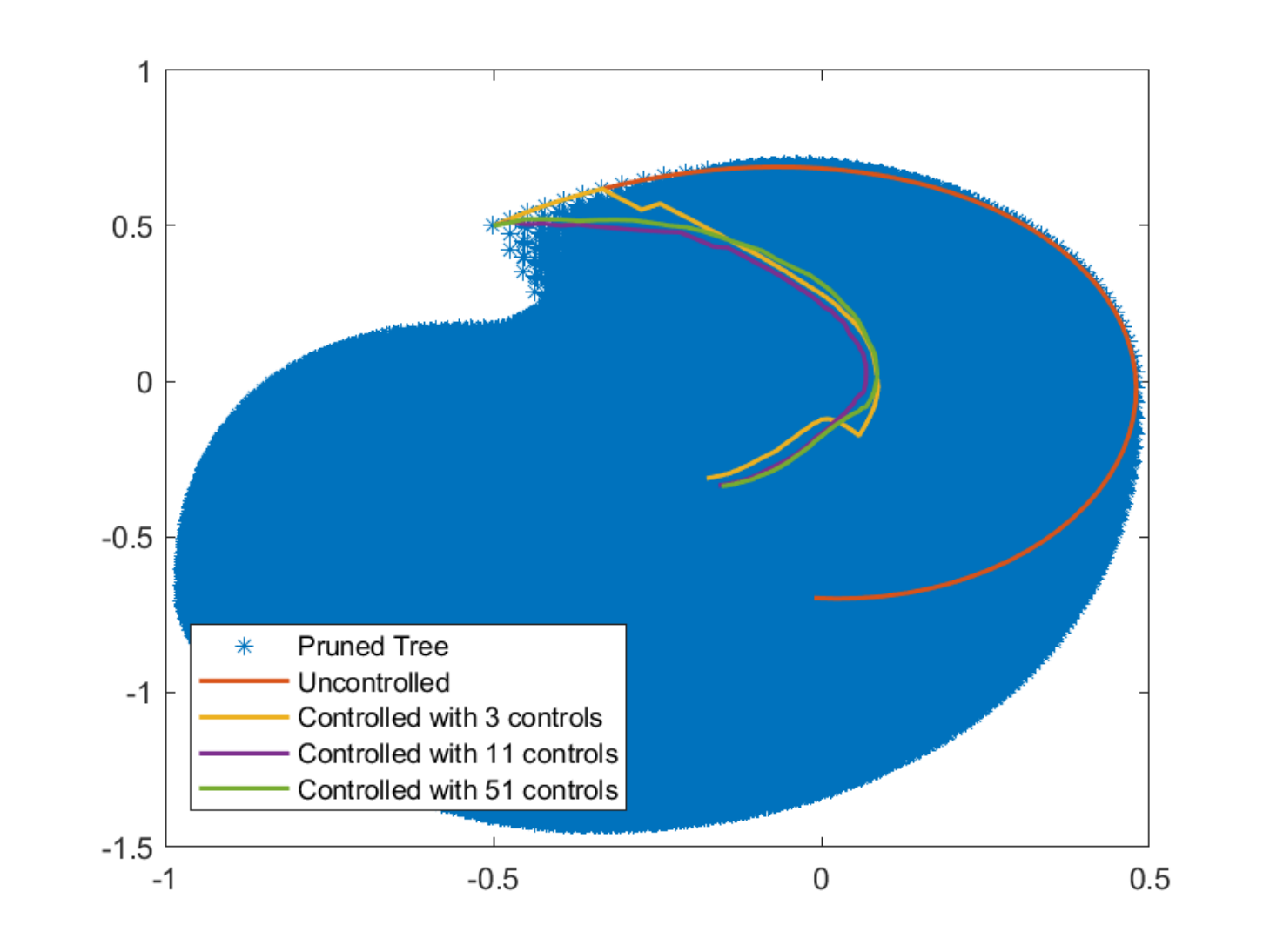} 
	  \includegraphics[scale=0.4]{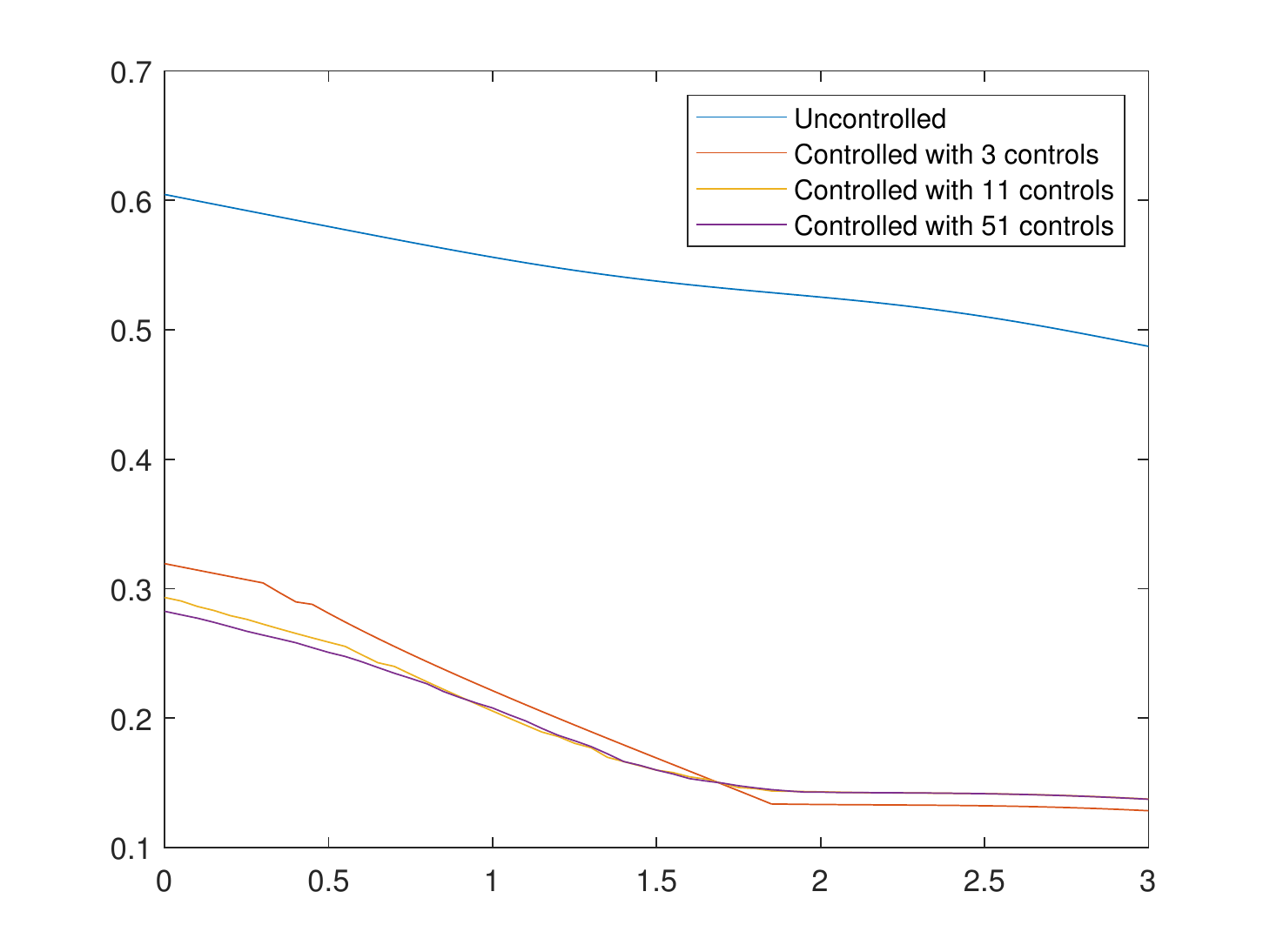}
       \caption{\revfirst{Test 3: Pruned tree with the uncontrolled and controlled dynamics (left) and comparison of the cost functional on time varying the number of discrete controls (right)}}
       \label{fig:nonauto}
\end{figure}
In the left panel of Figure \ref{fig:nonauto} we show the tree nodes and the optimal trajectory computed with Algorithm \ref{alg_1} and the uncontrolled solution. To show the quality of the controlled solution we evaluate the cost functional for each time step as shown in the right panel of Figure \ref{fig:nonauto}. As expected the controlled trajectory is always below the uncontrolled one. \revfirst{In order to further show the effectiveness of the pruning criteria we have increased the number of controls up to $M=51$ and the horizon up to $T=3$. Again, this would not be possible without a pruning criteria due to the dimension of the tree}. 
\subsection{Test 4: Heat equation}
The fourth example concerns the control of a PDE. In the first three examples \revfirst{we showed the accuracy of our method} with respect to existing methods \revfirst{ for low-dimensional problems}. In what follows we would like to give an idea of how the proposed method can work in higher dimension.

We want to study the following heat equation:
\begin{equation}
\begin{cases}
y_t=\sigma y_{xx}+ y_0(x) u(t) & (x,t) \in \Omega \times [0,T] \;, \\
 y(x,t)=0 & (x,t) \in \partial \Omega \times [0,T] \;, \\
 y(x,0)=y_0(x) & x \in \Omega \;,
\end{cases}
\label{heat}
\end{equation}
where the state lies in an infinite-dimensional Hilbert space (see e.g. \cite{E10}). \revfirst{Here, we consider the term $y_0(x)u(t)$ to provide a spatial dependence to the control input. This is a particular choice, but the algorithm has no restrictions on more general shape functions.}
To write equation \eqref{heat} in the form \eqref{eq} we use the centered finite difference method which leads to the following ODEs system
\begin{equation}
\dot{y}(t)=A y(t) + B u(t),
\label{ODE}
\end{equation}
where the matrix $A\in\R^{d\times d}$ is the so called {\em stiffness} matrix whereas the vector $B\in\R^n$ is given by  $(B)_i=y_0(x_i)$ for $i=1,\ldots,n$ and $x_i$ is the spatial grid with constant step size $\Delta x$.
The cost functional we want to minimize reads: 
\begin{equation*}
J_{y_0,t}(u)= \int_{t}^T \left( \delta_1 \|y(s)\|_2^2 \,dx+ \gamma |u(s)|^2 \right)\, ds +  \|y(T)\|_2^2,
\end{equation*}

where $y(t)$ is the solution of \eqref{ODE}, $u(t)$ is taken in the admissible set of controls $\mathcal{U}=\{u:[0,T] \rightarrow [-1,1] \}$ and $\Omega = [0,1]$. \revsecond{We set $\delta_1=1$ and $\gamma=0.01$.}
\paragraph{Smooth initial condition} In the numerical approximation of \eqref{heat} we consider $y_0(x)=-x^2+x$, $\Delta x=10^{-3},\, \Delta t=0.05,\, T=1$ and $\sigma=0.1$. The dimension of the problem is $d=1000$. We use an implicit Euler scheme to integrate the system \eqref{ODE} and guarantee its stability. We note that the use of a one step implicit is straightforward even if we have introduced an explicit scheme in the previous sections. We refer to \cite{QV99} for more details about the method. 
\begin{figure}[htbp]	
\centering
	\includegraphics[scale=0.25]{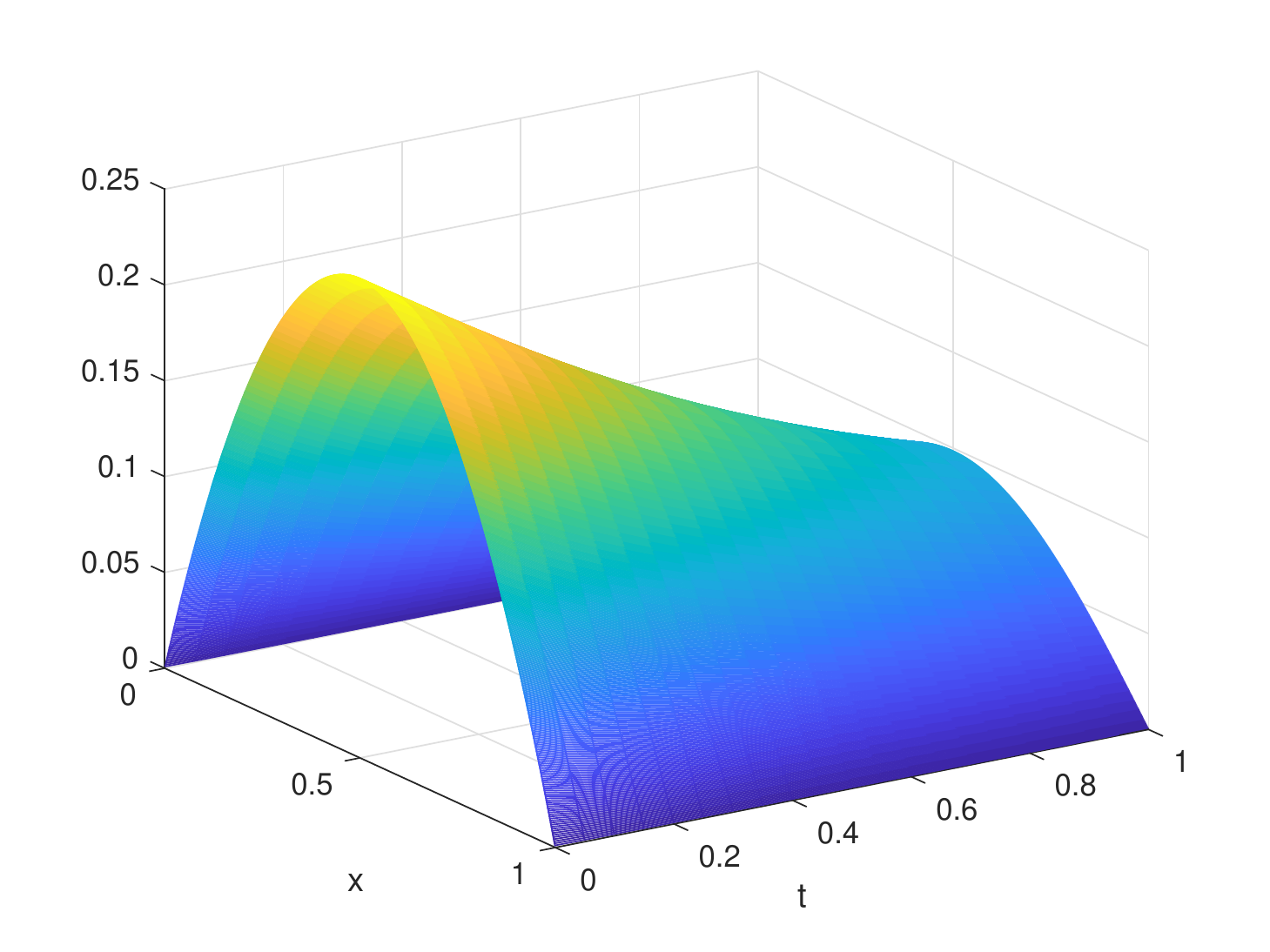}	
	\includegraphics[scale = 0.25]{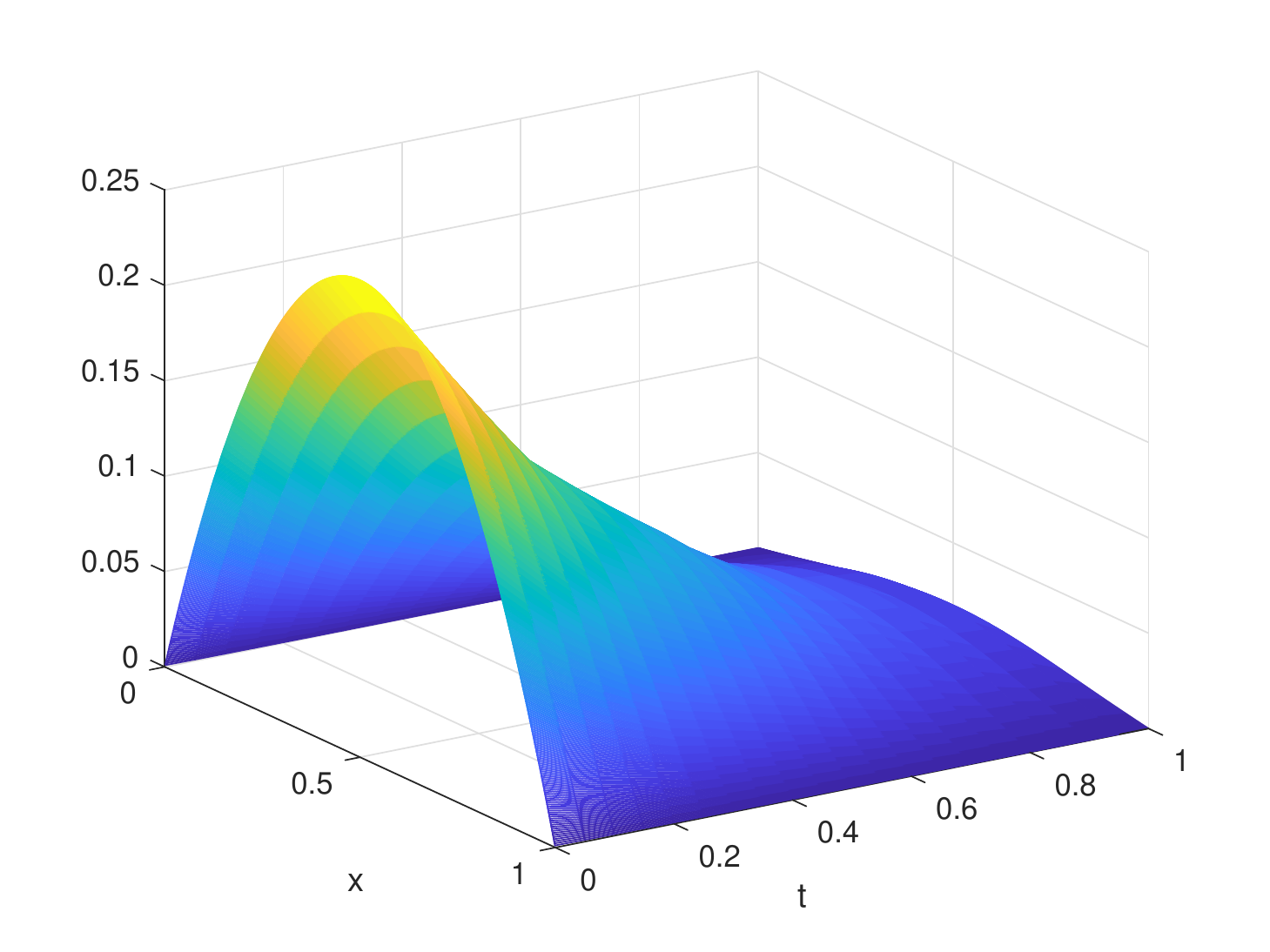}	
	\includegraphics[scale= 0.25]{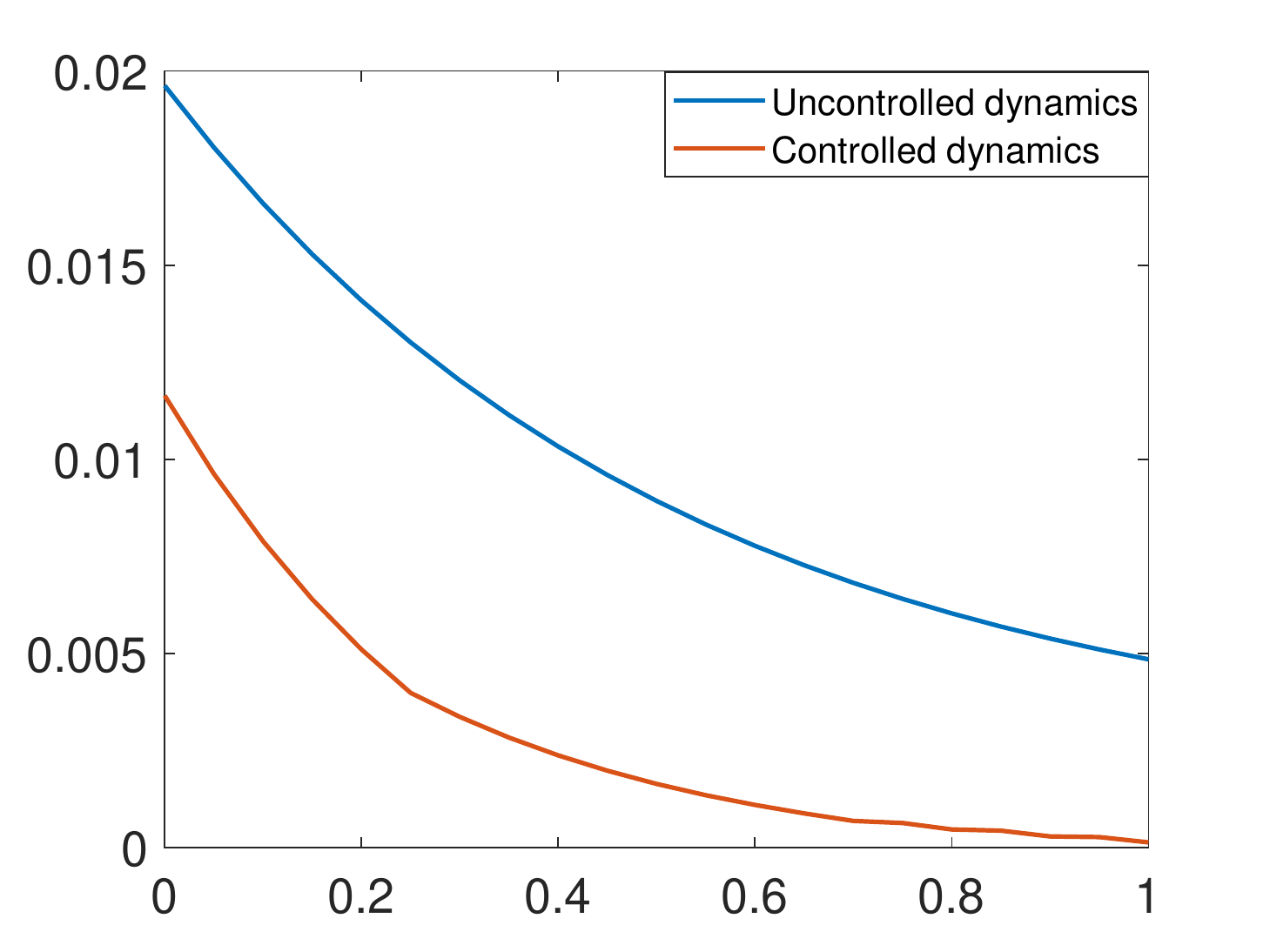}		
\caption{Test 4 (smooth initial condition): Uncontrolled solution (left), optimal control solution (middle),  time comparison of the cost functional of the uncontrolled solution and controlled solution (right).}	
\label{fig3:heat}
\end{figure}

The solution of the uncontrolled problem \eqref{heat} with $u(t)\equiv 0$ is shown in the left panel of Figure \ref{fig3:heat}. In the middle we show the solution of the controlled problem where the value function is computed with Algorithm \ref{alg_1} and the control is computed as explained in \eqref{feedback}. We note that feedback control was computed with the discrete control set $U=\{-1,0,1\}$ as for the value function. A refinement for the control set would require further investigation that we will address in the near future. However, it is extremely interesting to show that we are able to compute the value function for \eqref{heat} in dimension $1000$. This approach might substitute recent advances where the feedback for PDEs was computed by coupling the HJB equation with model order reduction techniques such as, e.g., Proper Orthogonal Decomposition \cite{KVX04}. Finally in the right panel of Figure \ref{fig3:heat} we show the time behaviour of the cost functional for the uncontrolled and the controlled solution. As expected, the cost functional of the latter is lower.

\paragraph{Non-smooth initial condition}

In this example we consider the following non-smooth initial $y_0(x)=\chi_{[0.25,0.75]}(x)$, where $\chi_{\omega}(x)$ is the characteristic function in the domain $\omega$, whereas the other parameters are set as in the previous case.

 \begin{figure}[htbp]	
\centering
	\includegraphics[scale=0.25]{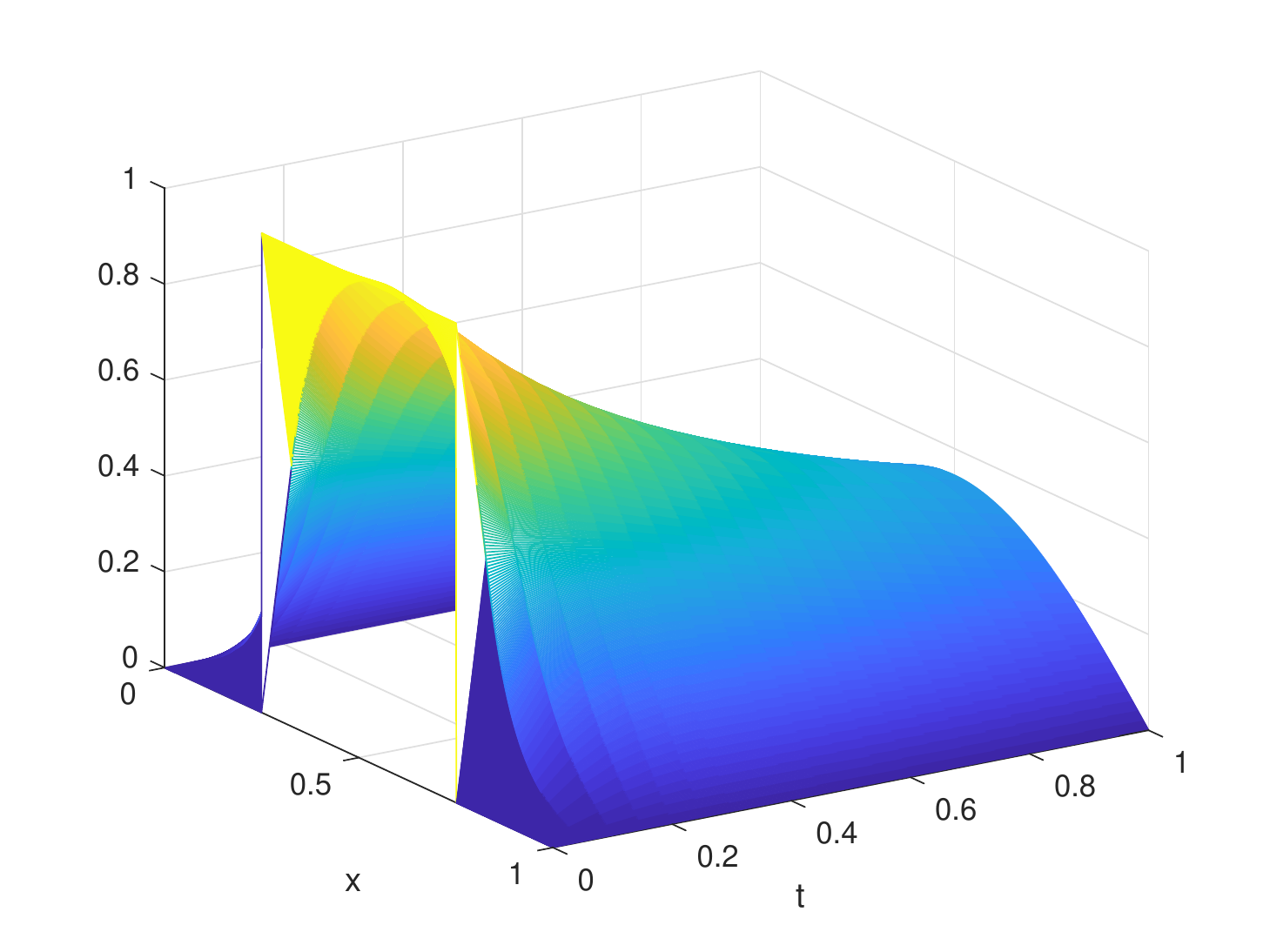}	
	\includegraphics[scale = 0.25]{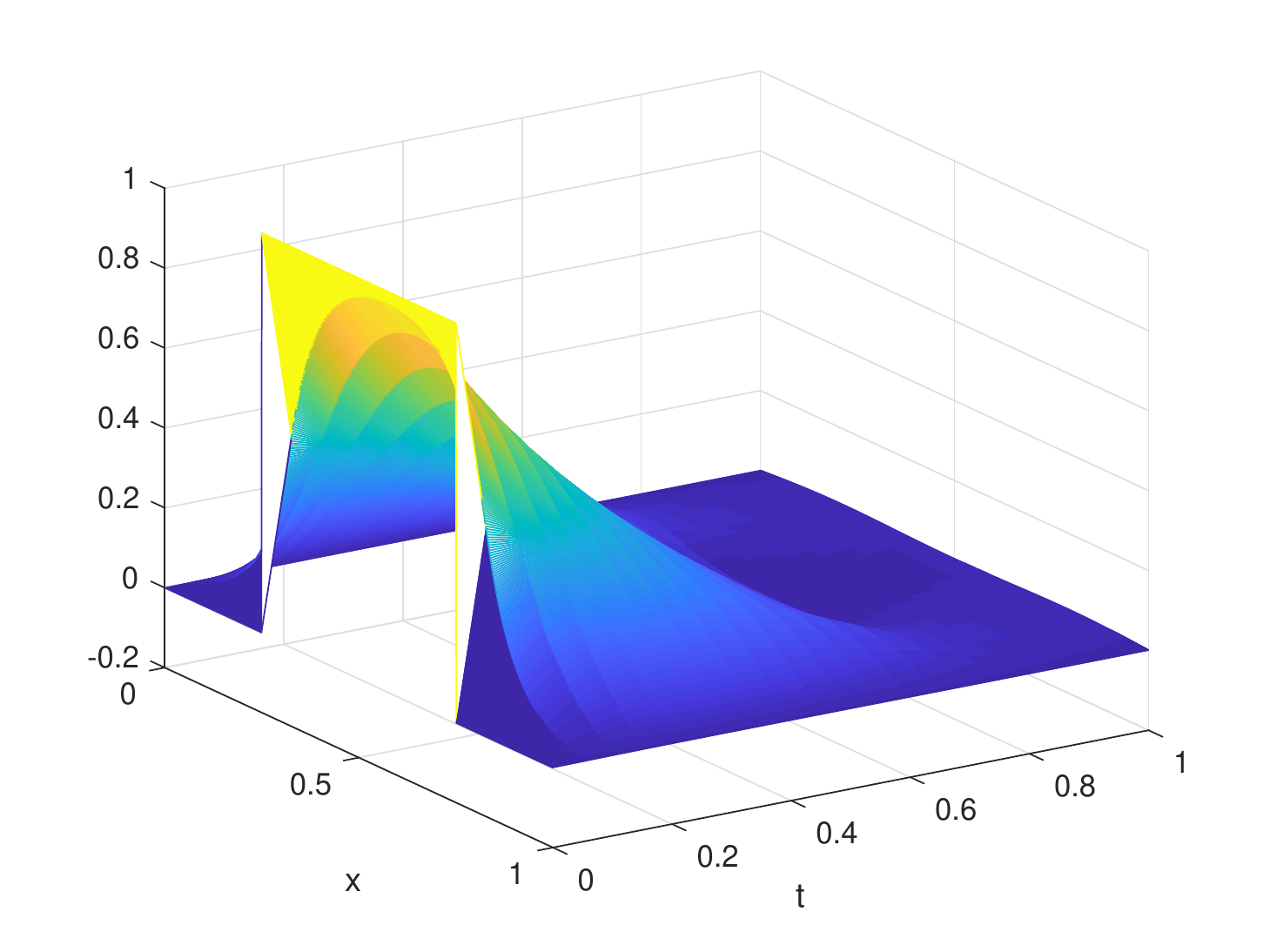}	
	\includegraphics[scale= 0.25]{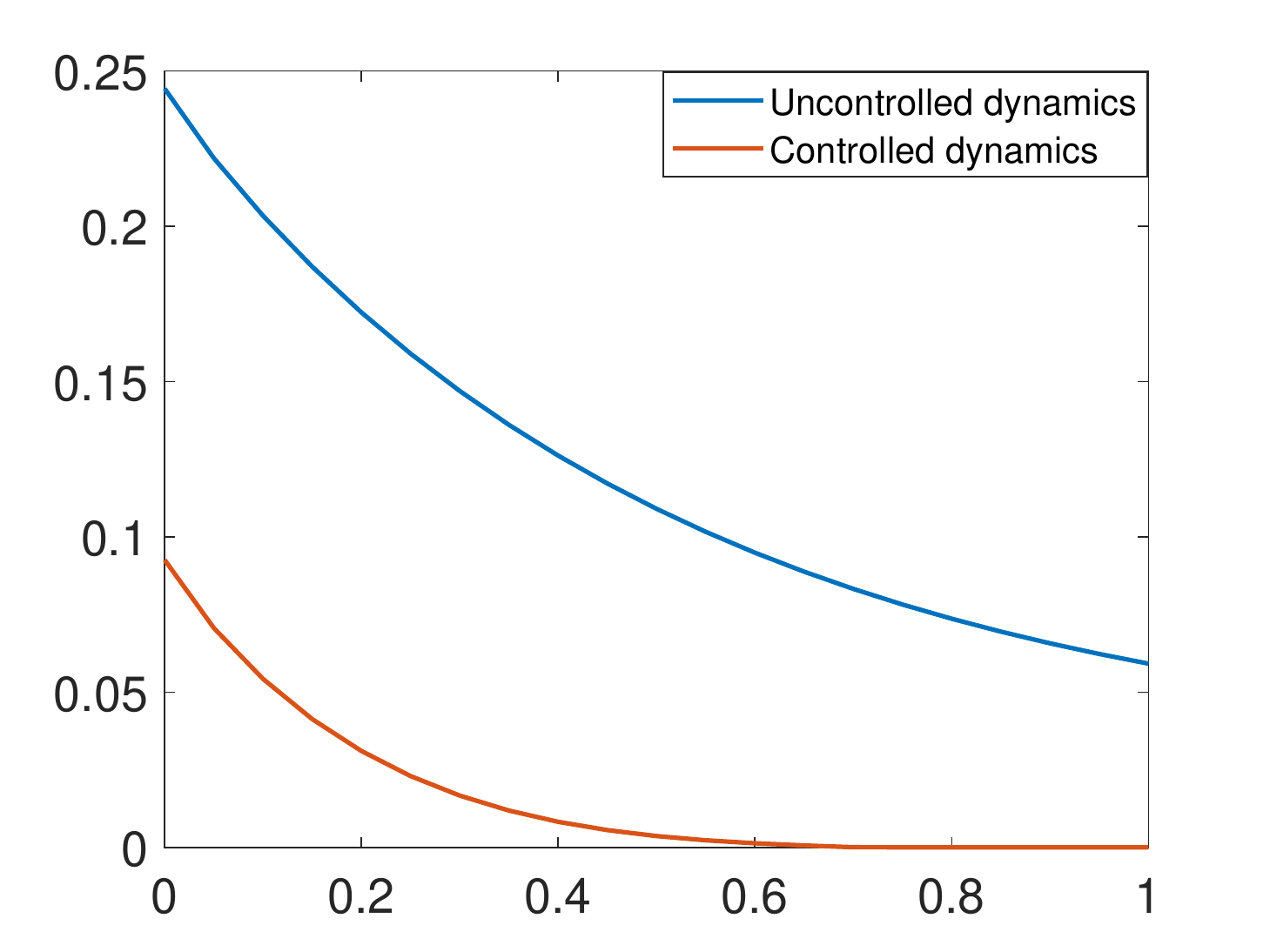}		
\caption{Test 4 (non-smooth initial condition): Uncontrolled solution (left), optimal control solution (middle),  time comparison of the cost functional of the uncontrolled solution and controlled solution (right).}	
\label{fig4:heat}
\end{figure}

As one can see from Figure \ref{fig4:heat}, we are able to approximate the control problem even if the initial condition is non-smooth. We note that, although the simple diffusive properties of the problem, a model reduction approach will not be able to reconstruct such initial condition with a few number of basis functions. Therefore it will not be possible to solve this problem with a classical approach. This again shows the effectiveness of the method.

\subsection{Test 5: Wave equation}
In this last example we consider a hyperbolic PDE, the wave equation which reads:
\begin{equation}
\begin{cases}
w_{tt}= c \, w_{xx}+ \chi_{\omega}(x) u(t) & (x,t) \in \Omega \times [0,T] \;, \\
 w(x,t)=0 & (x,t) \in \partial \Omega \times [0,T] \;, \\
 w(x,0)=w_0(x)\,, \; w_t(x,0)=w_1(x)  & x \in \Omega \;.
\end{cases}
\label{wave}
\end{equation}

where $\omega$ is a subset of $\Omega$. \revfirst{For all initial data $(w_0, w_1) \in H^1_0(\Omega)\times L^2(\Omega)$ and every $u(t) \in L^2(0, T)$, there exists a unique solution $w \in C^0(0, T; H^1_0(\Omega)) \cap C^1(0, T;L^2(\Omega)) \cap C^2(0, T; H^{-1}(\Omega))$ of the Cauchy problem (\ref{wave}).} We refer to \cite{E10} for more details about this equation.
We can rewrite the wave equation in the following compact form
$$
\dot{y}(t)=Ay(t) + Bu(t)\;,
$$
defining
\begin{equation}
	y(t)=
		  \begin{pmatrix}
         w(t) \\
         w_t(t)
          \end{pmatrix}, \quad
           A=
          \begin{pmatrix}
           0 \;\; I \\
           c\, \partial^2_x \;\; 0 
          \end{pmatrix} , \quad
          Bu(t)=
		  \begin{pmatrix}
         0 \\
         \chi_{\omega}(x) u(t)
          \end{pmatrix}.         
\end{equation}

\begin{figure}[htbp]	
\centering
	\includegraphics[scale=0.25]{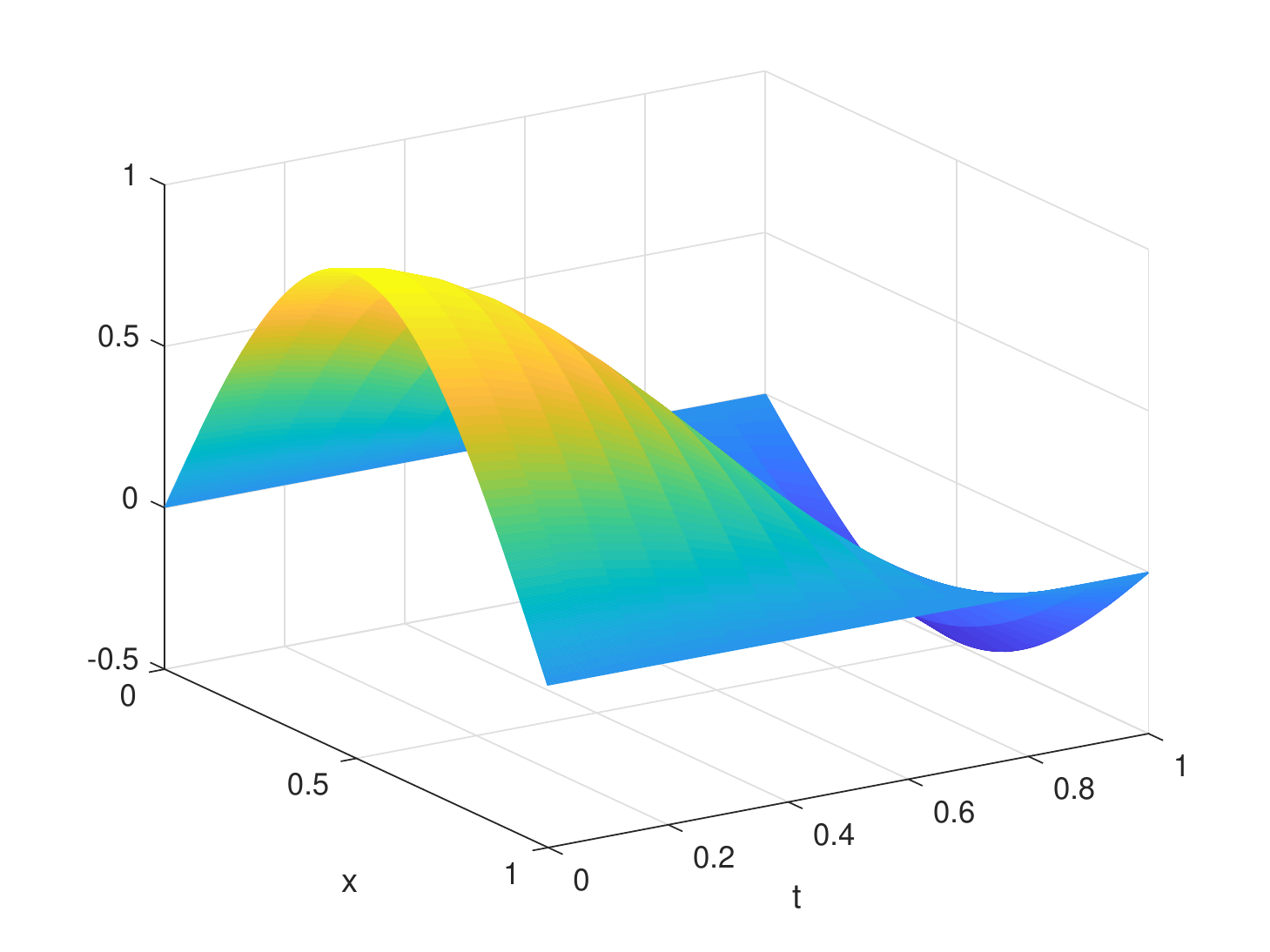}	
	\includegraphics[scale = 0.25]{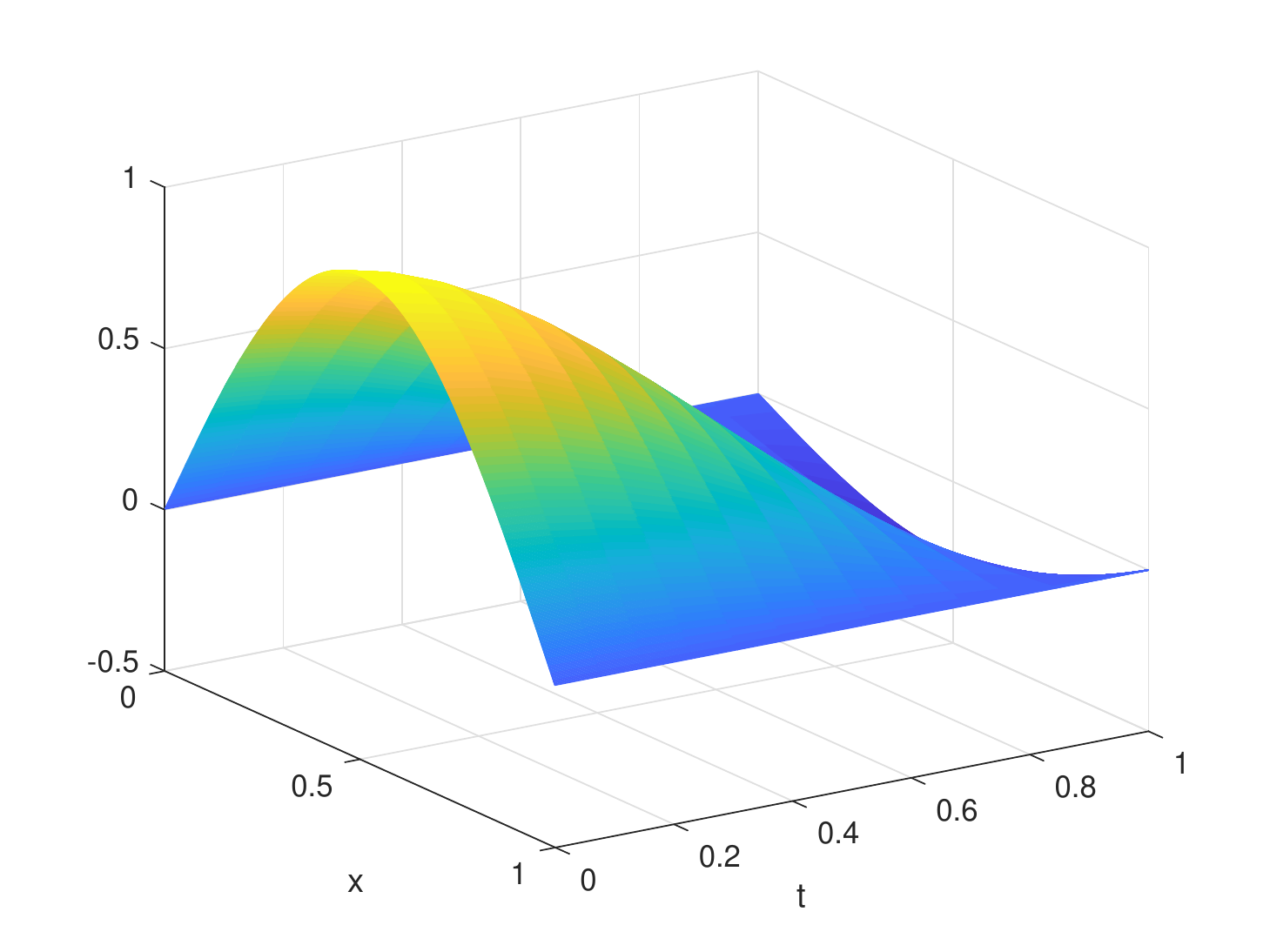}	
	\includegraphics[scale= 0.25]{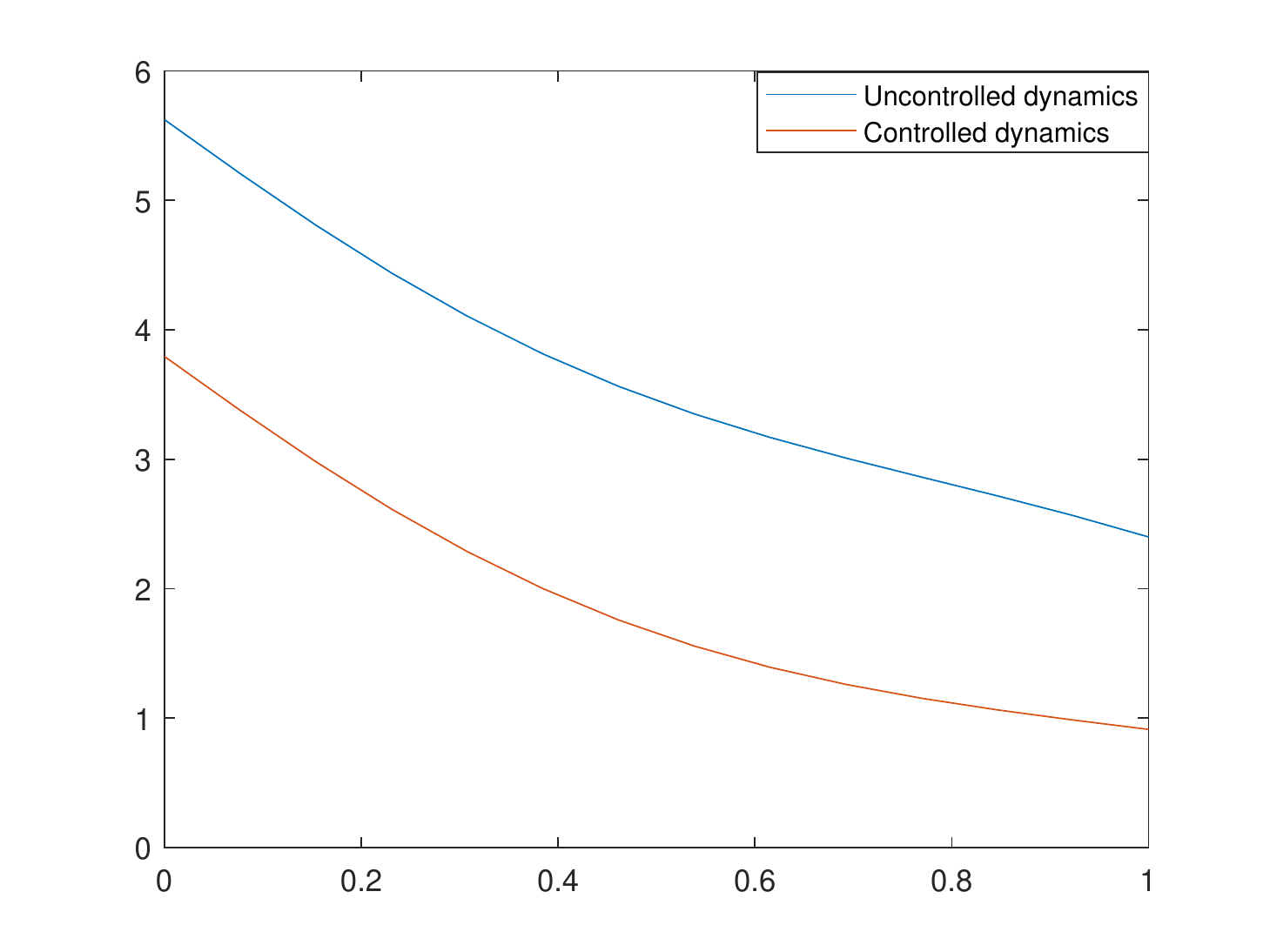}		
\caption{Test 5: Uncontrolled solution (left), optimal control solution (middle),  time comparison of the cost functional of the uncontrolled solution and controlled solution (right).}	
\label{fig5:wave}
\end{figure}

Again we apply an implicit Euler scheme to avoid narrow CFL conditions. We want to minimize the following cost functional 
 $$
J_{y_0,t}(u)= \int_0^T \left( \varphi(\Vert y(s) \Vert^2_2) + \gamma |u(s)|^2 \right)\, ds + \varphi(\Vert y(T)\Vert^2_2) \;,
$$
with $w_0(x)=\sin(\pi x),\, w_1(x)=0$, $\gamma=0.01,\, T=1,\, c=0.5,\, \Omega=(0,1)$ and $\omega=(0.4,0.6)$, $\Delta x=10^{-3}, \Delta t=0.05$. We note that the dimension of the semi-discrete problem is $d=2000$.
 \paragraph{Quadratic cost functional} 
 We first consider a standard tracking problem e.g. $\varphi(x) = x$ in the cost functional. In Figure \ref{fig5:wave} we show the uncontrolled solution in the left panel and the controlled solution in the middle. A comparison of the evaluations of the cost functional is given in the right panel. As expected the controlled solution is below the uncontrolled one  for each time instance. This shows the capability of the method for high dimensional problem even for hyperbolic equations.

\begin{figure}[htbp]	
\centering
	\includegraphics[scale=0.25]{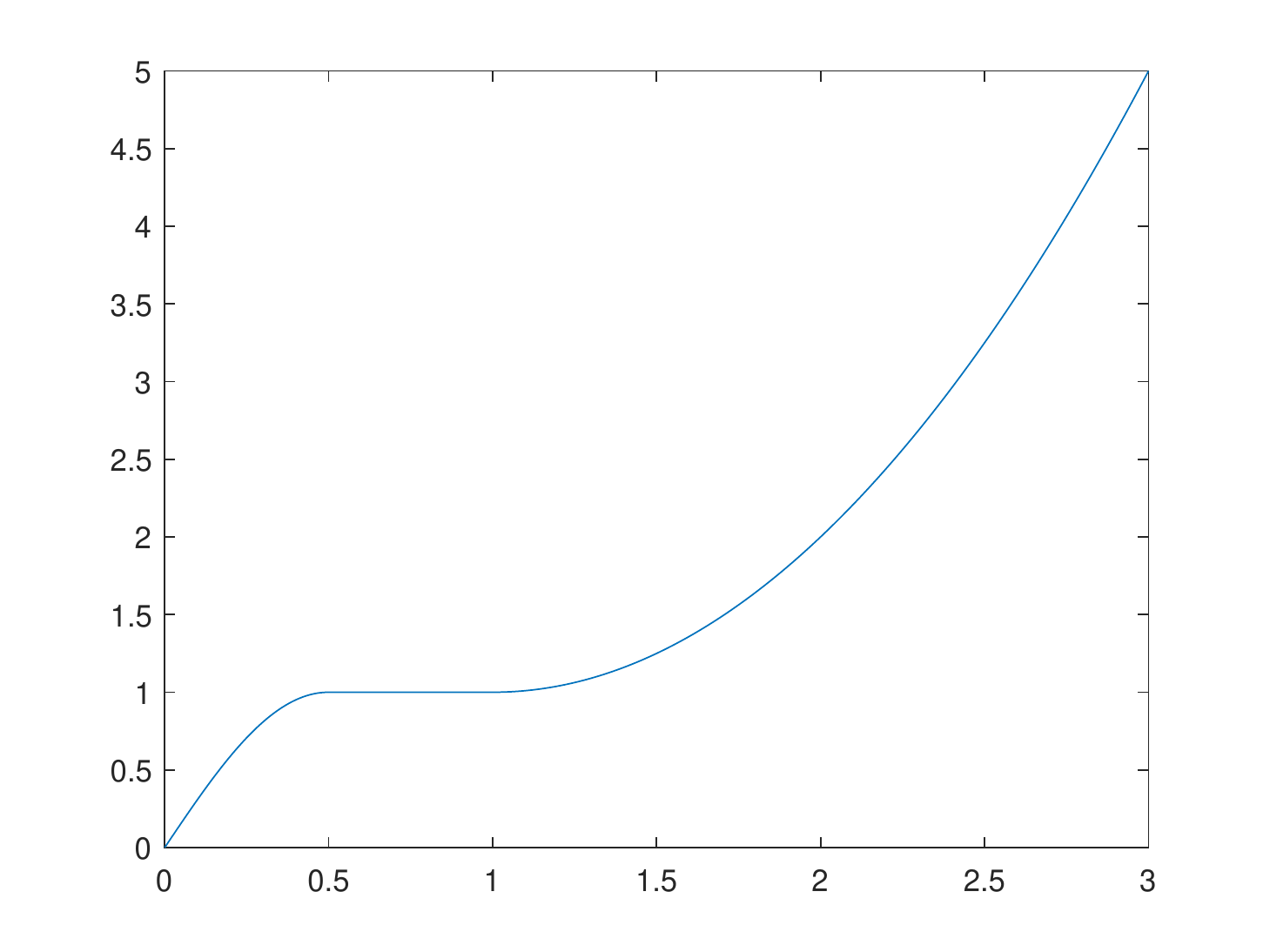}	
	\includegraphics[scale = 0.25]{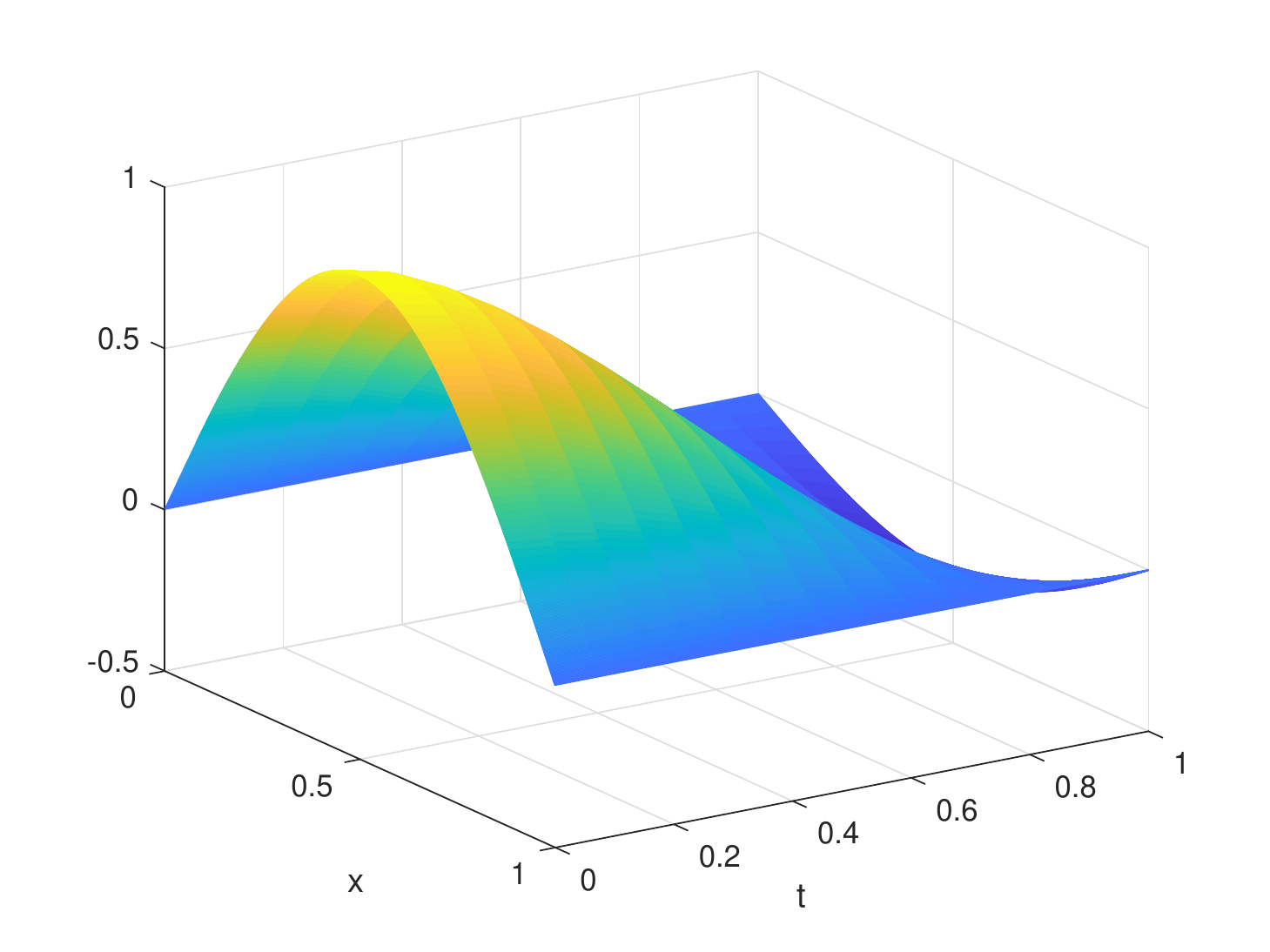}	
	\includegraphics[scale= 0.25]{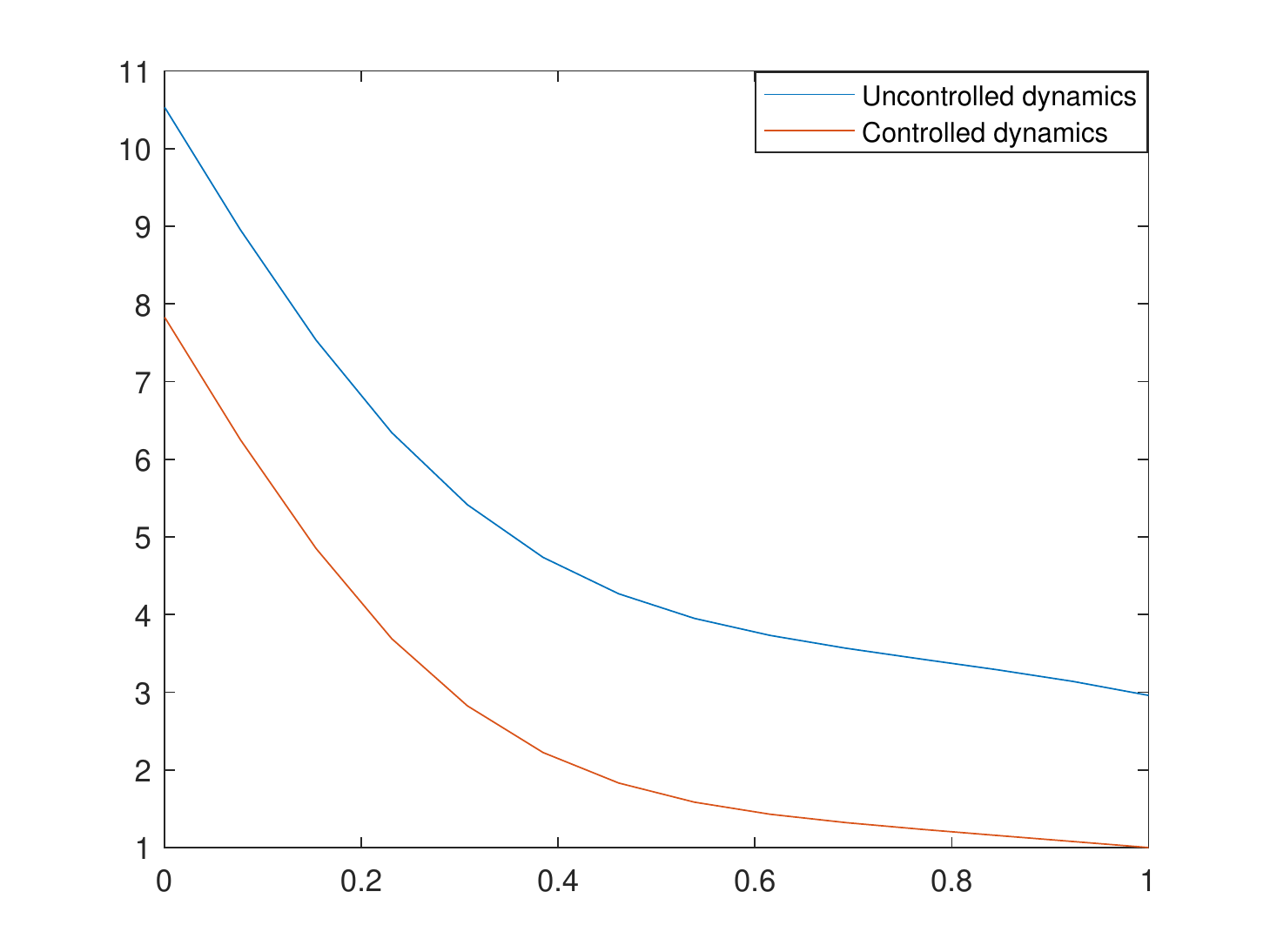}		
\caption{Test 5 (Non-quadratic cost functional): Graphics of $\varphi(x)$ (left), optimal control solution (middle), time comparison of the cost functional of the uncontrolled solution and controlled solution (right).}	
\label{fig6:wave}
\end{figure}

 \paragraph{Non-quadratic cost functional} \revfirst{Now, we consider a more complicated example which deals with a non-quadratic cost functional.} Let us consider for example the following cost functional where 
$$
\varphi(x)= \begin{cases} 
sin(\pi |x|) & |x| \le 0.5 \;,\\
1 & 0.5<|x|\le 1\;, \\
(|x|-1)^2 +1 & |x|>1   \;,
\end{cases}    
$$
as shown in the left panel of Figure \ref{fig6:wave}.
We consider the same parameters as in the previous case, which lead to the same uncontrolled solution as shown in the left panel of Figure \ref{fig5:wave}. 
In the middle of Figure \ref{fig6:wave} one can see the uncontrolled solution and in the right panel a comparison of the evaluation of the cost functional. Again, here we would like to stress the capability of the method to work with high dimensional problem and with non-smooth cost functionals.

\section{Conclusions and future works}\label{sec:con}

We have proposed a novel method to approximate time dependent HJB equations via DP scheme on a tree structure. The proposed algorithm creates the tree structure according to all the possible directions of the controlled dynamical system for a finite set of controls. This procedure has several advantages with respect to the DP algorithm based on the classical time and space discretization. The first advantage is that we do not have to build a space grid and
a local space interpolation. Furthermore, TSA does not require an a-priori choice of a numerical domain $\Omega$ to set the numerical scheme and, consequently, there is no need to impose boundary conditions. The construction of the tree is made step-by-step, via the pruning rule. Thus, the complexity of the problem is drastically reduced cutting all the branches laying in a small neighbourhood. After pruning the tree the efficiency of TSA is greatly improved in terms of CPU time.
This approach allows to apply the DP method to high-dimensional problems  as it has been shown in the numerical section for both ODEs and PDEs, in some test problems we solved an optimal control problem in dimension $2000.$ \\
We note that the method could be easily extended to second order approximation schemes using e.g. Heun method for the dynamics and it is also possible to reduce the CPU time via a parallel version of the method. Although the numerical results are promising, some issues are still open in the analysis of the method. The first is to derive error estimates in agreement with the order of convergence shown in Table \ref{test1:tab_nosel}. Furthermore, we would like to extend the method to the control of nonlinear PDEs coupling TSA with model order reduction methods as discussed in \cite{AFV17} and taking advantage of the theoretical results found there. These extensions could open the way to the application of DP techniques for real industrial problems.  

\bibliographystyle{siamplain}


%
%
%
%
%
%
%
%
%
%
%
%
%

\end{document}